\documentclass[1 [leqno,11pt]{amsart}
\usepackage{amsfonts}
\usepackage{amssymb}
\usepackage{amssymb, amsmath}
\def\inte#1{
\displaystyle\mathop{#1\kern0pt}^\circ }

\let\pa=\partial

\let\f=\frac

\def\pa{\partial}

\def\virgp{\raise 2pt\hbox{,}}
\def\cdotpv{\raise 2pt\hbox{;}}

\def\eqdefa{\buildrel\hbox{\footnotesize def}\over =}

\def\C{\mathop{\mathbb C\kern 0pt}\nolimits}
\def\DD{\mathop{\mathbb D\kern 0pt}\nolimits}
\def\EE{\mathop{{\mathbb E \kern 0pt}}\nolimits}
\def\K{\mathop{\mathbb K\kern 0pt}\nolimits}
\def\N{\mathop{\mathbb N\kern 0pt}\nolimits}
\def\Q{\mathop{\mathbb Q\kern 0pt}\nolimits}
\def\R{\mathop{\mathbb R\kern 0pt}\nolimits}
\def\SS{\mathop{\mathbb S\kern 0pt}\nolimits}
\def\ZZ{\mathop{\mathbb Z\kern 0pt}\nolimits}
\def\TT{\mathop{\mathbb T\kern 0pt}\nolimits}
\def\P{\mathop{\mathbb P\kern 0pt}\nolimits}

\def\na{\nabla}

\newcommand{\ef}{ \hfill $ \blacksquare $ \vskip 3mm}

\newcommand{\beq}{\begin{equation}}
\newcommand{\eeq}{\end{equation}}
\newcommand{\ben}{\begin{eqnarray}}
\newcommand{\een}{\end{eqnarray}}
\newcommand{\beno}{\begin{eqnarray*}}
\newcommand{\eeno}{\end{eqnarray*}}

\newtheorem{thm}{Theorem}[section]
\newtheorem{lem}{Lemma}[section]
\newtheorem{rmk}{Remark}[section]
\newtheorem{col}{Corollary}[section]
\newtheorem{prop}{Proposition}[section]
\renewcommand{\theequation}{\thesection.\arabic{equation}}

\begin{document}
\title[Smoothing effects for Boltzmann equation]
{Smoothing effect for   Boltzmann equation with full-range
interactions}
 \author[Y. Chen]{Yemin Chen}%
\address[Y. Chen]
 { Department of Mathematics, Beijing Institute
of Technology\\
100081 Beijing, P. R.  China.} \email{yeminchen@bit.edu.cn}
\author[L. He]{Lingbing He}
\address[L. He]{Department of Mathematical Sciences, Tsinghua University\\
Beijing 100084,  P. R.  China.} \email{lbhe@math.tsinghua.edu.cn}

\maketitle
\begin{abstract}
In this work, we are concerned with the regularities of the
solutions to   Boltzmann equation with the physical collision
kernels for the full range of intermolecular repulsive potentials,
$r^{-(p-1)}$ with $p>2$. We
  give the new and constructive upper and lower bounds  for the
collision operator in terms of standard fractional Sobolev norm. As
an application, we prove that the strong solutions obtained by
Desvillettes \& Mouhot \cite{dm} to homogeneous Boltzmann equation
and classical solutions obtained by Gressman-Strain \cite{gs1,gs2}
or Alexandre-Morimoto-Ukai-Xu-Yang \cite{amuxy3,amuxy5} for the
inhomogeneous  Boltzmann equation become immediately smooth with
respect to all variables. And as another application, we obtain  the
global entropy dissipation estimate  which is a little stronger than
the one of Alexandre-Desvillettes-Villani-Wennberg \cite{advw}.

\end{abstract}
\tableofcontents

\noindent {\sl Keywords:} Hypoellipticity, homogeneous and
inhomogeneous Boltzmann
equation, full-range interactions, Littlewood-Paley analysis \\

\vskip 0.2cm

\noindent {\sl AMS Subject Classification (2010):} 35Q20, 35R11, 75P05.  \\

\renewcommand{\theequation}{\thesection.\arabic{equation}}
\setcounter{equation}{0}

\section{Introduction}
In the present work, we continue the study on the  smoothness of the
solutions to the Boltzmann equation with the collision kernels for
the inverse intermolecular potentials $r^{-(p-1)}$ with $p>2$. It is
well known that the Boltzmann equation is a fundamental equation in
statistical physics. The readers can refer to \cite{ce1, ce2, LL,
vi2} and the references therein for the physical background of the
equation and also for the mathematical theories for this equation.
Mathematically, the Boltmzann equation reads:
\begin{eqnarray}
\partial _t f+v\cdot\nabla  _x f=Q(f,f),
 \label{bol}
\end{eqnarray}
where $f(t,x,v)\geq 0$ is the (spatially periodic) distribution
function in the phase space of collision particles which at time
$t\geq 0$ and point $x\in \TT^3=[-\pi ,\pi ]^3$ move with velocity
$v\in\R^3$. The Boltzmann collision operator $Q$ is a bilinear
operator which acts only on the velocity variables $v$, that is,
\beno Q(g,f)(v)\eqdefa
\int_{\R^3}\int_{\SS^{2}}B(v-v_*,\sigma)(g'_*f'-g_*f)d\sigma dv_*.
\eeno Here we use the standard shorthand $f=f(v)$, $g_*=g(v_*)$,
$f'=f(v')$, $g'_*=g(v'_*)$ where $v'$, $v_*'$ are given by
\begin{eqnarray}\label{e3}
v'=\frac{v+v_{*}}{2}+\frac{|v-v_{*}|}{2}\sigma\ ,\ \ \
v'_{*}=\frac{v+v_{*}}{2}-\frac{|v-v_{*}|}{2}\sigma\
,\qquad\sigma\in\SS^{n-1}.
\end{eqnarray}
We stress that the representation follows the parametrization of the
set of solutions of the physical law of elastic collision: \beno
v+v_*&=&v'+v_*',\\ |v|^2+|v_*|^2&=&|v'|^2+|v'_*|^2. \eeno

The nonnegative function $B(v-v_*,\sigma)$  in the collision
operator is called the Boltzmann collision kernel. It is always
  assumed to depend only on $|v-v_{*}|$ and $\langle\frac{v-v_{*}}{|v-v_{*}|},\sigma
  \rangle$. Usually, we introduce the angle variable $\theta$ through
  $\cos\theta=\langle\frac{v-v_{*}}{|v-v_{*}|},\sigma
  \rangle$.  Without loss of generality, we may
assume that $B(v-v_{*},\sigma)$ is supported in the set
 $0\leq\theta\leq\frac{\pi}{2}$ , i.e,
$\langle\frac{v-v_{*}}{|v-v_{*}|},\sigma
  \rangle\ge0
$, for otherwise $B$ can be replaced by its symmetrized form:\beno
\bar{B}(v-v_{*},\sigma)=[B(v-v_{*},\sigma)+B(v-v_{*},-\sigma)]\mathbf{1}_{\langle\frac{v-v_{*}}{|v-v_{*}|},\sigma
  \rangle\ge0}.\eeno
Above, $\mathbf{1}_A$ is the characteristic function of the set $A$.
The typical example we have in mind is the case that the interaction
potential obeys the inverse repulsive potential which takes the
form of \beno \phi(r)=r^{-(p-1)}.\eeno

 According to these
potentials, in this paper, we consider the  collision kernel
satisfying the following assumptions:

{\bf Assumption A}
\begin{itemize}
\item  The cross-section $B(v-v_{*},\sigma)$ takes a product form as
\ben \label{e6} B(v-v_{*},\sigma)=\Phi(|v-v_*|)b(\cos\theta),\een
 where both $\Phi$ and $b$ are nonnegative functions.
\item The angular function $b(t)$ is not locally integrable and it
satisfies \begin{eqnarray}\label{e7} K\theta^{-1-2s} \le\sin\theta
b(\cos\theta)\le K^{-1}
 \theta^{-1-2s},  \quad
\mbox{with}\   0<s<1,\ K>0.
\end{eqnarray}
\item The kinetic factor $\Phi$ takes the form \ben\label{e15}
\Phi(|v-v_*|)=|v-v_{*}|^{\gamma},\een
 where the parameter $\gamma$   verifies that  $\gamma+2s>-1$.
\end{itemize}

We remark that for inverse repulsive potential, there holds that
$\gamma=\frac{p-5}{p-1}$ and $s=\frac{1}{p-1}$. It is easy to check
that $\gamma+4s=1$ which gives the sense of the assumption
$\gamma+2s>-1$. Generally, the case $\gamma>0$, $\gamma=0$, and
$\gamma<0$ correspond to so-called hard, maxwellian, and soft
potentials.

There are lots of literatures on the well-posedness problem of the
Boltzmann equation, and we will start off by mentioning a brief few.
As for the  case of Grad's angular cut-off, in 1989, DiPerna and
Lions \cite{dl} proved the celebrated result: the global existence
of renormalized solution to  the inhomogeneous Boltzmann equation
with arbitrary initial data. Thanks to this breakthrough, based on
the new definition of weak solution, the hydrodynamic limit from
Boltzmann equation to the equations of fluid mechanics can be
considered afterwards, see \cite{fs, lm1, lm2}. Another direction to
obtain the global solution is due to the work by Guo
\cite{guo1,guo2} and Liu-Yang-Yu \cite{lyy} who introduce the
nonlinear energy method to construct the classical solutions near
the equilibrium. We point out that their approach relies heavily on
the analysis of the linearized Boltzmann operator.

As for the non cut-off theory which means physically relevant
effects of the angular singularities are considered, it has been
made big progress in these years, see \cite{al1}. In 1995,
Desvillettes in \cite{ld1} first showed that the solution of the
spatially homogeneous non cut-off Kac equation becomes very regular
with respect to the velocity variable as soon as the time is
strictly positive. This testified the conjecture that when the cross
section is concentrating on the grazing collisions, the nonlinear
collision operator should behave like a fractional Laplacian in the
velocity variable $v$, see \cite{DTTSP1, goudon, vi3}. Later on,
  Alexandre in \cite{al2} formally showed that the smoothness estimates could indeed be
 deduced from the entropy dissipation
  $D(g,f)$ defined as \beno D(g,
f)=-\int_{\R^6}\int_{\SS^2}B(|v-v_*,\sigma)(g_*'f'-g_*f)\log fd
\sigma dvdv_*.\eeno Lions in \cite{lions} proved a functional
inequality of the form \beno  \|\sqrt{f}\|_{\dot{H}^\alpha(|v|\le
R)}\le C_R\|f\|_{L^1_v}^{\theta}(\|f\|_{L^1_v}+D(f,
f))^{1-\theta}\eeno for $\alpha<s$. Shortly after, the optimal
Sobolev exponent $s$ was achieved   by   Villani \cite{vi1} and
Alexandre \cite{al3} but at the price that solution is required to
be locally bounded below. In 2000, the work of
Alexandre-Desvillettes-Villani-Wennberg \cite{advw} showed that
usual estimate on the entropy dissipation automatically entails
regularization effects:\ben\label{co1} \|\sqrt{f}\|^2_{H^s(|v|\le
R)}\le C_{g,R}[D(g, f)+\|g\|_{L^1_2(\R^3_v)}\|f\|_{L^1_2(\R^3_v)}].
\een   which indicates that for a given $g\in L^1_v$, the Boltzmann
operator behaves as: \ben\label{be} Q(g,f)=-C_g(-\triangle)^s f +
\mbox{lower\,order\,terms}.\een Moreover, two basic formula such as
the cancellation lemma and sub-elliptic coercivity estimates are
also given there. Thanks to this breakthrough, Alexandre-Villani in
\cite{av,av2} first
 generalized the renormalized solution with defect measure for
Boltzmann equation with long-range interaction and then  gave the
rigorously justification to the Landau approximation. Another
application of the basic tools is to demonstrate the smoothing
effect of the classical solutions to the spatially homogeneous
Boltzmann equation with regularized potentials \cite{dw, muxy, alsa,
hmuy}. We mention that smoothing behavior is radically different
from that of the Boltzmann equation with angular cutoff (see for
example \cite{parme, dly} and the references therein for precise
statements). In this last case, propagation of regularity as well as
singularities (in the variable $v$) occurs, thanks to the properties
of the positive part of Boltzmann operator (Cf. \cite{w}, \cite{bd1}
and \cite{lu}). In 2009, Desvillettes and Mouhot in \cite{dm} proved
some a priori estimates for the stability and uniqueness for spatial
homogeneous Boltzmann equation with long-range interaction, and they
also showed the existence for moderate angular singularities.

 Recently, Alexandre-Morimoto-Ukai-Xu-Yang
\cite{amuxy1} introduced the pseudo-differential operator and
harmonic analysis to build so-called uncertainty principle to study
the hypoellipticity of the kinetic equation. And as the application,
they showed the regularizing effects for the linearized Boltzmann
equation with non cut-off and linearized Landau equation. Later on,
in \cite{amuxy2}, for the modified kinetic factor, that is, \beno
\Phi(|v|)=(1+|v|^2)^{\f12},\eeno based on the pseudo-differential
calculus and generalized Bobylev formula(see \cite{bob}), they
developed the methods to sharpen the upper bound estimate for the
Boltzmann collision operator (see also \cite{al, alhe}) which helped
them not only to establish the local existence theory for the non
cut-off inhomogeneous Boltzmann equation with arbitrary initial data
and but also to prove the instantaneous smoothness of the solutions.
More recently, Gressman-Strain \cite{gs1,gs2,gs3} and
Alexandre-Morimoto-Ukai-Xu-Yang \cite{amuxy3,amuxy4,amuxy5}
independently established the global existence of the classical
solutions to the Boltzmann equation with long-range of inverse power
intermolecular potentials, $r^{-p+1} $ with $p>2$ when the initial
data  are  close to the equilibrium. Both of the methods rely on the
estimate for the linearized collision operator.  Let us give some
comments on the work of the coercive estimate. In \cite{gs1}, the
authors showed that at the linearized level, the collision operator
can be regarded as a fractional Laplacian on a manifold and this
manifold depends in an essential way on the collision geometry. More
recently, in \cite{gs3},  they provide sharp constructive upper and
lower bound estimates for the Boltzmann collision operator. It is
shown that under  the assumption of  high regularity and
sufficiently rapid growth of the weight at infinity   on the
function $g$, there holds that \ben\label{gs}
C_g^1\|f\|^2_{\dot{N}^{s,\gamma}}\lesssim -\langle
Q(g,f),f\rangle_v+\|f\|_{L^2_\gamma}^2\lesssim C_g^2
\|f\|^2_{N^{s,\gamma}},\een  with
 \beno   \|f\|^2_{N^{s,\gamma}}\eqdefa \|f\|_{L^2_{\gamma+2s}}^2+\|f\|^2_{\dot{N}^{s,\gamma}}\eeno and \beno
  \|f\|^2_{\dot{N}^{s,\gamma}}\eqdefa  \int_{\R^6} (\langle v\rangle \langle
v'\rangle)^{\f{\gamma+2s+1}{2}}\frac{( f' - f
)^2}{d(v',v)^{3+2s}}\mathbf{1}_{d(v',v)\le1}dvdv'.\eeno Here the
non-isotropic metric $d(v,v')$ is defined on the "lifted"
paraboloid: \beno d(v,v')\eqdefa
\sqrt{|v-v'|^2+\f14(|v|^2-|v'|^2)^2}.\eeno Moreover, under the same
assumption, they prove the global entropy production estimates which
is \ben\label{gsc}  D(g,f)\ge \|\sqrt{f}\|^2_{\dot{N}^{s,\gamma}}-
C_g^3\|f\|_{L^1_\gamma}. \een We remark that the norm of $
\dot{N}^{s,\gamma} $ is a semi-norm. While in the work of
Alexandre-Morimoto-Ukai-Xu-Yang \cite{amuxy3, amuxy4}, they gave
another way to understand the coercivity of the collision operator.
Precisely, in contrast to \cite{advw},  they regarded the quantity
\beno \int_{\R^6}\int_{\SS^2}|v-v_*|^\gamma b(\cos
\theta)e^{-\frac{|v_*|^2}{2}}(f'-f)^2d\sigma dvdv_*\eeno as the new
norm instead of standard fractional Sobolev norm to bound the
linearized Boltzmann operator. This is key point to construct the
global classical solutions of the Boltzmann equation when the
initial data are near equilibrium.

In the present work, we are going to investigate the regularities of
the solutions to both homogeneous and inhomogeneous Boltzmann
equation with the physical collision kernels for the full range of
intermolecular repulsive potentials. It can be viewed as a
continuation of the recent work \cite{cdh} where they demonstrated
the $C^\infty $ regularizing effect for the full Landau equation. As
we known, the main difficulty  to prove the smoothing effect for the
 nonlinear Boltzmann equation comes from the upper and lower bound
for the collision operator. The main reason lies in the fact that
the Boltzmann operator only involves singular integral behaving like
a fractional differential operator but no explicit derivative or
pseudo-differential operator occurs.

To overcome the difficulty, motivated by the collision geometry  and
the standard Littlewood-Paley decomposition, we carry out the new
strategy to bound the dual form $\langle Q(g, h), f\rangle_v$.
Roughly speaking, in contrast to the previous work \cite{amuxy2}, by
denoting $\mathcal {G}=\langle v\rangle^{N_1}g, \mathcal {H}=\langle
v\rangle^{N_2}h$ and $ \mathcal {F}=\langle v\rangle^{N_3}f$,  we
first transform  $\langle Q(g, h), f\rangle_v$ to the new functional
$\langle \mathcal
{Q}(\mathcal{G},\mathcal{H}),\mathcal{F}\rangle_v$. The most
convenience of the transformation is that the new factors $ \langle
v_*\rangle^{-N_1}, \langle v\rangle^{-N_2}$ and $\langle
v'\rangle^{-N_3}$  which are inside the new functional will absorb
the weight  coming from the cross-section. Thanks to this design,
now we can apply the Littlewood-Paley decomposition to the functions
$\mathcal{H}$ and $\mathcal{F}$ to make full use of the cancellation
between the different frequency part of them. Combined with the
Bernstein's inequality and the proper cut off for the angular, the
upper bound estimate in terms of standard fractional Sobolev norm
for the functional is finally obtained which also implies   the
upper bound for the collision operator  by duality. One may check
the details in section 2.

Another contribution of the paper lies in the new estimation for the
coercivity of the Boltzmann collision operator. We show that for the
non Maxwellian potentials, the global sub-elliptic estimate  with
some weight can be obtained. Roughly, if $\gamma+2s>0$ and
$\gamma\le 2$, the estimate for regularizing effect \eqref{co1} can
be improved as: \ben\label{co2}
 \|\sqrt{f}\langle v\rangle
^{\frac{\gamma}{2}}\|_{H^s(\R^3)}^2&\le& C_{g}[D(g, f)+(\|g
\|_{L^1_2(\R^3_v)} +1)^{C_s}\|f \|_{L^1_2(\R^3_v)} ]. \een While for
the case of $\gamma+2s\le 0$, the similar estimate as \eqref{co2}
still can be  obtained but at the cost that we have to impose the
condition of high integrability (for instance, $L^{\f32+\delta}$
with $\delta>0$) on the function $g$ which is also observed for the
estimate to collision operator (one may check the corresponding
theorems for details). We remark that the critical value
$\gamma=-2s$ corresponds to the threshold below which there is no
spectral gap for the linearized Boltzmann operator.
   We also point out that comparing to the estimate
\eqref{gsc}, we only use the conserved quantities of Boltzmann
equation to capture the smoothing effect in \eqref{co2} in the case
of $\gamma+2s>0$. One may check the corresponding section for
details.

With in hand the upper and lower bound for the collision operator,
 now we are in a position to state our main results. The first one is
concerned with the spatial homogeneous Boltzmann equation which
means the distribution function does not depend on the spatial
variables, i.e, \beno \pa_t f= Q(f,f).\eeno

\begin{thm}\label{smh} Let the collision kernel $B(|v-v_*|, \sigma)$ verify the assumption A, and   $f$ be  the unique solution of the
homogeneous Boltzmann equation satisfying the infinite moment
estimates, that is,  for any $l\in \R^+$, \ben\label{am1} \|f\langle
v\rangle^l\|_{L^\infty([0,\infty);L^1(\R^3_v))}<\infty, \quad
\mbox{if}\quad \gamma+2s>0; \een   or  \ben\label{am2} \|f\langle
v\rangle^l\|_{L^\infty([0,\infty);L^2(\R^3_v))}<\infty, \quad
\mbox{if}\quad \gamma+2s\le0. \een  Then for all $t_0>0$, the
solution $f$ lies in $L^\infty([t_0,\infty); \mathcal {S})$ .
\end{thm}

\begin{rmk} Noting the global existence result (for the case of $\gamma+2s>0$) and local existence
result (for the case of $\gamma+2s\le0$)  for moderate angular
singularity (which means $s<\f12$) by Desvillettes-Mouhot \cite{dm},
it shows that the result of the Theorem \ref{smh} is not empty.
Actually, following the proof of the Theorem \ref{smh}, we can show
that the regularity of the strong solution  constructed in Theorem
1.3 by Desvillettes-Mouhot \cite{dm} can be propagated which implies
that the strong solution  is exactly the classical solution  when we
impose the regularity on the initial datum.\end{rmk}

\begin{rmk}
To our knowledge, it is the first time to prove the smoothing effect
of the homogeneous Boltzmann equation for the "true" hard potentials
and "true" moderately soft potentials.
 We also mention that the assumption \eqref{am2} for the case
of $\gamma+2s$ can  be weakened  by  \ben\label{am4} \|f\langle
v\rangle^l\|_{L^\infty([0,\infty);L^\f32(\R^3_v))}<\infty, \quad
\mbox{if}\quad \gamma+2s\le0. \een The main reason lies in the upper
and lower bounds for the Boltzmann collision operator. We omit the
details here and one may check the corresponding parts in Section 4.
\end{rmk}

 Let us give some comments on the
difference between our result with the previous work \cite{dw, muxy,
 alsa, hmuy}. In their work,  they  actually deal with the case of modified
kinetic factor $\Phi$ which usually takes the form of
$(1+|v-v_*|^2)^{\f{\gamma}{2}}$. We stress out that this
mollification plays the key role in the proof to the smoothing
effect of the homogeneous Boltzmann equation. In fact, it will bring
them both upper and lower bounds for the collision operator. For the
upper bound,   the mollification makes it possible   to use
integration by parts with respect to $v_*$. Roughly speaking, by
Bobylev's formula, one has \beno \widehat{ -\triangle Q(g,h)}(\xi)=
\int_{\R^3}\int_{\SS^2}
\bigg[\hat{g}(\xi^-+\xi_*)\hat{f}(\xi^+-\xi_*)-\hat{g}(\xi_*)\hat{f}(\xi-\xi_*)\bigg]b(\frac{\xi}{|\xi|}\cdot
\sigma)\hat{\Phi}(\xi_*)|\xi|^2d\sigma d\xi,\eeno where
$\xi^+=\frac{\xi+|\xi|\sigma}{2}, \xi^-=\frac{\xi-|\xi|\sigma}{2}$.
Observing the fact \beno |\xi|^2=|\xi-\xi_*|^2+|\xi_*|^2+2(\xi-
\xi_*)\cdot \xi_* \eeno and \beno
|\xi|^2=|\xi^+-\xi_*|^2+|\xi_*|^2+2(\xi^+- \xi_*)\cdot
\xi_*+|\xi^-|^2, \eeno one may expect that the derivative required
for $g$ can be transferred to the kinetic factor $\Phi$ which leads
to the optimal upper bound for the collision operator, that is,
\beno \|Q(f,g)\|_{H^m_N}\lesssim
\|f\|_{L^1_{N^++(\gamma+2s)^+}}\|g\|_{H^{m+2s}_{(N+\gamma+2s)^+}}.\eeno
   While for the
lower bound, thanks to the inequality \beno \langle v\rangle ^\gamma
\le \langle v-v_*\rangle^\gamma\langle v_*\rangle^{|\gamma|}, \eeno
  the coercivity estimate    of
the collision operator for the case of hard potentials and soft
potentials can be concluded to the case of Maxwellian potential.
Thus the lower bound for the Boltzmann operator can be easily
obtained due to the work by Alexandre-Desvillettes-Villani-Wennberg
\cite{advw}. Since now  the collision kernel only verifies the
assumption A, one has to find another approach to give the estimates
to the upper and lower bound for the collision operator. And these
are exactly what we do in this paper.

For the inhomogeneous Boltzmann equation, to achieve our goal, we
still have to bypass the problem how to get the regularity with
respect to $x, v$. Thanks to the upper bound estimate for  the
collision operator, we show that $Q(g, f)$ belongs to the space
$L^2_{t,x}(H^{-s}_v)$. This means the hypo-elliptic estimate in
\cite{bou} for the kinetic equation can be applied. One may treat
the Boltzmann equation as Alexandre-Morimoto-Ukai-Xu-Yang done in
\cite{amuxy1} by employing the generalized uncertainty principle.
Here we opt for another approach which mainly comes from the work
\cite{cdh}: once the fractional derivatives(with respect to $x$) are
gained, to avoid estimating the commutator,  one may continue to
perform the energy estimates (and the estimates based on the
averaging lemmas) for weighted finite differences of derivatives of
$f$. By iteration,  we finally can obtain the full one derivative
with respect to $x$ and $v$. One has
\begin{thm}\label{smih} Let the collision kernel $B(|v-v_*|, \sigma)$ verifies the assumption A, and
 $f$ be  the unique classical solution of the
inhomogeneous Boltzmann equation satisfying that  for any $l\in
\R^+$, \ben\label{am3} \|f\langle
v\rangle^l\|_{L^\infty([0,\infty);H^5_{x,v})}<\infty, \een and there
exists a universal constant $C_l$ and $C_u$ such that for any
$x\in\TT^3$, \ben\label{am5}
0<C_l<\|f\|_{L^1(\R^3_v)}<C_u<\infty.\een Then the solution $f$ lies
in $ W^{\infty,\infty}([t_0,\infty);H^{\infty,l}_{x,v}) $ for all
$t_0>0$.
\end{thm}

Let $\mu=\frac{1}{(2\pi)^\f32}e^{-|v|^2}$ be the normalized
Maxwillian and  $F=F(t,x,v)$ be the standard perturbation with
respect to $\mu$ as
\begin{equation*}
f=\mu+\sqrt{\mu}F.
\end{equation*}
Then by Theorem 1 of \cite{gs2}, we get for any $l\geq 0$ and any
integer $N\geq 5$,
\begin{equation*}
\|f(t)\langle v\rangle^l\|_{H^N_{x,v}}\leq C_1+C_2\|F_0\langle
v\rangle^l\|_{H^N_{x,v}},
\end{equation*}
where $C_1, C_2$ depends   on $l$ and the constants appeared in
Theorem 1 of \cite{gs2}. Moreover, simple calculation gives that for
any $x\in\TT^3$, \beno
\|\mu\|_{L^1(\R^3)}-\|F\|_{L^\infty}\|\sqrt{\mu}\|_{L^1(\R^3)}\le
\|f\|_{L^1(\R^3)}\le
\|\mu\|_{L^1(\R^3)}+\|F\|_{L^\infty}\|\sqrt{\mu}\|_{L^1(\R^3)}.
\eeno  Choose $\|F\|_{L^\infty}$ small enough and then we can obtain
the estimate \eqref{am5} which implies that the Theorem \ref{smih}
can be applied to the solutions constructed by Gressman-Strain
\cite{gs1,gs2} or Alexandre-Morimoto-Ukai-Xu-Yang
\cite{amuxy3,amuxy5} for the inhomogeneous Boltzmann equation.

 The rest of the
paper will be organized as follows. First of all, in section 2, we
will use Littlewood-Paley analysis to study the upper bound estimate
for the collision operator. Moreover, the estimate for the
commutator between weight and Boltzmann operator is also given
there. In section 3, we will give the proof to the improved
coercivity estimate of the collision operator. Then in next two
sections, the regularizing effect for homogeneous and inhomogeneous
Boltzmann equation will be proven under the initial regularity
assumption on the solution. In the appendix, we shall give the proof
to some useful interpolation inequality.

\medskip

Let us complete this section by the function spaces and notations,
which we shall use throughout the paper. For notational simplicity,
we omit the integrating domains $\TT^3$ and $\R^3$, which correspond
to variables $x$ and variable $v$ respectively. For example, we
write $L^2_{x,v}$ instead of $L^2_x(\TT^3;L^2_v(\R^3))$. For
integer $N\geq 0$, we define the Sobolev space
\begin{equation*}H^N_{x,v}=\bigg\{f(x,v): \sum_{|\alpha |+|\beta |\leq N}\|\partial _x^{\alpha }
\partial _v^{\beta }f\|_{L^2_{x,v}}<+\infty\bigg\},
\end{equation*}
and for integer $N\geq 0$ and real number $l\geq 0$, we define the weighted Sobolev space
\begin{equation*}
H^{N,l}_{x,v}=\bigg\{f(x,v): \sum_{|\alpha |+|\beta |\leq
N}\|(\partial _x^{\alpha }\partial _v^{\beta }f) \,
 \langle v\rangle^l\|_{L^2_{x,v}}<+\infty\bigg\},
\end{equation*}
where the multi-index $\alpha =(\alpha _1,\alpha _2,\alpha _3)$,
$|\alpha |=\alpha _1+\alpha _2+\alpha _3$ and $\partial _x^{\alpha
}=\partial _{x_1}^{\alpha _1}\partial _{x_2}^{\alpha _2}\partial
_{x_3}^{\alpha _3}$ with $x=(x_1,x_2,x_3)$, and $ \langle
v\rangle=(1+|v|^2)^{\frac{1}{2}}$. The notations for $\beta $ are
the same. It is obvious that $H^{N,0}_{x,v}=H^N_{x,v}$. We also
define $H^{\infty}_{x,v}$ and $H^{\infty,l}_{x,v}$ by
\begin{equation*}
H^{\infty}_{x,v}=\bigcap_{N\geq 0}H^N_{x,v},\qquad
H^{\infty,l}_{x,v}=\bigcap _{N\geq 0}H^{N,l}_{x,v}.
\end{equation*}

We also introduce the standard notations \beno
\|f\|_{L^p_\alpha}=\big(\int_{\R^3} |f(v)|^p\langle v
\rangle^{\alpha p}dv\big)^{\f1{p}},\quad \|f\|_{L\log L}=\int_{\R^3}
f\log (1+f)dv,\eeno and \beno  \|f\|_{H^s_l}=\|f\langle v
\rangle^{l}\|_{H^s}.\eeno $a\lesssim b$, we mean that there is a
uniform constant $C,$ which may be different on different lines,
such that $a\leq Cb$. $a\sim b$ if both $a\lesssim b$ and $b\lesssim
a$.

\section{Upper bound on the collision operator}
In this section, we shall give the upper bound estimate for the
collision operator. Our main motivation comes from the singularity
of the cross-section and collision geometry which allow us to apply
Littlewood-Paley analysis to the boundedness of the collision
operator in terms of weighted fractional Sobolev spaces. It is  one
of the key steps to prove the smoothing effect of the non cut-off
Boltzmann equation. We remark that the variables $(t, x)$ are
considered as parameter for the all the estimates in this section.

\begin{thm}\label{ub} Let $0<s<1$ and   $N_2, N_3\in \R$. Suppose $ N_1=|N_2|+|N_3|$ and $\tilde{N}_1 =N_2+N_3$ with $\tilde{N}_1\ge
\gamma+2s $. Then for nonnegative and smooth functions $g, h$ and
$f$, there hold

\begin{enumerate}
\item if $\gamma+2s>0$,     \ben\label{u1} |\langle Q(g,h),f\rangle_v|\lesssim
\|g\|_{L^1_{N_1}(\R^3_v)}\|h\|_{H^s_{N_2}(\R^3_v)}\|f\|_{H^s_{N_3}(\R^3_v)};
 \een
 \item if $\gamma+2s\le 0$,     \ben\label{u2} |\langle Q(g,h),f\rangle_v|\lesssim
(\|g\|_{L^1_{N_1}(\R^3_v)}+\|g\|_{L^\frac{3}{2}_{N_1}(\R^3_v)})\|h\|_{H^s_{N_2}(\R^3_v)}\|f\|_{H^s_{N_3}(\R^3_v)}.
 \een
\end{enumerate}
\end{thm}

\medskip

 Let us give some comments on the main result of the theorem. First of all, \eqref{u1} and \eqref{u2} can be regarded as another proof to the
fact that the Boltzmann operator takes the form of \eqref{be}.
Secondly, we stress that the weight in $v$ comes  not only from the
kinetic factor $\Phi$ but also from the integration with respect to
the angular. Thirdly, in the case of $\gamma+2s<0$,  the additional
$L^\f32$ bound for the function $g$ results from the strong
singularity caused by the kinetic factor $\Phi$. Fourthly, by
duality, one may take $N_3=0$ to obtain the upper bound for the
collision operator which will be very useful in the next section.
Last we would point out that our proof relies only on the trick of
change of variables and cancellation lemma.

\medskip

 \noindent {\bf Proof of the Theorem \ref{ub}:}  By change of
variables, one may obtain \beno \langle
Q(g,h),f\rangle_v&=&\int_{\R^3}dv\int_{\R^3}dv_*\int_{\SS^{2}}B(v-v_*,\sigma)g_*h(f'-f)d\sigma\\
&=&\int_{\R^3}dv\int_{\R^3}dv_*\int_{\SS^{2}}B(v-v_*,\sigma)\langle
v_*\rangle^{-N_1}\langle v\rangle^{-N_2}(\langle v_*\rangle^{N_1}
g_*)(\langle v\rangle^{N_2}h)\\&&\quad\times (\langle
v'\rangle^{-N_3}(\langle v'\rangle^{N_3}f')-\langle
v\rangle^{-N_3}(\langle v\rangle^{N_3}f))d\sigma.\eeno

Set $\mathcal {G}=\langle v\rangle^{N_1}g, \mathcal {H}=\langle
v\rangle^{N_2}h$ and $ \mathcal {F}=\langle v\rangle^{N_3}f$. Then
we can rewrite the above equality as \beno  && \langle
Q(g,h),f\rangle_v\\&=&
\int_{\R^3}dv\int_{\R^3}dv_*\int_{\SS^{2}}B(v-v_*,\sigma)\langle
v_*\rangle^{-N_1}\langle v\rangle^{-N_2}\mathcal{G}_*\mathcal{H}
(\langle v'\rangle^{-N_3}\mathcal{F}'-\langle
v\rangle^{-N_3}\mathcal{F})d\sigma\\
&\eqdefa& \langle \mathcal
{Q}(\mathcal{G},\mathcal{H}),\mathcal{F}\rangle_v.\eeno In the
following analysis, we will turn our attention to the new defined
functional involving   the Boltzmann collision operator. Let us give
some comments on the new defined functional. It will bring us two
convenience: the first one is that the weight inside the integration
will absorb the polynomial of $|v-v_*|$ which probably comes from
the cross-section; the second one is that we can use the
Littlewood-Paley decomposition for the functions $\mathcal{H}$ and
$\mathcal{F}$ which is key for the upper bound estimates.

 Set $B\eqdefa\{\xi\in\R^3,\  |\xi|\leq
\f43\}$ and  $\mathcal {C}\eqdefa\{\xi\in\R^3,\ \f34\leq |\xi|\leq
\f83\}$. In view of Littlewood-Paley decomposition, one may
introduce two cut off functions $\phi\in C_c^\infty(B)$ and  $
\varphi\in C_c^\infty(\mathcal{C} )$  which satisfy \beno
\phi(\xi)+\sum_{j\ge0} \varphi(2^{-j}\xi)=1,\quad \xi\in \R^3.\eeno

  We denote $h\eqdefa \mathcal{F}^{-1}\varphi$ and $\tilde{h}\eqdefa
  \mathcal{F}^{-1}\phi$, then the
dyadic operators $\triangle_j$  can be defined as follows \beno
&&\triangle_jf
=2^{3j}\int_{\R^3}h(2^jy)f(x-y)dy,\quad\mbox{for}\quad
j\ge 0. \nonumber\\
&&\triangle_{-1}f=\int_{\R^3}\tilde{h}(y)f(x-y)dy.\eeno
 Then the new defined functional can be presented as \ben \langle \mathcal
{Q}(\mathcal{G},\mathcal{H}),\mathcal{F}\rangle_v=\sum_{j<k}\langle
\mathcal {Q}(\mathcal{G},\mathcal{H}_k),\mathcal{F}_j\rangle_v
+\sum_{j\le k}\langle \mathcal
{Q}(\mathcal{G},\mathcal{H}_j),\mathcal{F}_k\rangle_v,\een where
$\mathcal{H}_k=\triangle_k \mathcal{H}$ and
$\mathcal{F}_j=\triangle_j\, \mathcal{F}$. Now we will perform the
estimate for the Boltzmann collision operator. Since the proof is a
little bit longer, we shall divide it into two steps.

\medskip

 \noindent{\bf  Step 1: Frequency dominated by the function
$\mathcal{H} $.} We first treat with the case that the frequency of
function $\mathcal{H}_k$ prevails over  the one of the function
$\mathcal{F}_j$ which means  $j<k$. Introduce the smooth function
$\phi$ defined as before, set $\phi_j(w)=\phi(2^jw)$ and  one has
that \beno \langle \mathcal {Q}(\mathcal {G},\mathcal {H}_k),
\mathcal {F}_j\rangle_v =\sum_{i=1}^3 \Gamma_i,\eeno where \beno
\Gamma_1&\eqdefa&\int_{\R^3}dv\int_{\R^3}dv_*\int_{\SS^{2}}B(v-v_*,\sigma)\langle
v_*\rangle^{-N_1}\langle v\rangle^{-N_2}\mathcal{G}_*\mathcal{H}_k
\phi_j(|v-v'|) \\&& \mathbf{1}_{2^{-j}|v-v_*|^{-1}\le
c\f{\pi}{4}}(\langle v'\rangle^{-N_3}\mathcal{F}'_j-\langle
v\rangle^{-N_3}\mathcal{F}_j)d\sigma,\\
\Gamma_2&\eqdefa&\int_{\R^3}dv\int_{\R^3}dv_*\int_{\SS^{2}}B(v-v_*,\sigma)\langle
v_*\rangle^{-N_1}\langle v\rangle^{-N_2}\mathcal{G}_*\mathcal{H}_k
[1-\phi_j(|v-v'|)] \\&& \mathbf{1}_{2^{-j}|v-v_*|^{-1}\le
c\f{\pi}{4}}(\langle v'\rangle^{-N_3}\mathcal{F}'_j-\langle
v\rangle^{-N_3}\mathcal{F}_j)d\sigma,\eeno and \beno
\Gamma_3&\eqdefa&\int_{\R^3}dv\int_{\R^3}dv_*\int_{\SS^{2}}B(v-v_*,\sigma)\langle
v_*\rangle^{-N_1}\langle v\rangle^{-N_2}\mathcal{G}_*\mathcal{H}_k
  \\&& \mathbf{1}_{2^{-j}|v-v_*|^{-1}\ge
c\f{\pi}{4}}(\langle v'\rangle^{-N_3}\mathcal{F}_j'-\langle
v\rangle^{-N_3}\mathcal{F}_j)d\sigma.\eeno

We remark that the above decomposition comes from the collision
geometry and the singularity caused by the angular.

 To overcome the strong
singularity caused by the collision kernel, motivated by the
cancelation lemma, we shall use standard Taylor expansion.
Precisely, let \ben\label{F}\mathbb{F}_j=\langle
v\rangle^{-N_3}\mathcal{F}_j,\een then one has \ben\label{ta1}
\mathbb{F}_j(v')-\mathbb{F}_j(v)=(v'-v)\cdot\na_v\mathbb{F}_j(v)+\int_0^1
(v'-v)\otimes(v'-v):\na_v^2 \mathbb{F}(\gamma(\kappa))d\kappa, \een
where $\gamma(\kappa)=\kappa v'+(1-\kappa)v$.

We stress that we only give the proof to the estimate  in the case
of $s\ge\f12$ and one may follow the same procedure to prove the
case of $s<\f12$.

\begin{lem}\label{le1} If $\gamma+2s>0$, there holds \ben
|\,\Gamma_1|\lesssim 2^{2js} \|\mathcal{G}\|_{L^1(\R^3_v)}
\|\mathcal{H}_k\|_{L^2(\R^3_v)}\|\mathcal{F}_j\|_{L^2(\R^3_v)}, \een
and if $\gamma+2s\le 0$, there holds \ben |\,\Gamma_1|\lesssim
2^{2js}
(\|\mathcal{G}\|_{L^1(\R^3_v)}+\|\mathcal{G}\|_{L^\frac{3}{2}(\R^3_v)})
\|\mathcal{H}_k\|_{L^2(\R^3_v)}\|\mathcal{F}_j\|_{L^2(\R^3_v)}. \een
\end{lem}
\noindent{ \bf Proof:}  $\Gamma_1$ can be  split into two parts
$\Gamma_{1,1}$ and $\Gamma_{1,2}$ which separately contain  the term
in the righthand side of \eqref{ta1}. Notice that \beno
&&\int_{\SS^2}  b(\langle
\frac{v-v_*}{|v-v_*|},\sigma\rangle)(v-v')\phi_j(|v-v'|)d\sigma\\
&&\quad=\int_{\SS^2}  b(\langle
\frac{v-v_*}{|v-v_*|},\sigma\rangle)\frac{v-v'}{|v-v'|}\cdot\frac{v-v_*}{|v-v_*|}|v-v'|\phi_j(|v-v'|)\frac{v-v_*}{|v-v_*|}d\sigma\\
&&\quad=\int_{\SS^2}  b(\langle \frac{v-v_*}{|v-v_*|},\sigma\rangle)
(\f{1-\langle
\frac{v-v_*}{|v-v_*|},\sigma\rangle}{2})\phi_j(|v-v'|)d\sigma
(v-v_*)
 \eeno

Since  $\sin\frac{\theta}{2}=\frac{|v-v'|}{|v-v_*|}$, one may obtain
that \beno &&|\int_{\SS^2}  b(\langle
\frac{v-v_*}{|v-v_*|},\sigma\rangle)(v-v')\phi_j(|v-v'|)d\sigma|\\
&&\lesssim \int_{\sqrt{\f{1-\langle
\frac{v-v_*}{|v-v_*|},\sigma\rangle}{2}}\lesssim 2^{-j}|v-v_*|^{-1}}
b(\langle \frac{v-v_*}{|v-v_*|},\,\sigma\rangle) \f{1-\langle
\frac{v-v_*}{|v-v_*|},\sigma\rangle}{2}d\sigma
|v-v_*|\\
&&\lesssim \int_0^{2^{-j}|v-v_*|^{-1}} \theta^{1-2s}d\theta |v-v_*| \\
&&\lesssim (2^{-j}|v-v_*|^{-1})^{2-2s}|v-v_*|.
 \eeno
 Observing that
$|v-v_*|\sim|v_*-v'|\sim|v_*-\gamma(\kappa)|$,  we may deduce that
 there holds \ben\label{wei} \langle v_*\rangle^{-N_1}\langle
v\rangle^{-N_2}\langle \gamma(\kappa)\rangle^{-N_3}\lesssim \langle
v_*-\gamma(\kappa)\rangle^{-\tilde{N}_1}.\een  The reader may check
it directly by the definition of $N_1$ and $\tilde{N}_1$.

 Then \beno
|\Gamma_{1,1}|&=&\bigg|\int_{\R^3}dv\int_{\R^3}dv_*\int_{\SS^{2}}B(v-v_*,\sigma)\langle
v_*\rangle^{-N_1}\langle v\rangle^{-N_2}\mathcal{G}_*\mathcal{H}_k
\phi_j(|v-v'|) \\&& \mathbf{1}_{2^{-j}|v-v_*|^{-1}\le
c\f{\pi}{4}} (v'-v)\cdot\na_v\mathbb{F}_j(v)d\sigma\bigg|\\
&\lesssim& \int_{\R^3}dv\int_{\R^3}dv_*
(2^{-j}|v-v_*|^{-1})^{1-2s}|v-v_*|^{\gamma+1}\langle v-
v_*\rangle^{-\tilde{N}_1}\mathcal{G}_*|\mathcal{H}_k|(|\mathcal{F}_j|+|\na_v\mathcal{F}_j|)
\eeno where we use the fact that $2^{-j}|v-v_*|^{-1}\le c\f{\pi}{4}$
and \ben\label{F_1}|\na_v\mathbb{F}_j(v)|\lesssim \langle
v\rangle^{-N_3}(|\mathcal{F}_j|+|\na_v\mathcal{F}_j|).\een For the
case of $\gamma+2s>0$, one may take $\tilde{N}_1\ge \gamma+2s$ and
get  \beno \Gamma_{1,1}&\lesssim&
2^{2js}2^{-j}\|\mathcal{G}\|_{L^1(\R^3_v)}\|\mathcal{H}_k\|_{L^2(\R^3_v)}(\|\mathcal{F}_j\|_{L^2(\R^3_v)}+\|\na_v\mathcal{F}_j\|_{L^2(\R^3_v)})\\
&\lesssim&2^{2js}
\|\mathcal{G}\|_{L^1(\R^3_v)}\|\mathcal{H}_k\|_{L^2(\R^3_v)}\|\mathcal{F}_j\|_{L^2(\R^3_v)}.\eeno
Here we use the Berstein's inequality.

For the case of $-1<\gamma+2s\le 0$, one may take
$\tilde{N}_1\ge-|\gamma+2s|$ and obtain the similar estimate
  \beno  \Gamma_{1,1}&\lesssim&\int_{\R^3}dv\int_{\R^3}dv_*
 2^{-j}2^{2js} (1+|v-v_*|^{\gamma+2s}\mathbf{1}_{|v-v_*|\le1}) \mathcal{G}_*|\mathcal{H}_k|(|\mathcal{F}_j|+|\na_v\mathcal{F}_j|)\\
&\lesssim&2^{2js}(\|\mathcal{G}\|_{L^1(\R^3_v)}+\|\mathcal{G}\|_{L^\frac{3}{2}(\R^3_v)})\|\mathcal{H}_k\|_{L^2(\R^3_v)}\|\mathcal{F}_j\|_{L^2(\R^3_v)}.\eeno

Next, we focus on the estimate for $\Gamma_{1,2}$. Notice \beno
\Gamma_{1,2}&\le& \int_0^1d\kappa
\int_{\R^3}dv\int_{\R^3}dv_*\int_{\SS^{2}}B(v-v_*,\sigma)\langle
v_*\rangle^{-N_1}\langle v\rangle^{-N_2}\mathcal{G}_*|\mathcal{H}_k|
\phi_j(|v-v'|) \\&& \mathbf{1}_{2^{-j}|v-v_*|^{-1}\le c\f{\pi}{4}}
|v'-v|^2|\na^2_v\mathbb{F}_j(\gamma(\kappa))|d\sigma.\eeno  Note
that \ben \label{F_2}|\na^2_v\mathbb{F}_j (v)|\lesssim \langle
v\rangle^{-N_3}(|\mathcal{F}_j|+|\na_v\mathcal{F}_j|+|\na^2_v
\mathcal{F}_j|).\een Then by Cauchy-Schwartz inequality, it suffices
to bound the quantities \beno I&\eqdefa& \int_0^1d\kappa
\int_{\R^3}dv\int_{\R^3}dv_*\int_{\SS^{2}}B(v-v_*,\sigma)\langle
v_*\rangle^{-N_1}\langle v\rangle^{-N_2}\mathcal{G}_* \phi_j(|v-v'|)
\\&& \mathbf{1}_{2^{-j}|v-v_*|^{-1}\le c\f{\pi}{4}}
|v'-v|^2\langle \gamma(\kappa)\rangle^{-N_3}(\sum_{i=0}^2
|\na_v^i\mathcal{F}_j(\gamma(\kappa))| )^2d\sigma \eeno and \beno
II&\eqdefa& \int_0^1d\kappa
\int_{\R^3}dv\int_{\R^3}dv_*\int_{\SS^{2}}B(v-v_*,\sigma)\langle
v_*\rangle^{-N_1}\langle v\rangle^{-N_2}\mathcal{G}_* \phi_j(|v-v'|)
\\&& \mathbf{1}_{2^{-j}|v-v_*|^{-1}\le c\f{\pi}{4}}
|v'-v|^2\langle \gamma(\kappa)\rangle^{-N_3}| \mathcal{H}_k
|^2d\sigma \eeno

Next we follow the well-known change of variables
$u=\gamma(\kappa)=\kappa v'+(1-\kappa)v$, which changes $v$ to $u$.
Thanks to the fact that the Jacobian is \ben\label{ja}
\bigg|\frac{du_m}{dv_n}\bigg|=(1-\frac{\kappa}{2})^2\{(1-\f{\kappa}{2})+\f{\kappa}{2}\langle\frac{v-v_*}{|v-v_*|},
\sigma\rangle\},\een  one has  \beno I&\lesssim&  \int_0^1d\kappa
\int_{\R^3}du\int_{\R^3}dv_*\int_{\SS^{2}}B(v-v_*,\sigma)\langle
v_*-u\rangle^{-\tilde{N}_1}\mathcal{G}_*(|\mathcal{F}_j(u)|+|\na_v\mathcal{F}_j(u)|\\&&+|\na^2_v
\mathcal{F}_j(u)|)^2\mathbf{1}_{2^{-j}|u-v_*|^{-1}\le
c\f{\pi}{4}}\phi_j(|v-v'|) |v'-v|^2d\sigma,\eeno where we use
\eqref{wei}. Due to the fact that \ben\label{ks}
|v'-v|^2&=&|v-v_*|^2\f{1-\langle
\frac{v-v_*}{|v-v_{*}|},\sigma\rangle}{2}\\&\lesssim&
|u-v_*|^2\f{1-\langle
\frac{v-v_*}{|v-v_{*}|},\sigma\rangle}{2},\nonumber\een one may
arrive at \beno I&\le&
\int_{\R^3}du\int_{\R^3}dv_*|v_*-u|^{\gamma+2}\langle
v_*-u\rangle^{-\tilde{N}_1}\mathcal{G}_*(|\mathcal{F}_j(u)|+|\na_v\mathcal{F}_j(u)|\\&&+|\na^2_v
\mathcal{F}_j(u)|)^2\mathbf{1}_{2^{-j}|u-v_*|^{-1}\le
c\f{\pi}{4}}\int_{\SS^2} b(\langle
\frac{v-v_*}{|v-v_*|},\sigma\rangle)\f{1-\langle
\frac{v-v_*}{|v-v_{*}|},\sigma\rangle}{2} \phi_j(|v-v'|)d\sigma.
\eeno Noting that \beno &&\int_{\SS^2}b(\langle
\frac{v-v_*}{|v-v_*|},\sigma\rangle)\f{1-\langle
\frac{v-v_*}{|v-v_{*}|},\sigma\rangle}{2}
\phi_j(|v-v'|)d\sigma\\&&\lesssim \int_{\sqrt{\f{1-\langle
\frac{v-v_*}{|v-v_*|},\sigma\rangle}{2}}\lesssim 2^{-j}|u-v_*|^{-1}}
b(\langle \frac{v-v_*}{|v-v_*|},\,\sigma\rangle) \f{1-\langle
\frac{v-v_*}{|v-v_*|},\sigma\rangle}{2}d\sigma\\&&\quad \lesssim
\int_0^{2^{-j}|u-v_*|^{-1}} \theta^{1-2s}d\theta\lesssim
(2^{-j}|u-v_*|^{-1})^{2-2s},\eeno we get \beno I&\le& 2^{-2j}2^{2js}
\int_{\R^3}du\int_{\R^3}dv_*|v_*-u|^{\gamma+2s}\langle
v_*-u\rangle^{-\tilde{N}_1}\mathcal{G}_*(|\mathcal{F}_j(u)|+|\na_v\mathcal{F}_j(u)|\\&&+|\na^2_v
\mathcal{F}_j(u)|)^2,\eeno which implies that for the case of
$\gamma+2s>0$ and $\tilde{N}_1\ge \gamma+2s $, \beno I
&\lesssim&2^{2js}2^{-2j} \|\mathcal{G}\|_{L^1(\R^3_v)}
\|\mathcal{F}_j\|_{H^2(\R^3_v)}^2,\eeno and for the case of
$\gamma+2s\le0$ and $\tilde{N}_1\ge -|\gamma+2s|$,  \beno I
&\lesssim&2^{2js}2^{-2j}
(\|\mathcal{G}\|_{L^1(\R^3_v)}+\|\mathcal{G}\|_{L^\frac{3}{2}(\R^3_v)})
\|\mathcal{F}_j\|_{H^2(\R^3_v)}^2.\eeno The similar estimate can be
applied to $II$ and it gives that for the case of $\gamma+2s>0$,
\beno II &\lesssim&2^{2js}2^{-2j} \|\mathcal{G}\|_{L^1(\R^3_v)}
\|\mathcal{H}_k\|_{L^2(\R^3_v)}^2,\eeno and for the case of
$\gamma+2s\le0$,  \beno II &\lesssim&2^{2js}2^{-2j}
(\|\mathcal{G}\|_{L^1(\R^3_v)}+\|\mathcal{G}\|_{L^\frac{3}{2}(\R^3_v)})
\|\mathcal{H}_k\|_{L^2(\R^3_v)}^2.\eeno Then the fact
 $\Gamma_{1,2}\le I^\f12 II^\f12$ and the Berstein's inequality lead to the estimate for
 $\Gamma_{1,2}$ which is exactly as the same as the case of
 $\Gamma_{1,1}$. We complete the proof to the Lemma. \ef

\medskip

 \begin{lem}\label{le2} If $\gamma+2s>0$, there holds \ben
|\Gamma_2|\lesssim 2^{2js} \|\mathcal{G}\|_{L^1(\R^3_v)}
\|\mathcal{H}_k\|_{L^2(\R^3_v)}\|\mathcal{F}_j\|_{L^2(\R^3_v)}, \een
and if $\gamma+2s\le 0$, there holds \ben |\Gamma_2|\lesssim 2^{2js}
 (\|\mathcal{G}\|_{L^1(\R^3_v)}+\|\mathcal{G}\|_{L^\frac{3}{2}(\R^3_v)})
\|\mathcal{H}_k\|_{L^2(\R^3_v)}\|\mathcal{F}_j\|_{L^2(\R^3_v)}. \een
\end{lem}

\noindent {\bf Proof:} Due to the cut-off function $\phi_j$, it is
easy to check that there is no singularity caused by the collision
kernel. Then one may estimate $\Gamma_{2}$ directly and we only
present the proof to bound the quantity \beno \Gamma_{2,1}
&\eqdefa&\int_{\R^3}dv\int_{\R^3}dv_*\int_{\SS^{2}}B(v-v_*,\sigma)\langle
v_*\rangle^{-N_1}\langle v\rangle^{-N_2}\mathcal{G}_*\mathcal{H}_k
[1-\phi_j(|v-v'|)] \\&& \mathbf{1}_{2^{-j}|v-v_*|^{-1}\le
c\f{\pi}{4}}\langle v'\rangle^{-N_3}\mathcal{F}'_jd\sigma,\eeno By
Cauchy-Schwartz inequality, it can be reduced to bound

\beno III
&\eqdefa&\int_{\R^3}dv\int_{\R^3}dv_*\int_{\SS^{2}}B(v-v_*,\sigma)\langle
v_*\rangle^{-N_1}\langle v\rangle^{-N_2}\mathcal{G}_*\mathcal{H}_k^2
[1-\phi_j(|v-v'|)]
\\&& \mathbf{1}_{2^{-j}|v-v_*|^{-1}\le c\f{\pi}{4}}\langle
v'\rangle^{-N_3} d\sigma,\eeno and \beno IV
&\eqdefa&\int_{\R^3}dv\int_{\R^3}dv_*\int_{\SS^{2}}B(v-v_*,\sigma)\langle
v_*\rangle^{-N_1}\langle v\rangle^{-N_2}\mathcal{G}_*
[1-\phi_j(|v-v'|)] \\&& \mathbf{1}_{2^{-j}|v-v_*|^{-1}\le
c\f{\pi}{4}}\langle v'\rangle^{-N_3}|\mathcal{F}'_j|^2d\sigma. \eeno

We choose to bound the quantity $IV$ since it is a little more
complicated than $III$. Fixed $\sigma$ and $v_*$, we perform the
change of variables $v\rightarrow v'$ and by a direct calculation,
its  Jacobian determinant is $\f18(1+\langle\frac{v-v_*}{|v-v_*|},
\sigma\rangle)$ which corresponds to the case $\kappa=1$ in
\eqref{ja}. From which together with \eqref{wei}, we arrive at \beno
IV &\lesssim&\int_{\R^3}dv'\int_{\R^3}dv_*|v_*-v'|^\gamma\langle
v_*-v'\rangle^{-\tilde{N}_1}\mathcal{G}_*
\mathbf{1}_{2^{-j}|v'-v_*|^{-1}\le c\f{\pi}{4}}
\\&&\times|\mathcal{F}'_j|^2\int_{\SS^{2}}b(\langle\frac{v-v_*}{|v-v_*|},
\sigma\rangle)[1-\phi_j(|v-v'|)]d\sigma. \eeno Noting that \beno
&&\int_{\SS^2}b(\langle\frac{v-v_*}{|v-v_*|},
\sigma\rangle)[1-\phi_j(|v-v'|)]d\sigma\mathbf{1}_{2^{-j}|v'-v_*|^{-1}\le
c\f{\pi}{4}}\\&&\quad\lesssim \int_{\sqrt{\f{1-\langle
\frac{v-v_*}{|v-v_*|},\sigma\rangle}{2}}\gtrsim 2^{-j}|v'-v_*|^{-1}}
b(\langle \frac{v-v_*}{|v-v_*|},\,\sigma\rangle)
 d\sigma\\&&\quad\lesssim\int_{c2^{-j}|v'-v_*|^{-1}}^{\f{\pi}2}\theta^{-1-2s}d\theta
\lesssim (2^{-j}|v'-v_*|^{-1})^{-2s},\eeno one has \beno IV\lesssim
2^{2js} \int_{\R^3}dv'\int_{\R^3}dv_*|v_*-v'|^{\gamma+2s}\langle
v_*-v'\rangle^{-\tilde{N}_1}\mathcal{G}_*
  |\mathcal{F}'_j|^2. \eeno
When $\gamma+2s>0$, there holds \beno IV\lesssim 2^{2js}
\|\mathcal{G}\|_{L^1(\R^3_v)}
\|\mathcal{F}_j\|_{L^2(\R^3_v)}^2.\eeno While $\gamma+2s\le0$, we
obtain that
 \beno IV\lesssim 2^{2js}
(\|\mathcal{G}\|_{L^1(\R^3_v)}+\|\mathcal{G}\|_{L^\frac{3}{2}(\R^3_v)})
\|\mathcal{F}_j\|_{L^2(\R^3_v)}^2. \eeno  The fact that
$|\Gamma_{2,1}|\lesssim III^{\f12}IV^{\f12}$ implies $\Gamma_{2,1}$
enjoys the same estimate as in the Lemma and it is enough to show
that  the Lemma also holds true. \ef

\begin{lem}\label{le3} If $\gamma+2s>0$, there holds \ben
|\Gamma_3|\lesssim 2^{2js} \|\mathcal{G}\|_{L^1(\R^3_v)}
\|\mathcal{H}_k\|_{L^2(\R^3_v)}\|\mathcal{F}_j\|_{L^2(\R^3_v)}, \een
and if $\gamma+2s\le 0$, there holds \ben |\Gamma_3|\lesssim
2^{2js}(\|\mathcal{G}\|_{L^1(\R^3_v)}+\|\mathcal{G}\|_{L^\frac{3}{2}(\R^3_v)})
 \|\mathcal{H}_k\|_{L^2(\R^3_v)}\|\mathcal{F}_j\|_{L^2(\R^3_v)}. \een
\end{lem}

\noindent {\bf Proof:} Similarly, to overcome the strong singularity
caused by collision kernel, we divide $\Gamma_3$ into two parts:
$\Gamma_{3,1}$ and $\Gamma_{3,2}$ which defined as \beno
\Gamma_{3,1}
&\eqdefa&\int_{\R^3}dv\int_{\R^3}dv_*\int_{\SS^{2}}B(v-v_*,\sigma)\langle
v_*\rangle^{-N_1}\langle v\rangle^{-N_2}\mathcal{G}_*\mathcal{H}_k
  \\&& \mathbf{1}_{2^{-j}|v-v_*|^{-1}\ge
c\f{\pi}{4}}(v'-v)\na_v\mathbb{F}_jd\sigma,\eeno and \beno
\Gamma_{3,2}&\eqdefa&\int_0^1d\kappa\int_{\R^3}dv\int_{\R^3}dv_*\int_{\SS^{2}}B(v-v_*,\sigma)\langle
v_*\rangle^{-N_1}\langle v\rangle^{-N_2}\mathcal{G}_*\mathcal{H}_k
  \\&& \mathbf{1}_{2^{-j}|v-v_*|^{-1}\ge
c\f{\pi}{4}}(v'-v)\otimes(v'-v):\na_v^2\mathbb{F}_j(\gamma(\kappa))d\sigma,\eeno
where we use the notation \eqref{F}. Observing that
 \beno
&&|\int_{\SS^2}b(\langle\frac{v-v_*}{|v-v_*|},
\sigma\rangle)(v'-v)d\sigma|\mathbf{1}_{2^{-j}|v-v_*|^{-1}\ge
c\f{\pi}{4}}\\&&\quad\lesssim |v-v_*|\int_0^
\f{\pi}2\theta^{1-2s}d\theta\mathbf{1}_{2^{-j}|v-v_*|^{-1}\ge
c\f{\pi}{4}}\lesssim |v-v_*|(2^{-j}|v-v_*|^{-1})^{2-2s}.\eeno
 then  for the term $\Gamma_{3,1}$, one has
 \beno |\Gamma_{3,1}|&\lesssim& 2^{-2j}2^{2js}\int_{\R^3}dv\int_{\R^3}dv_*|v-v_*|^{\gamma+2s-1}\langle
v-v_*\rangle^{-\tilde{N}_1}\mathcal{G}_*|\mathcal{H}_k|\\&&\quad\times
\mathbf{1}_{2^{-j}|v-v_*|^{-1}\ge
c\f{\pi}{4}}(|\mathcal{F}_j|+|\na_v\mathcal{F}_j|).
 \eeno
For the case of $\gamma+2s>0$, one may obtain that there exists a
constant $\delta\in (0,1)$ such that\beno |\Gamma_{3,1}|&\lesssim&
2^{-2j}2^{2js}\int_{\R^3}dv\int_{\R^3}dv_*|v-v_*|^{-1+\delta}
\mathcal{G}_*|\mathcal{H}_k| \mathbf{1}_{2^{-j}|v-v_*|^{-1}\ge
c\f{\pi}{4}}\\&&\quad\times
2^{-\f{j}2}|v-v_*|^{-\f12}(|\mathcal{F}_j|+|\na_v\mathcal{F}_j|)\\&
 \lesssim& 2^{2js}2^{-\f52j} \|\mathcal{G}\|_{L^1(\R^3_v)}
\|\mathcal{H}_k\|_{L^2(\R^3_v)}(\|\mathcal{F}_j\|_{L^\infty(\R^3_v)}+\|\na_v\mathcal{F}_j\|_{L^\infty(\R^3_v)})
\\&
 \lesssim& 2^{2js}  \|\mathcal{G}\|_{L^1(\R^3_v)}
\|\mathcal{H}_k\|_{L^2(\R^3_v)}\|\mathcal{F}_j\|_{L^2(\R^3_v)},
 \eeno
where  we use the Bernstein's inequality in the last inequality.

While for the case of  $\gamma+2s\le0$, we arrives at
 \beno
|\Gamma_{3,1}|&\lesssim& 2^{-2j}2^{2js}\int_{\R^3}dv\int_{\R^3}dv_*
|v-v_*|^{\gamma+2s-1}  \mathcal{G}_*|\mathcal{H}_k|\\&&\quad\times
\mathbf{1}_{2^{-j}|v-v_*|^{-1}\ge c\f{\pi}{4}}
 (|\na_v\mathcal{F}_j|+|\mathcal{F}_j|)\\&
 \lesssim& 2^{2js}2^{-2j}
  \|\mathcal{G}\|_{L^\frac{3}{2}(\R^3_v)}
\|\mathcal{H}_k\|_{L^2(\R^3_v)}
(\|\na_v\mathcal{F}_j\|_{L^6(\R^3_v)}+\|
\mathcal{F}_j\|_{L^6(\R^3_v)})
\\&
 \lesssim& 2^{2js}   \|\mathcal{G}\|_{L^\frac{3}{2}(\R^3_v)}
\|\mathcal{H}_k\|_{L^2(\R^3_v)}\|\mathcal{F}_j\|_{L^2(\R^3_v)}.
 \eeno

Next we shall focus on the estimate for $\Gamma_{3,2}$. Thanks to
the estimate \eqref{F_2}, one may obtain \beno
|\Gamma_{3,2}|&\lesssim&\int_0^1d\kappa\int_{\R^3}dv\int_{\R^3}dv_*\int_{\SS^{2}}B(v-v_*,\sigma)\langle
v_*\rangle^{-N_1}\langle v\rangle^{-N_2}\mathcal{G}_*|\mathcal{H}_k|
  \\&& \mathbf{1}_{2^{-j}|v-v_*|^{-1}\ge
c\f{\pi}{4}}|v'-v|^2\langle\gamma(\kappa)\rangle^{-N_3}
 \sum_{i=0}^2|\na_v^i\mathcal{F}_j(\gamma(\kappa))|d\sigma.\eeno

Due to the Cauchy-Schwartz's inequality, we only have to bound the
quantity as \beno
V&\eqdefa&\int_0^1d\kappa\int_{\R^3}dv\int_{\R^3}dv_*\int_{\SS^{2}}B(v-v_*,\sigma)\langle
v_*\rangle^{-N_1}\langle v\rangle^{-N_2}\mathcal{G}_*
  \\&& \mathbf{1}_{2^{-j}|v-v_*|^{-1}\ge
c\f{\pi}{4}}|v'-v|^2\langle \gamma(\kappa)\rangle^{-N_3}
(\sum_{i=0}^2|\na_v^i\mathcal{F}_j(\gamma(\kappa))|)^2d\sigma,\eeno
By change of variables from $v\rightarrow \gamma(\kappa)=u$ and the
fact \eqref{ja} and \eqref{ks}, the term $V$ can be bounded as \beno
V&\lesssim&
 \int_{\R^3}du\int_{\R^3}dv_*|v_*-u|^{\gamma+2}\langle
u-v_*\rangle^{-\tilde{N}_1} \mathcal{G}_*
   \mathbf{1}_{2^{-j}|u-v_*|^{-1}\ge
c\f{\pi}{4}}\\&&\times
(\sum_{i=0}^2|\na_v^i\mathcal{F}_j(u)|)^2\int_{\SS^{2}}b(\langle\frac{v-v_*}{|v-v_*|},
\sigma\rangle)\f{1-\langle\frac{v-v_*}{|v-v_*|},
\sigma\rangle}{2}d\sigma\\&\lesssim&
 \int_{\R^3}du\int_{\R^3}dv_*|v_*-u|^{\gamma+2}\langle
u-v_*\rangle^{-\tilde{N}_1} |\mathcal{G}_*|
   \mathbf{1}_{2^{-j}|u-v_*|^{-1}\ge
c\f{\pi}{4}}\\&&\times (\sum_{i=0}^2|\na_v^i\mathcal{F}_j(u)|)^2
(2^{-j}|u-v_*|^{-1})^{2-2s}.\eeno

For the case of $\gamma+2s>0$, one may get \beno V\lesssim
2^{-2j}2^{2js}\|\mathcal{G}\|_{L^1(\R^3_v)}
 \|\mathcal{F}_j\|_{H^2(\R^3_v)}^2.\eeno
While for the case of $-1<\gamma+2s\le0$, we obtain \beno
V&\lesssim&2^{-2j}2^{2js}\int_{\R^3}du\int_{\R^3}dv_*(1+|v_*-u|^{\gamma+2s}\mathbf{1}_{|v_*-u|\le
1}) \mathcal{G}_*     |\na_v^2\mathcal{F}_j(u)|)^2
 \\&\lesssim& 2^{-2j}2^{2js}
(\|\mathcal{G}\|_{L^1(\R^3_v)}+\|\mathcal{G}\|_{L^\frac{3}{2}(\R^3_v)})
 \|\mathcal{F}_j\|_{H^2(\R^3_v)}^2.\eeno
From which, we can deduce that the case of $\gamma+2s>0$, there
holds \beno |\Gamma_{3,2}|&\lesssim& 2^{2js}
\|\mathcal{G}\|_{L^1(\R^3_v)}
\|\mathcal{H}_k\|_{L^2(\R^3_v)}\|\mathcal{F}_j\|_{L^2(\R^3_v)},
 \eeno and in the case of $\gamma+2s\le0$,
there holds
 \beno
|\Gamma_{3,2}|&\lesssim&  2^{2js}
(\|\mathcal{G}\|_{L^1(\R^3_v)}+\|\mathcal{G}\|_{L^\frac{3}{2}(\R^3_v)})
\|\mathcal{H}_k\|_{L^2(\R^3_v)}\|\mathcal{F}_j\|_{L^2(\R^3_v)}.
 \eeno
From which, we complete the proof to the Lemma. \ef

Putting Lemma \ref{le1}, Lemma \ref{le2} together with Lemma
\ref{le3} , we easily deduce that for the case of $\gamma+2s>0$,
there holds \ben \sum_{j<k}|\langle \mathcal
{Q}(\mathcal{G},\mathcal{H}_k),\mathcal{F}_j\rangle_v|&\lesssim&
\sum_{j<k} 2^{2js} \|\mathcal{G}\|_{L^1(\R^3_v)}
\|\mathcal{H}_k\|_{L^2(\R^3_v)}\|\mathcal{F}_j\|_{L^2(\R^3_v)}\nonumber\\
&\lesssim&
\|\mathcal{G}\|_{L^1(\R^3_v)}\sum_{j<k}2^{(j-k)s}\|\mathcal{H}_k\|_{H^s(\R^3_v)}\|\mathcal{F}_j\|_{H^s(\R^3_v)}\nonumber
\\
&\lesssim&
\|\mathcal{G}\|_{L^1(\R^3_v)}\|\mathcal{H}\|_{H^s(\R^3_v)}\|\mathcal{F}\|_{H^s(\R^3_v)}.\label{q1}
\een  Similarly, for the case of $\gamma+2s\le0$, there holds there
holds \ben  \sum_{j<k}|\langle {Q}(g,h_k),f_j\rangle_v|&\lesssim&
(\|\mathcal{G}\|_{L^1(\R^3_v)}+\|\mathcal{G}\|_{L^\frac{3}{2}(\R^3_v)})\|\mathcal{H}\|_{H^s(\R^3_v)}\|\mathcal{F}\|_{H^s(\R^3_v)}.\label{q2}
\een

 \noindent{\bf  Step 2: Frequency dominated by    the function
$\mathcal{F}$.} Since the Boltzmann collision operator can be
regrades as  the bilinear operator, we believe that in some sense
there exists symmetric  structure inside the operator. We recall
that \beno
 &&\langle\mathcal {Q}(\mathcal{G},\mathcal{H}),\mathcal{F}\rangle_v\\&=&
\int_{\R^3}dv\int_{\R^3}dv_*\int_{\SS^{2}}B(v-v_*,\sigma)\langle
v_*\rangle^{-N_1}\langle v\rangle^{-N_2}\mathcal{G}_*\mathcal{H}
\bigg(\langle v'\rangle^{-N_3}\mathcal{F}'-\langle
v\rangle^{-N_3}\mathcal{F}\bigg)d\sigma.\eeno It is easy to check
that \beno \langle\mathcal
{Q}(\mathcal{G},\mathcal{H}),\mathcal{F}\rangle_v=\langle\mathcal
{Q}_1(\mathcal{G},\mathcal{H}),\mathcal{F}\rangle_v +\langle\mathcal
{Q}_2(\mathcal{G},\mathcal{H}),\mathcal{F}\rangle_v,\eeno  where
\beno &&\langle\mathcal
{Q}_1(\mathcal{G},\mathcal{H}),\mathcal{F}\rangle_v\\&=&
\int_{\R^3}dv\int_{\R^3}dv_*\int_{\SS^{2}}B(v-v_*,\sigma)\langle
v_*\rangle^{-N_1}\mathcal{G}_*\bigg(\langle
v\rangle^{-N_2}\mathcal{H}-\langle
v'\rangle^{-N_2}\mathcal{H}'\bigg)\langle
v'\rangle^{-N_3}\mathcal{F}' d\sigma,\eeno and \beno
&&\langle\mathcal
{Q}_2(\mathcal{G},\mathcal{H}),\mathcal{F}\rangle_v\\&=&
\int_{\R^3}dv\int_{\R^3}dv_*\int_{\SS^{2}}B(v-v_*,\sigma)\langle
v_*\rangle^{-N_1}\mathcal{G}_*\bigg(\langle
v'\rangle^{-N_2}\mathcal{H}'\langle
v'\rangle^{-N_3}\mathcal{F}'-\langle v\rangle^{-N_2}
\mathcal{H}\langle v\rangle^{-N_3}\mathcal{F}\bigg).\eeno

We remark that in the case of $j\le k$, $\langle\mathcal
{Q}_1(\mathcal{G},\mathcal{H}_j),\mathcal{F}_k\rangle_v$ enjoys some
similar  structure as $\langle\mathcal
{Q}(\mathcal{G},\mathcal{H}_k),\mathcal{F}_j\rangle_v$. Precisely,
one may follow the same procedure to handle the inner product
 $\langle\mathcal {Q}_1(\mathcal{G},\mathcal{H}_j),\mathcal{F}_k\rangle_v$ as we did
 for $\langle\mathcal  {Q}(\mathcal{G},\mathcal{H}_k),\mathcal{F}_j\rangle_v$.
We point out that the main difference  lies in the Taylor expansion.
If we set $\mathbb{H}_j=\langle v\rangle^{-N_2}\mathcal{H}_j $, then
in this case, the Taylor expansion should be taken as
\ben\label{ta2}
\mathbb{H}_j(v)-\mathbb{H}_j(v')=(v-v')\cdot\na_v\mathbb{H}_j(v')+\int_0^1
(v-v')\otimes(v-v'):\na_v^2 \mathbb{H}(\gamma(\kappa))d\kappa, \een
where $\gamma(\kappa)=\kappa v'+(1-\kappa)v$.

 Another difference comes from the following  fact. For each $\sigma$
and $v_*$, let $\psi_\sigma(v')$ represents the inverse transform
$v'\rightarrow \psi_\sigma(v')=v$(see \cite{advw}).  Then due to
\eqref{ja}, one has \ben\label{ja1}
\big|\frac{dv'}{dv}\big|=\f{\langle\frac{v'-v_*}{|v'-v_*|},\sigma\rangle^2}{4}.\een
 Thus for fixed $v_*$ and smooth  function $\phi$, one has
 \ben\label{van}
&&\int_{\R^3}dv
\int_{\SS^{2}}B(v-v_*,\sigma)g(v')\phi(|v-v'|)(v-v')d\sigma\nonumber\\&&=
4\int_{\R^3}dv \int_{\langle\frac{v-v_*}{|v-v_*|},\sigma\rangle\ge
\frac{1}{\sqrt{2}}}B(\psi_\sigma(v)-v_*,\sigma)g(v)\phi(|\psi_\sigma(v)-v|)\f{\psi_\sigma(v)-v}{\langle\frac{v-v_*}{|v-v_*|},\sigma\rangle^2}
d\sigma\nonumber\\&&=0, \een  where we use the symmetric property of
$\psi_\sigma(v)$ with respect to $\sigma$. The fact will give the
reduction for proof to the corresponding terms such as
$\Gamma_{1,1}$ and $\Gamma_{3,1}$.

Fortunately, the differences mentioned before do harmless to the
proof for ${Q}_1(\mathcal{G},\mathcal{H}_j),\mathcal{F}_k\rangle_v$
as we did in the step 1. We omit the details here and obtain that
for the case of $\gamma+2s>0$, there holds \ben\label{q3} \sum_{j\le
k}|\langle \mathcal
{Q}_1(\mathcal{G},\mathcal{H}_j),\mathcal{F}_k\rangle_v| &\lesssim&
\|\mathcal{G}\|_{L^1(\R^3_v)}\|\mathcal{H}\|_{H^s(\R^3_v)}\|\mathcal{F}\|_{H^s(\R^3_v)},
\een  and for the case of $\gamma+2s\le0$, there holds
\ben\label{q4} \sum_{j\le k}|\langle \mathcal
{Q}_1(\mathcal{G},\mathcal{H}_j),\mathcal{F}_k\rangle_v|&\lesssim&
(\|\mathcal{G}\|_{L^1(\R^3_v)}+\|\mathcal{G}\|_{L^\frac{3}{2}(\R^3_v)})\|\mathcal{H}\|_{H^s(\R^3_v)}\|\mathcal{F}\|_{H^s(\R^3_v)}.
\een

Now we turn to estimate
${Q}_2(\mathcal{G},\mathcal{H}_j),\mathcal{F}_k\rangle_v$. Actually,
by change of variables and cancellation lemma(see \cite{advw}), it
can be written as \beno &&\langle\mathcal
{Q}_2(\mathcal{G},\mathcal{H}_j),\mathcal{F}_k\rangle_v\\&=& |\SS^1|
\int_{\R^6}\int_0^\frac{\pi}{2}
\sin\theta\bigg(\f1{\cos^3\frac{\theta}{2}}B(\frac{|v-v_*|}{\cos\frac{\theta}{2}},\cos\theta)-B(|v-v_*|,\cos\theta)\bigg)\\&&\times
\langle v_*\rangle^{-N_1}\mathcal{G}_* \langle v\rangle^{-N_2}
\mathcal{H}\langle v\rangle^{-N_3}\mathcal{F} d\theta dv dv_*.\eeno

It is easy to check that  \beno
&&B(\frac{|v-v_*|}{\cos^3\frac{\theta}{2}},\cos\theta)-B(|v-v_*|,\cos\theta)\\&&=b(\cos\theta)
|v-v_*|^\gamma\{(\cos^{-1}\f{\theta}2)^{\gamma+3}-1\}.\eeno

Using the fact  that \beno x^m-y^m= m\int_y^x z^{m-1}dz,\eeno we
obtain that for $\theta\in[0,\frac{\pi}{2}],$ \beno
(\cos^{-1}\f{\theta}2)^{\gamma+3}-1\lesssim
\sin^2\frac{\theta}{2},\eeno which immediately implies that \ben
&&\big|\langle\mathcal
{Q}_2(\mathcal{G},\mathcal{H}_j),\mathcal{F}_k\rangle_v\big|\nonumber\\&\lesssim&
  \int_{\R^6}\int_0^\frac{\pi}{2} \theta^{1-2s}|v-v_*|^\gamma
 \langle v- v_*\rangle^{-\tilde{N}_1}\mathcal{G}_*|\mathcal{H}_j \mathcal{F}_k| d\theta dv dv_*\nonumber\\&\lesssim&
  \int_{\R^6}  |v-v_*|^\gamma
 \langle v- v_*\rangle^{-\tilde{N}_1}\mathcal{G}_*|\mathcal{H}_j \mathcal{F}_k| dv dv_*\\&\eqdefa&
 \Xi_1+\Xi_2,\nonumber\een
where \beno \Xi_1= \int_{\R^6} |v-v_*|^\gamma
 \langle v- v_*\rangle^{-\tilde{N}_1}\mathcal{G}_*|\mathcal{H}_j \mathcal{F}_k |\mathbf{1}_{2^{-k}|v-v_*|^{-1}\ge c\frac{\pi}{4}} dv
 dv_*,\eeno
and \beno  \Xi_2= \int_{\R^6} |v-v_*|^\gamma
 \langle v- v_*\rangle^{-\tilde{N}_1}\mathcal{G}_*|\mathcal{H}_j \mathcal{F}_k |\mathbf{1}_{2^{-k}|v-v_*|^{-1}\le c\frac{\pi}{4}}
dv dv_*\eeno

Here we emphasize that  the analogue decomposition for the case
$0<s<\f12$ is radically different from the case of $\f12\le s <1$.
The quantity $2^{-k}|v-v_*|^{-1}$ inside the decomposition should be
replaced by $2^{-j}|v-v_*|^{-1}$.  We remark that here $j$
represents the low frequency   and $k$ represents the high
frequency.

As for the term $\Xi_1$, one has \beno \Xi_1&\lesssim&\int_{\R^6}
|v-v_*|^\gamma (2^{-k}|v-v_*|^{-1})^{2-2s}
 \langle v- v_*\rangle^{-\tilde{N}_1}\mathcal{G}_*|\mathcal{H}_j \mathcal{F}_k| \mathbf{1}_{2^{-k}|v-v_*|^{-1}\ge c\frac{\pi}{4}}  dv dv_*
 \\&\lesssim& 2^{2ks}2^{-2k}\int_{\R^6}
  |v-v_*|^{\gamma+2s-2}
 \langle v- v_*\rangle^{-\tilde{N}_1}\mathcal{G}_*|\mathcal{H}_j \mathcal{F}_k| \mathbf{1}_{2^{-k}|v-v_*|^{-1}\ge c\frac{\pi}{4}} dv dv_*\eeno

For the case of $\gamma+2s>0$, one has for some small $\delta\in
(0,2)$, \beno |v-v_*|^{\gamma+2s-2}
 \langle v- v_*\rangle^{-\tilde{N}_1}\le |v-v_*|^{\delta-2},\eeno
which leads to \beno \Xi_1&\lesssim&  2^{2ks}2^{-2k}\int_{\R^6}
  |v-v_*|^{\delta-2}
 \mathcal{G}_*|\mathcal{H}_j \mathcal{F}_k|
\mathbf{1}_{2^{-k}|v-v_*|^{-1}\ge c\frac{\pi}{4}} dv
dv_*\\&\lesssim&2^{2ks}2^{-2k}\|\mathcal{G}\|_{L^1(\R^3_v)}\|\mathcal{H}_j\|_{L^6(\R^3_v)}\|\mathcal{F}_k\|_{L^6(\R^3_v)}\\
&\lesssim&
2^{-(k-j)(1-s)}\|\mathcal{G}\|_{L^1(\R^3_v)}\|\mathcal{H}_j\|_{H^s(\R^3_v)}\|\mathcal{F}_k\|_{H^s(\R^3_v)},\eeno
where we use the fact $2^{-k}|v-v_*|^{-1}\ge c\frac{\pi}{4}$ and the
Bernstein's inequality.

While for the case of $-1<\gamma+2s\le 0$, one has \beno
\Xi_1&\lesssim& 2^{2ks}2^{-2k}\int_{\R^6}
   |v-v_*|^{\gamma+2s-2}
\mathcal{G}_*|\mathcal{H}_j \mathcal{F}_k|
\mathbf{1}_{ |v-v_*| \le c^{-1}\frac{4}{\pi}} dv dv_*\\
&\lesssim&
2^{2ks}2^{-2k}\|\mathcal{G}\|_{L^\f32(\R^3_v)}\|\mathcal{H}_j\|_{L^6(\R^3_v)}\|\mathcal{F}_k\|_{L^6(\R^3_v)}
 \\&\lesssim&
2^{-(k-j)(1-s)}\|\mathcal{G}\|_{L^2(\R^3_v)}\|\mathcal{H}_j\|_{H^s(\R^3_v)}\|\mathcal{F}_k\|_{H^s(\R^3_v)}.
 \eeno

We turn to bound the quantity $\Xi_2$.   In the case of
$2^{-k}|v-v_*|^{-1}\le c\frac{\pi}{4}$, one may deduce that \beno
1\lesssim (2^{-k}|v-v_*|^{-1})^{1-2s}.\eeno Then we arrive at \beno
\Xi_2&\lesssim& 2^{2ks}2^{-k} \int_{\R^6}
  |v-v_*|^{\gamma+2s-1}
 \langle v- v_*\rangle^{-\tilde{N}_1}\mathcal{G}_*|\mathcal{H}_j \mathcal{F}_k|\mathbf{1}_{2^{-j}|v-v_*|^{-1}\ge c\frac{\pi}{4}}  dv dv_*
  \eeno

For the case of $\gamma+2s>0$, one has for some small $\delta\in
(0,1)$, \beno |v-v_*|^{\gamma+2s-1}
 \langle v- v_*\rangle^{-\tilde{N}_1}\le |v-v_*|^{\delta-1},\eeno
which leads to \beno \Xi_2&\lesssim& 2^{2ks}2^{-k}\int_{\R^6}
  |v-v_*|^{\delta-1}
 \mathcal{G}_*|\mathcal{H}_j \mathcal{F}_k |dv
dv_*\\&\lesssim&2^{2ks}2^{-k} \|\mathcal{G}\|_{L^1(\R^3_v)}\|\mathcal{H}_j\|_{L^6(\R^3_v)}\|\mathcal{F}_k\|_{L^2(\R^3_v)}\\
&\lesssim&
2^{-(k-j)(1-s)}\|\mathcal{G}\|_{L^1(\R^3_v)}\|\mathcal{H}_j\|_{H^s(\R^3_v)}\|\mathcal{F}_k\|_{H^s(\R^3_v)}.\eeno

While for the case of $-1<\gamma+2s\le 0$, one has \beno
|v-v_*|^{\gamma+2s-1}
 \langle v- v_*\rangle^{-\tilde{N}_1}\le 1+|v-v_*|^{\gamma+2s-1}\mathbf{1}_{|v-v_*|\le 1}.\eeno Then it gives \beno
\Xi_2&\lesssim&  2^{2ks}2^{-k}
(|\mathcal{G}\|_{L^1(\R^3_v)}+\|\mathcal{G}\|_{L^\f32(\R^3_v)})(\|\mathcal{H}_j\|_{L^6(\R^3_v)}+
\|\mathcal{H}_j\|_{L^2(\R^3_v)})\|\mathcal{F}_k\|_{L^2(\R^3_v)}\\&\lesssim&
 2^{-(k-j)(1-s)}(|\mathcal{G}\|_{L^1(\R^3_v)}+\|\mathcal{G}\|_{L^\f32(\R^3_v)}) \|\mathcal{H}_j\|_{H^s(\R^3_v)}\|\mathcal{F}_k\|_{H^s(\R^3_v)}.
 \eeno

Patch together the estimates before, we finally obtain that for the
case of $\gamma+2s>0$, there holds \ben\label{q5} \sum_{j\le
k}|\langle \mathcal
{Q}_2(\mathcal{G},\mathcal{H}_j),\mathcal{F}_k\rangle_v| &\lesssim&
\|\mathcal{G}\|_{L^1(\R^3_v)}\|\mathcal{H}\|_{H^s(\R^3_v)}\|\mathcal{F}\|_{H^s(\R^3_v)},
\een  and for the case of $\gamma+2s\le0$, there holds
\ben\label{q6} \sum_{j\le k}|\langle \mathcal
{Q}_2(\mathcal{G},\mathcal{H}_j),\mathcal{F}_k\rangle_v|&\lesssim&
(\|\mathcal{G}\|_{L^1(\R^3_v)}+\|\mathcal{G}\|_{L^\f32(\R^3_v)})\|\mathcal{H}\|_{H^s(\R^3_v)}\|\mathcal{F}\|_{H^s(\R^3_v)}.
\een

Thanks to the \eqref{q1}, \eqref{q3} and \eqref{q5}, we obtain the
first estimate in the Theorem \ref{ub}. Similarly,  \eqref{q2},
\eqref{q4} and \eqref{q6} imply the second one. \ef

In order to get the regularizing effect of the solutions for the
inhomogeneous Boltzmann equation, we still need the estimate for
some commutator. Fortunately, we can follow the same idea of the
proof to the Theorem and then obtain the corollary:
\begin{col}\label{comm} Let $ N_1=|N_2|+|N_3|+\max\{|l-2|, |l-1|\}$ and
$\tilde{N}_1=N_2+N_3$ with $N_2, N_3, l \in\R$. Then if
$\tilde{N}_1\ge l+\gamma$ and $s<\f12$, one has \ben
 &&|\langle Q(g,h)\langle v\rangle^{l},f\rangle_v-\langle
Q(g,h\langle v\rangle^{l}),f\rangle_v|\\\quad&&\lesssim
 \|g\|_{L^1_{N_1}(\R^3_v)} \|h\|_{H^{\varrho}_{N_2}(\R^3_v)}
\|f\|_{H^{\varrho}_{N_3}(\R^3_v)},\nonumber\een where $\varrho<s$.

 When $\tilde{N}_1\ge l-1+\gamma+2s$ and $s\ge\f12$,
one has that in the case of $\gamma+2s>0$, there holds \ben
 &&|\langle Q(g,h)\langle v\rangle^{l},f\rangle_v-\langle
Q(g,h\langle v\rangle^{l}),f\rangle_v|\\\quad&&\lesssim
 \|g\|_{L^1_{N_1}(\R^3_v)} \|h\|_{H^s_{N_2}(\R^3_v)}
\|f\|_{L^2_{N_3}(\R^3_v)}.\nonumber\een While in the case of
$\gamma+2s\le0$, there holds \ben
 &&|\langle Q(g,h)\langle v\rangle^{l},f\rangle_v-\langle
Q(g,h\langle v\rangle^{l}),f\rangle_v|\\\quad&&\lesssim
 (\|g\|_{L^1_{N_1}(\R^3_v)}+\|g\|_{L^\f32_{N_1}(\R^3_v)}) \|h\|_{H^s_{N_2}(\R^3_v)}
\|f\|_{L^2_{N_3}(\R^3_v)}.\nonumber\een
\end{col}

\noindent{\bf Proof:} Direct calculation  gives that \beno &&\langle
Q(g,h)\langle v\rangle^{l},f\rangle_v-\langle Q(g,h\langle
v\rangle^{l}),f\rangle_v\\&=&
\int_{\R^6}dvdv_*\int_{\SS^{2}}B(v-v_*,\sigma)\langle
v_*\rangle^{-N_1}\langle v\rangle^{-N_2}\mathcal{G}_*\mathcal{H}
\langle v'\rangle^{-N_3}\mathcal{F}'(\langle v'\rangle^{l}-\langle
v\rangle^{l})d\sigma.\eeno

\noindent{\bf Step 1: Bounds in the case of $s<\f12$.}  It is easy
to check \beno &&|\langle c Q(g,h)\langle
v\rangle^{l},f\rangle_v-\langle Q(g,h\langle
v\rangle^{l}),f\rangle_v|\\&\lesssim&
 \int_{\R^6}dvdv_*\int_{\SS^{2}} |v-v_*|^{\gamma+1}b(\cos\theta)\sin
\theta \langle v-v_* \rangle^{-\tilde{N}_1+l-1}\mathcal{G}_*
\mathcal{H} \mathcal{F}' d\sigma.\eeno

By Cauchy-Schwartz's inequality and change of variables, the desired
estimate can be reduced to the boundness of the quantity  \beno
\mathcal{K}=\int_{\R^6}dvdv_* |v-v_*|^{\gamma+1} \langle v-v_*
\rangle^{-\tilde{N}_1+l-1}\mathcal{G}_*\mathcal{H}^2.\eeno Since
$\tilde{N}_1\ge l+\gamma$,  we deduce that  in the case of
$\gamma+1\ge0$, there holds \ben\label{com3} \mathcal{K}\lesssim
\|\mathcal{G}\|_{L^1(\R^3_v)}
 \|\mathcal{H}\|_{L^2(\R^3_v)}^2. \een

While in the case of $\gamma+1<0$, one has \beno \mathcal{K}\lesssim
\int_{\R^6}dvdv_*(1+ |v-v_*|^{\gamma+1}\mathbf{1}_{|v-v_*|\le 1})
 \mathcal{G}_*\mathcal{H}^2.\eeno
Thanks to the assumption that $\gamma+2s+1>0$, there exist a
positive $\delta$ such that $\gamma+1=-(2s-\delta)$ which implies
that
 $\chi(v)=\mathbf{1}_{|v|\le 1}|v|^{\gamma+1}\in L^{\frac{3-\epsilon}{2s-\delta}}$. Then by H\"{o}lder's inequality, one may
obtain that \beno  \mathcal{K} &\lesssim& \|\mathcal
{G}\|_{L^1(\R^3_v)}(\|\mathcal {H}\|_{L^2(\R^3_v)}^2+\|\mathcal
{H}^2\|_{L^q(\R^3_v)}\|\chi\|_{L^p(\R^3_v)}), \eeno where
$\f1p+\f1q=1$ and $p=\frac{3-\epsilon}{2s-\delta}$. By standard
Sobolev's embedding theorem, one has \beno \|\mathcal
{H}^2\|_{L^q(\R^3_v)}\le \|\mathcal {H}\|_{H^{\varrho}(\R^3_v)}^2,
\eeno where $\varrho=\frac{6s-3\delta}{6-2\epsilon}$. Choose
$\epsilon <\frac{3\delta}{2s}$, then we deduce that $\varrho<s$. And
we arrive at \ben\label{com4}  \mathcal{K}&\lesssim& \|\mathcal {G}
\|_{L^1(\R^3_v)}\|\mathcal {H} \|_{H^{\varrho}(\R^3_v)}^2.\een

Combing \eqref{com3}, \eqref{com4} and Cauchy-Schwartz's inequality,
we easily obtain the first estimate in the corollary.

\medskip

\noindent{\bf Step 2: Bounds in the case of $s\ge \f12$.} One has
the following decomposition: \beno &&\langle Q(g,h)\langle
v\rangle^{l},f\rangle_v-\langle Q(g,h\langle
v\rangle^{l}),f\rangle_v\\
&\eqdefa&\mathfrak{R}_1+\mathfrak{R}_2,\eeno where \beno
\mathfrak{R}_1=\int_{\R^6}dvdv_*\int_{\SS^{2}}B(v-v_*,\sigma)\langle
v_*\rangle^{-N_1}\langle
v'\rangle^{-\tilde{N}_1}\mathcal{G}_*\mathcal{H}'
 \mathcal{F}'\bigg(\langle v'\rangle^{l}-\langle v\rangle^{l}\bigg)d\sigma
\eeno and \beno
\mathfrak{R}_2=\int_{\R^6}dvdv_*\int_{\SS^{2}}B(v-v_*,\sigma)\langle
v_*\rangle^{-N_1}\langle v'\rangle^{-N_3}\mathcal{G}_*\mathcal{F}'
 \bigg(\langle v\rangle^{-N_2}\mathcal{H}-\langle v'\rangle^{-N_2}\mathcal{H}'\bigg)\bigg(\langle v'\rangle^{l}-\langle v\rangle^{l}
 \bigg)d\sigma.
 \eeno

Firstly, as for the $\mathfrak{R}_1$, we have the following lemma:
\begin{lem}\label{le4}
Let $\tilde{N}_1\ge l-1+ \gamma+2s $.     In the case of
$\gamma+2s>0$, there holds \beno  \mathfrak{R}_{1 }  \lesssim
\|\mathcal{G}\|_{L^1(\R^3_v)}
\|\mathcal{H}\|_{H^s(\R^3_v)}\|\mathcal{F}\|_{L^2(\R^3_v)}. \eeno
While in the case of $\gamma+2s\le 0$, there holds \beno
 \mathfrak{R}_{1 }  \lesssim
(\|\mathcal{G}\|_{L^1(\R^3_v)}+\|\mathcal{G}\|_{L^\f32(\R^3_v)})
\|\mathcal{H}\|_{H^s(\R^3_v)}\|\mathcal{F}\|_{L^2(\R^3_v)}. \eeno
\end{lem}

\noindent{\bf Proof:} To overcome the singularity caused by the
collision kernel, we use the Taylor expansion formula up to the
order 2: \beno \langle v\rangle^{l}-\langle
v'\rangle^{l}=(v-v')\cdot\na_v(\langle
v\rangle^{l})|_{v=v'}+\int_0^1 (v-v')\otimes(v-v'):\na_v^2 (\langle
v \rangle^{l})|_{v=\gamma(\kappa)}d\kappa, \eeno where
$\gamma(\kappa)=\kappa v'+(1-\kappa)v$. Set
\beno\mathfrak{R}_{1,1}=\int_{\R^6}dvdv_*\int_{\SS^{2}}B(v-v_*,\sigma)\langle
v_*\rangle^{-N_1}\langle
v'\rangle^{-\tilde{N}_1}\mathcal{G}_*\mathcal{H}'
 \mathcal{F}'(v'-v)\na_{v}(\langle v'\rangle^{l}) d\sigma. \eeno
Then we claim that $\mathfrak{R}_{1,1}=0$. In fact, by change of
variables form $v\rightarrow v'$ and fact \eqref{ja1}, one has \beno
\mathfrak{R}_{1,1}=\int_{\R^6}dv'dv_*\langle
v_*\rangle^{-N_1}\langle
v'\rangle^{-\tilde{N}_1}\mathcal{G}_*\mathcal{H}'
 \mathcal{F}'\na_{v}(\langle v'\rangle^{l}) \int_{\SS^{2}}\f{4}{\langle\frac{v'-v_*}{|v'-v_*|},\sigma\rangle^2}B(v-v_*,\sigma)(v'-v)d\sigma.\eeno
Thanks to \eqref{van}, we prove the claim that
$\mathfrak{R}_{1,1}=0$. Next we shall focus on the estimate to the
following term: \beno
\mathfrak{R}_{1,2}&=&-\int_0^1d\kappa\int_{\R^6}dvdv_*\int_{\SS^{2}}B(v-v_*,\sigma)\langle
v_*\rangle^{-N_1}\langle
v'\rangle^{-\tilde{N}_1}\\&&\times\mathcal{G}_*\mathcal{H}'
 \mathcal{F}' (v-v')\otimes(v-v'):\na_v^2(\langle v\rangle^{l})|_{v=\gamma(\kappa)} d\sigma. \eeno
It is easy to see that \beno
\mathfrak{R}_{1,2}&\lesssim&\int_0^1d\kappa\int_{\R^6}dvdv_*\int_{\SS^{2}}B(v-v_*,\sigma)\langle
v_*\rangle^{-N_1}\langle
v'\rangle^{-\tilde{N}_1}\mathcal{G}_*\mathcal{H}'
 \mathcal{F}'|v-v'|^2\langle \gamma(\kappa)  \rangle^{l-2} d\sigma\\
 &\lesssim& \int_{\R^6}dv'dv_*\int_0^\f{\pi}{2} \theta^{1-2s}d\theta\langle
v_*-v'\rangle^{-\tilde{N}_1+l-2}\mathcal{G}_*\mathcal{H}'
 \mathcal{F}'|v_*-v'|^{\gamma+2} , \eeno
where we change the variables from $v$ to $v'$ and use the fact
\ben\label{ks1} |v-v'|^2\lesssim
|v'-v_*|^2\frac{1-\langle\frac{v-v_*}{|v-v_*|},\sigma\rangle}{2} \een
and \beno  \langle v_*\rangle^{-N_1}\langle
v'\rangle^{-\tilde{N}_1}\langle \gamma(\kappa) \rangle^{l-2}\lesssim
\langle v_*-v'\rangle^{-\tilde{N}_1+l-2}. \eeno

Then we deduce that \beno \mathfrak{R}_{1,2}&\lesssim&\sum_{j
}\mathfrak{R}_{1,2}^{j } , \eeno where \beno
\mathfrak{R}_{1,2}^{j}&\eqdefa& \int_{\R^6}dv'dv_* \langle
v_*-v'\rangle^{-\tilde{N}_1+l-2}\mathcal{G}_*|\mathcal{H}'_j|
 \mathcal{F}'|v_*-v'|^{\gamma+2}. \eeno

To bound the term $\mathfrak{R}_{1,2}^{j}$,
  we first observe that in the region of $2^{-j}|v'-v_*|^{-1}\ge \pi/4$,
there holds \beno \mathfrak{R}_{1,2}^{j}&\lesssim&
\int_{\R^6}dv'dv_* (2^{-j}|v'-v_*|^{-1})^{2-2s}\langle
v_*-v'\rangle^{-\tilde{N}_1+l-2} \mathcal{G}_*|\mathcal{H}'_j|
 \mathcal{F}' |v_*-v'|^{\gamma+2}\\&\lesssim& 2^{-2j(1-s)}
\int_{\R^6}dv'dv_* \langle v_*-v'\rangle^{-\tilde{N}_1+l-2}
\mathcal{G}_*|\mathcal{H}'_j|
 \mathcal{F}' |v_*-v'|^{\gamma+2s}. \eeno

Choose $\tilde{N}_1\ge l-1+\gamma+2s$. Then for the case of
$\gamma+2s>0$, we obtain that \beno \mathfrak{R}_{1,2}^{j
}&\lesssim&
2^{-2(1-s)j}\|\mathcal{G}\|_{L^1(\R^3_v)}\|\mathcal{H}_j\|_{L^2(\R^3_v)}\|\mathcal{F}
\|_{L^2(\R^3_v)}.\eeno While in the case of $\gamma+2s\le 0$,  we
deduce that \beno \mathfrak{R}_{1,2}^{j }&\lesssim&
2^{-2(1-s)j}(\|\mathcal{G}\|_{L^1(\R^3_v)}+\|\mathcal{G}\|_{L^\f32(\R^3_v)})\|\mathcal{H}_j\|_{L^2(\R^3_v)}
\|\mathcal{F} \|_{L^2(\R^3_v)}.\eeno

Next we treat the case of $2^{- j}|v'-v_*|^{-1}\le \pi/4$. One has
  \beno \mathfrak{R}_{1,2}^{j }&\lesssim&
\int_{\R^6}dv'dv_* (2^{-j}|v'-v_*|^{-1})^{1-2s}\langle
v_*-v'\rangle^{-\tilde{N}_1+l-2} \mathcal{G}_*|\mathcal{H}'_j|
 \mathcal{F}'|v_*-v'|^{\gamma+2}\\&\lesssim& 2^{2js}2^{-j}
\int_{\R^6}dv'dv_* \langle v_*-v'\rangle^{-\tilde{N}_1+l-2}
\mathcal{G}_*|\mathcal{H}'_j|
 \mathcal{F}'|v_*-v'|^{\gamma+2s+1}, \eeno
which implies \beno \mathfrak{R}_{1,2}^{j }&\lesssim& 2^{-(1-s)j}
\|\mathcal{G}\|_{L^1(\R^3_v)}
\|\mathcal{H}_j\|_{H^s(\R^3_v)}\|\mathcal{F} \|_{L^2(\R^3_v)}.\eeno

Patch together all the estimates for $\mathfrak{R}_{1,2}^{j}$ and
choose $\tilde{N}_1\ge l-1+\gamma+2s$, we deduce that in the case of
$\gamma+2s>0$, there holds \beno \sum_{j }\mathfrak{R}_{1,2}^{j}
\lesssim  \|\mathcal{G}\|_{L^1(\R^3_v)}
\|\mathcal{H}\|_{H^s(\R^3_v)}\|\mathcal{F}\|_{L^2(\R^3_v)}. \eeno
While in the case of $\gamma+2s\le 0$, there holds \beno \sum_{j
}\mathfrak{R}_{1,2}^{j} \lesssim
(\|\mathcal{G}\|_{L^1(\R^3_v)}+\|\mathcal{G}\|_{L^\f32(\R^3_v)})
\|\mathcal{H}\|_{H^s(\R^3_v)}\|\mathcal{F}\|_{L^2(\R^3_v)}. \eeno
\ef

Now we turn to focus on the term $\mathfrak{R}_2$. One has
\begin{lem}\label{le5} Let $\tilde{N}_1\ge
l-1+\gamma+2s$. Then in the case of $\gamma+2s>0$, there holds \beno
\mathfrak{R}_{2 } &\lesssim&
 \|\mathcal{G}\|_{L^1(\R^3_v)}\|\mathcal{H}\|_{H^s(\R^3_v)}\|\mathcal{F}\|_{L^2(\R^3_v)}.\eeno
While in the case of $\gamma+2s\le 0$,  there holds \beno
\mathfrak{R}_{2 }  &\lesssim&
(\|\mathcal{G}\|_{L^1(\R^3_v)}+\|\mathcal{G}\|_{L^2(\R^3_v)})\|\mathcal{H}\|_{H^s(\R^3_v)}\|\mathcal{F}\|_{L^2(\R^3_v)}.
\eeno
\end{lem}
\noindent{\bf Proof:} We introduce the Littlewood-Paley
decomposition and set \beno
\mathfrak{R}_2^j&=&\int_{\R^6}dvdv_*\int_{\SS^{2}}B(v-v_*,\sigma)\langle
v_*\rangle^{-N_1}\langle
v'\rangle^{-N_3}\\&&\times\mathcal{G}_*\mathcal{F}'
 \bigg(\langle v\rangle^{-N_2}\mathcal{H}_j-\langle v'\rangle^{-N_2}\mathcal{H}'_j\bigg)\bigg(\langle v'\rangle^{l}-\langle v\rangle^{l}
 \bigg)d\sigma.
 \eeno Then $\mathfrak{R}_2=\sum_j \mathfrak{R}_2^j $. As done before, we also introduce the angular cut-off
function $\phi_j$ and split $\mathfrak{R}_2^j$ into two parts
$\mathfrak{R}_{2,1}^j$ and $\mathfrak{R}_{2,2}^j$ which are defined
as \beno
\mathfrak{R}_{2,1}^j&=&\int_{\R^6}dvdv_*\int_{\SS^{2}}B(v-v_*,\sigma)\langle
v_*\rangle^{-N_1}\langle v'\rangle^{-N_3}\mathcal{G}_*\mathcal{F}'
 \bigg(\langle v\rangle^{-N_2}\mathcal{H}_j-\langle v'\rangle^{-N_2}\mathcal{H}'_j\bigg)\\&&\quad\times
 \phi_j(|v-v'|)\bigg(\langle v'\rangle^{l}-\langle v\rangle^{l}
 \bigg)d\sigma
 \eeno
and \beno
\mathfrak{R}_{2,2}^j&=&\int_{\R^6}dvdv_*\int_{\SS^{2}}B(v-v_*,\sigma)\langle
v_*\rangle^{-N_1}\langle v'\rangle^{-N_3}\mathcal{G}_*\mathcal{F}'
 \bigg(\langle v\rangle^{-N_2}\mathcal{H}_j-\langle v'\rangle^{-N_2}\mathcal{H}'_j\bigg)\\&&\quad\times
 [1-\phi_j(|v-v'|)]\bigg(\langle v'\rangle^{l}-\langle v\rangle^{l}
 \bigg)d\sigma.
 \eeno

We first treat with $\mathfrak{R}_{2,1}^j$. One may check that \beno
\mathfrak{R}_{2,1}^j&\lesssim&
\int_{\R^6}dvdv_*\int_{\SS^{2}}B(v-v_*,\sigma)\langle
v_*\rangle^{-N_1}\langle
v'\rangle^{-N_3}\mathcal{G}_*\mathcal{F}'|v-v'|^2\phi_j(|v-v'|)\\&&\quad\times
 \sum_{i=0}^1| \na_v^i \mathcal{H}_j(\gamma(\kappa_1))| \langle \gamma(\kappa_1)\rangle^{-N_2} \langle \gamma(\kappa_2)\rangle^{l-1} d\sigma.  \eeno
where $\kappa_1, \kappa_2\in [0,1]$. Noticing that \ben\label{wei1}
\langle v_*\rangle^{-N_1}\langle v'\rangle^{-N_3} \langle
\gamma(\kappa_1)\rangle^{-N_2} \langle
\gamma(\kappa_2)\rangle^{l-1}\lesssim \langle v_*-
\gamma(\kappa_1)\rangle^{-\tilde{N}_1+l-1},\een  we deduce that
\beno \mathfrak{R}_{2,1}^j&\lesssim&
\int_{\R^6}dvdv_*\int_{\SS^{2}}B(v-v_*,\sigma)\langle v_*-
\gamma(\kappa_1)\rangle^{-\tilde{N}_1+l-1}\mathcal{G}_*\mathcal{F}'|v-v'|^2\\&&\quad\times\phi_j(|v-v'|)
 \sum_{i=0}^1| \na_v^i \mathcal{H}_j(\gamma(\kappa_1))|   d\sigma.  \eeno
By Cauchy-Schwartz inequality,  $\mathfrak{R}_{2,1}^j$ can be
controlled by the quantities  \beno VI_1&\eqdefa&
\int_{\R^6}dvdv_*\int_{\SS^{2}}B(v-v_*,\sigma)\langle v_*-
\gamma(\kappa_1)\rangle^{-\tilde{N}_1+l-1}\mathcal{G}_*|v-v'|^2\\&&\quad\times
 \phi_j(|v-v'|)\sum_{i=0}^1| \na_v^i \mathcal{H}_j(\gamma(\kappa_1))|^2  d\sigma   \eeno
 and \beno VI_2&\eqdefa&
\int_{\R^6}dvdv_*\int_{\SS^{2}}B(v-v_*,\sigma)\langle v_*-
\gamma(\kappa_1)\rangle^{-\tilde{N}_1+l-1}\mathcal{G}_*|v-v'|^2\\&&\quad\times\phi_j(|v-v'|)
   |\mathcal{F}'|^2 d\sigma.   \eeno
We only need to show the estimate for $VI_1$.   One has \beno
VI_1&\lesssim& \int_{\R^6}dvdv_*\langle
v_*-\gamma(\kappa_1)\rangle^{-\tilde{N}_1+l-1}
\mathcal{G}_*|v_*-\gamma(\kappa_1)|^{\gamma+2}\sum_{i=0}^1| \na_v^i
\mathcal{H}_j(\gamma(\kappa_1))|^2\\&&\quad\times \int_{\SS^{2}}
 \phi_j(|v-v'|)b( \langle\frac{v-v_*}{|v-v_*|},
\sigma\rangle)\f{1-\langle\frac{v-v_*}{|v-v_*|}, \sigma\rangle}{2}
d\sigma.   \eeno Fixed $\sigma$ and $v_*$ and noting the fact
\eqref{ja}, we change of variables from $ v\rightarrow
u=\gamma(\kappa_1)$ and then it gives \beno VI_1&\lesssim&
\int_{\R^6}dudv_*\langle v_*-u\rangle^{-\tilde{N}_1+l-1}
\mathcal{G}_*|v_*-u|^{\gamma+2}\sum_{i=0}^1| \na_v^i
\mathcal{H}_j(u)|^2\\&&\quad\times \int_0^{c2^{-j}|u-v_*|^{-1}}
\theta^{1-2s} d\theta.\eeno

Since $\tilde{N}_1\ge l-1+\gamma+2s$, then for the case of
$\gamma+2s>0$, there holds \beno VI_1&\lesssim&
 \|\mathcal{G}\|_{L^1(\R^3_v)}\|\mathcal{H}_j\|_{H^s(\R^3_v)}^2.\eeno
While for the case of $\gamma+2s\le 0$, there holds \beno
VI_1&\lesssim&
 (\|\mathcal{G}\|_{L^1(\R^3_v)}+\|\mathcal{G}\|_{L^\f32(\R^3_v)})\|\mathcal{H}_j\|_{H^s(\R^3_v)}^2.\eeno

Similar calculation can be applied to $VI_2$ and it gives that in
the case of $\gamma+2s>0$, there holds \beno VI_2&\lesssim&
2^{-2j(1-s)}\|\mathcal{G}\|_{L^1(\R^3_v)}\|\mathcal{F}\|_{L^2(\R^3_v)}^2.\eeno
While for the case of $\gamma+2s\le 0$, there holds \beno
VI_2&\lesssim&
2^{-2j(1-s)}(\|\mathcal{G}\|_{L^1(\R^3_v)}+\|\mathcal{G}\|_{L^\f32(\R^3_v)})\|\mathcal{F}\|_{L^2(\R^3_v)}^2.\eeno

From which, we deduce that  in the case of $\gamma+2s>0$, there
holds \ben\label{com1} \sum_{j} \mathfrak{R}_{2,1}^j&\lesssim&
 \|\mathcal{G}\|_{L^1(\R^3_v)}\|\mathcal{H}\|_{H^s(\R^3_v)}\|\mathcal{F}\|_{L^2(\R^3_v)}.\een
While for the case of $\gamma+2s\le 0$, there holds \ben\label{com2}
\sum_{j} \mathfrak{R}_{2,1}^j &\lesssim&
(\|\mathcal{G}\|_{L^1(\R^3_v)}+\|\mathcal{G}\|_{L^\f32(\R^3_v)})\|\mathcal{H}\|_{H^s(\R^3_v)}\|\mathcal{F}\|_{L^2(\R^3_v)}.
\een

Now we treat with the term $\mathfrak{R}_{2,2}^j$. We first observe
that  the angular function $b(\theta)$ now is locally integrable
 thanks to the cut-off function $1-\phi_j(|v-v'|)$. Then the bound
 for $\mathfrak{R}_{2,2}^j$ can be reduced to the estimation of
\beno
VII&\eqdefa&\int_{\R^6}dvdv_*\int_{\SS^{2}}B(v-v_*,\sigma)\langle
v_*\rangle^{-N_1}\langle v'\rangle^{-N_3}\mathcal{G}_*\mathcal{F}'
  \langle v\rangle^{-N_2}|\mathcal{H}_j| \\&&\quad\times
 [1-\phi_j(|v-v'|)]\bigg(\langle v'\rangle^{l}-\langle v\rangle^{l}
 \bigg)d\sigma.
 \eeno

It is easy to check that \beno
VII&\lesssim&\int_{\R^6}dvdv_*\int_{\SS^{2}}B(v-v_*,\sigma)\langle
v_*\rangle^{-N_1}\langle v'\rangle^{-N_3}\mathcal{G}_*\mathcal{F}'
  \langle v\rangle^{-N_2}|\mathcal{H}_j||v-v'| \\&&\quad\times
 [1-\phi_j(|v-v'|)] \langle \gamma(\kappa_3)\rangle^{l-1} d\sigma,
 \eeno
where $\kappa_3\in [0,1]$. Thanks to \eqref{wei1}, by
Cauchy-Schwartz inequality, $VII$ can be controlled by $VII_1$ and
$VII_2$ defined as \beno
VII_1&\eqdefa&\int_{\R^6}dvdv_*\int_{\SS^{2}}B(v-v_*,\sigma)\langle
v_*-v\rangle^{-\tilde{N}_1+l-1}
\mathcal{G}_*\\&&\quad\times|\mathcal{H}_j |^2|v-v'|
 [1-\phi_j(|v-v'|)] d\sigma  \eeno
 and \beno
VII_2&\eqdefa&\int_{\R^6}dvdv_*\int_{\SS^{2}}B(v-v_*,\sigma)\langle
v_*-v'\rangle^{-\tilde{N}_1+l-1}
\mathcal{G}_*\\&&\quad\times|\mathcal{F}' |^2|v-v'|
 [1-\phi_j(|v-v'|)] d\sigma.  \eeno
 We only need to give the estimate to one of them. As for $VII_2$,
 by change of variables from $v$ to $v'$, one may has
\beno VII_2&\lesssim&\int_{\R^6}dv'dv_*\langle
v_*-v'\rangle^{-\tilde{N}_1+l-1}|v_*-v'|^{\gamma+1}
\mathcal{G}_*|\mathcal{F}' |^2\mathbf{1}_{c2^{-j}|v_*-v'|^{-1}<
 \frac{\pi}{4}}\\&&\quad\times\int_{ \sqrt{1-\langle
\frac{v-v_*}{|v-v_*|},\sigma \rangle }\ge c2^{-j}|v_*-v'|^{-1}}
 b(\langle
\frac{v-v_*}{|v-v_*|},\sigma \rangle )\sqrt{1-\langle
\frac{v-v_*}{|v-v_*|},\sigma \rangle } d\sigma,  \eeno where we use
the fact \eqref{ks1}. From which, we deduce that \beno
VII_2&\lesssim& 2^{-j}2^{2js}\int_{\R^6}dv'dv_*\langle
v_*-v'\rangle^{-\tilde{N}_1+l-1}|v_*-v'|^{\gamma+2s}
\mathcal{G}_*|\mathcal{F}' |^2, \eeno  which implies that in the
case of $\gamma+2s>0$, there holds \beno VII_2&\lesssim&
2^{-j(1-s)}2^{js}\|\mathcal{G}\|_{L^1(\R^3_v)}\|\mathcal{F}\|_{L^2(\R^3_v)}^2.\eeno
While for the case of $\gamma+2s\le 0$, there holds \beno
VII_2&\lesssim&
2^{-j(1-s)}2^{js}(\|\mathcal{G}\|_{L^1(\R^3_v)}+\|\mathcal{G}\|_{L^\f32(\R^3_v)})\|\mathcal{F}\|_{L^2(\R^3_v)}^2.\eeno

Similarly, one can obtain that in the case of $\gamma+2s>0$, there
holds \beno VII_1&\lesssim&
2^{-j(1-s)}2^{js}\|\mathcal{G}\|_{L^1(\R^3_v)}\|\mathcal{H}_j\|_{L^2(\R^3_v)}^2.\eeno
While for the case of $\gamma+2s\le 0$, there holds \beno
VII_1&\lesssim&
2^{-j(1-s)}2^{js}(\|\mathcal{G}\|_{L^1(\R^3_v)}+\|\mathcal{G}\|_{L^\f32(\R^3_v)})\|\mathcal{H}_j\|_{L^2(\R^3_v)}^2.\eeno

Thanks to the Bernstein's inequality, we arrive at in the case of
$\gamma+2s>0$ and $\tilde{N}_1\ge l-1+\gamma+2s $, there holds \beno
VII&\lesssim&
2^{-j(1-s)}\|\mathcal{G}\|_{L^1(\R^3_v)}\|\mathcal{F}\|_{L^2(\R^3_v)}\|\mathcal{H}_j\|_{H^s(\R^3_v)}.\eeno
While for the case of $\gamma+2s\le 0$, there holds \beno
VII&\lesssim&
2^{-j(1-s)}(\|\mathcal{G}\|_{L^1(\R^3_v)}+\|\mathcal{G}\|_{L^\f32(\R^3_v)})
\|\mathcal{F}\|_{L^2(\R^3_v)}\|\mathcal{H}_j\|_{H^s(\R^3_v)}.\eeno

Finally,  we   obtain that in the case of $\gamma+2s>0$, there holds
\beno \sum_{j} \mathfrak{R}_{2,2}^j&\lesssim&
 \|\mathcal{G}\|_{L^1(\R^3_v)}\|\mathcal{H}\|_{H^s(\R^3_v)}\|\mathcal{F}\|_{L^2(\R^3_v)}.\eeno
While for the case of $\gamma+2s\le 0$, there holds \beno \sum_{j}
\mathfrak{R}_{2,2}^j &\lesssim&
(\|\mathcal{G}\|_{L^1(\R^3_v)}+\|\mathcal{G}\|_{L^\f32(\R^3_v)})\|\mathcal{H}\|_{H^s(\R^3_v)}\|\mathcal{F}\|_{L^2(\R^3_v)}.
\eeno

From which together with \eqref{com1} and \eqref{com2}, we complete
the proof of the lemma. \ef

Thanks to  the Lemma \ref{le4} and Lemma \ref{le5}, we complete the
proof to the corollary. \ef

\section{Coercivity estimate for the collision operator}
In this section, we shall give the coercivity bound for  the
collision operator. Our main motivation comes from the sub-elliptic
estimate \eqref{co1} entailed by entropy dissipation and the upper
bound estimate \eqref{u1} and \eqref{u2} for the collision operator.
Roughly, our strategy is carried out as follows. We first use the
trick   to reformulate the functional $\langle -Q(g,f), f\rangle_v$
as $\langle -\mathcal {Q}(\mathcal {G},\mathcal {F}),\mathcal
{F}\rangle_v$ by introducing $\mathcal {G}=g\langle v\rangle^N$ and
$\mathcal {F}=f\langle v\rangle^N$. As a consequence, the kinetic
part $\Phi(|v-v_*|)$ will be cancelled by the additional   factor
$\langle v_*\rangle^{-N}$ and $\langle v\rangle^{-N}$ which means
the cases of the hard potential and soft potential can be reduced to
the Maxwellian case but at the price of occurring the lower order
terms. We point out that here  the lower order term means the lower
derivative term or the  term with lower weight. Thanks to the
estimates  \eqref{co1}, \eqref{u1} and \eqref{u2}, we finally obtain
the following theorem:

\begin{thm}\label{coc} Let the collision kernel $B(|v-v_*|, \sigma)$
satisfies the Assumption A. Suppose the function $g$ satisfies \beno
\|g\|_{L^1_2(\R^3_v)}+\|g\|_{L\log L( \R^3_v)}<\infty. \eeno Then
there exists a constant $C_g$ depending on $\|g\|_{L^1_1(\R^3_v)}$
and $\|g\|_{L\log L( \R^3_v)}$ such that in the case of
$\gamma+2s>0$, there holds \ben\label{coc1}  \langle -Q(g,f),
f\rangle_v &\gtrsim& C_g \|f\langle
v\rangle^{\frac{\gamma}{2}}\|_{H^s(\R^3_v)}^2- (C_g^{1-2s}\|g\langle
 v\rangle^{\tilde{\gamma}}\|_{L^1(\R^3_v)}^{2s}+
 \|g\langle
 v\rangle^{\tilde{\gamma}}\|_{L^1(\R^3_v)})\nonumber\\&&\quad\times  \|f\langle
v\rangle^{\frac{\gamma}{2}} \|_{L^2(\R^3_v)}^2-\|g\langle
v\rangle^{|\gamma|}\|_{L^1(\R^3_v)}
   \|f\langle v\rangle^\frac{\gamma}{2}\|_{H^{\varrho}(\R^3_v)}^2 ,\een  where
$\varrho<s$ and \beno \tilde{\gamma}=|\gamma+2|\mathbf{1}_{\gamma\le
0}+|\gamma-2|\mathbf{1}_{\gamma>0};\eeno and in the case of
$\gamma+2s\le 0$, there holds \ben\label{coc2} \langle -Q(g,f),
f\rangle_v &\gtrsim& C_g \|f\langle v\rangle^{\frac{\gamma}{2}}
\|_{H^s(\R^3_v)}^2- (C_g^{1-2s}\|g\langle
 v\rangle^{|\gamma+2|}\|_{L^1(\R^3_v)}^{2s}+
 \|g\langle
 v\rangle^{|\gamma+2|}\|_{L^1(\R^3_v)})\nonumber\\&&\, \times\|f\langle
v\rangle^{\frac{\gamma}{2}}\|_{L^2(\R^3_v)}^2-(\|g\langle
v\rangle^{|\gamma|}\|_{L^1(\R^3_v)}
 +\|g\langle
v\rangle^{|\gamma|}\|_{L^\f32(\R^3_v)})\|f\langle
v\rangle^\frac{\gamma}{2}\|_{H^{ \varrho}(\R^3_v)}^2
 ,\een with $\varrho<s$.
\end{thm}

\medskip

Before the proof, let us give some comments on the Theorem
\ref{coc}. First of all, the proof of theorem can be applied to
obtain the smoothing effect estimate \eqref{co2} which is entailed
by entropy dissipation. Secondly, comparing to the upper bound
estimates \eqref{u1} and \eqref{u2}, we lose  $2s$ order of weight
in $v$ in the coercivity estimate. The main reason may lie in the
fact that it is still not clear to the structure of the collision
operator.

\medskip

\noindent{\bf Proof of the Theorem:} It is easy to see that \beno
\langle  Q(g,f), f\rangle_v&=&\int_{\R^6}dv_*dv\int_{\SS^2}
B(|v-v_*|,\sigma)g_*f(f'-f)d\sigma \\
&=&\f12\int_{\R^6}dv_*dv\int_{\SS^2}
B(|v-v_*|,\sigma)g_*(f'^2-f^2)d\sigma\\&&\quad-\f12
\int_{\R^6}dv_*dv
\int_{\SS^2}B(|v-v_*|,\sigma)g_*(f'-f)^2d\sigma\\&=&
\mathfrak{I}_1-\mathfrak{I}_2.\eeno

\noindent{\bf Step 1: Upper bound for $\mathfrak{I}_1$. } By change
of variables, one has \beno \mathfrak{I}_1&=& |\SS^1|
\int_{\R^6}\int_0^\frac{\pi}{2}
\sin\theta\bigg(\f1{\cos^3\frac{\theta}{2}}B(\frac{|v-v_*|}{\cos\frac{\theta}{2}},\cos\theta)-B(|v-v_*|,\cos\theta)\bigg)g_*f^2d\theta
dv_*dv
\\&=&|\SS^1|\int_{\R^6}\int_0^\frac{\pi}{2}
\sin\theta b(\cos\theta)
|v-v_*|^\gamma\{(\cos^{-1}\f{\theta}2)^{\gamma+3}-1\}g_*f^2d\theta
dv_*dv. \eeno

Using the fact  that \beno x^m-y^m= m\int_y^x z^{m-1}dz,\eeno we
obtain that for $\theta\in[0,\frac{\pi}{2}],$ \beno
(\cos^{-1}\f{\theta}2)^{\gamma+3}-1\lesssim
\sin^2\frac{\theta}{2},\eeno which immediately implies \beno
\mathfrak{I}_1= C\int_{\R^6} |v-v_*|^{\gamma}g_*f^2 dvdv_*.\eeno
In the forthcoming argument, we shall give the different upper
bounds for $\mathfrak{I}_1$ with respect to the value $\gamma$.

 {\bf
Case 1: $\gamma>0$.} It is easy to check \beno
\mathfrak{I}_1\lesssim
\|g\|_{L^1_\gamma(\R^3_v)}\|f\|_{L^2_{\frac{\gamma}{2}}(\R^3_v)}^2
\eeno

{\bf Case 2: $\gamma<0$.} We may rewrite $\mathfrak{I}_1 $ as \beno
\mathfrak{I}_1= C\int_{\R^6} |v-v_*|^{\gamma}\langle
v\rangle^{-\gamma}\langle v_*\rangle^\gamma\mathcal {G}_*\mathcal
{F}^2 dvdv_*,\eeno where $\mathcal {G}_*=g_*\langle
v\rangle^{-\gamma}, \mathcal {F}=f\langle
v\rangle^{\frac{\gamma}{2}}$. It is easy to check that

\beno  \mathfrak{I}_1\lesssim \int_{\R^6}
(1+|v-v_*|^{\gamma}\mathbf{1}_{|v-v_*|\le 1})\mathcal {G}_*\mathcal
{F}^2 dvdv_*.\eeno

Now we first treat the case
 $\gamma+2s>0$. That is, there exist a positive $\delta$ such that $\gamma=-(2s-\delta)$ which implies that
 $\chi(v)=\mathbf{1}_{|v|\le 1}|v|^{\gamma}\in L^{\frac{3-\epsilon}{2s-\delta}}$. Then by H\"{o}lder's inequality, one may
obtain that \beno  \mathfrak{I}_1 &\lesssim& \|\mathcal
{G}\|_{L^1(\R^3_v)}(\|\mathcal {F}\|_{L^2(\R^3_v)}^2+\|\mathcal
{F}^2\|_{L^q(\R^3_v)}\|\chi\|_{L^p(\R^3_v)}), \eeno where
$\f1p+\f1q=1$ and $p=\frac{3-\epsilon}{2s-\delta}$. By standard
Sobolev's embedding theorem, one has \beno \|\mathcal
{F}^2\|_{L^q(\R^3_v)}\le \|\mathcal {F}\|_{H^{\varrho}(\R^3_v)}^2,
\eeno where $\varrho=\frac{6s-3\delta}{6-2\epsilon}$. Choose
$\epsilon <\frac{3\delta}{2s}$, then one may obtain that
$\varrho<s$. And we arrive at \beno  \mathfrak{I}_1&\lesssim&
\|g\langle v\rangle^{-\gamma} \|_{L^1(\R^3_v)}\|f\langle
v\rangle^{\frac{\gamma}{2}} \|_{H^{\varrho}(\R^3_v)}^2.\eeno

Secondly, we handle with the case  $-1<\gamma+2s<0$. Let
$\delta=1+\gamma+2s$. Then by the Assumption A, it gives that
$0<\delta<1$ and $\chi(v)\in L^{\frac{3-\epsilon }{1+2s-\delta}}$.
Choose $\epsilon<\f{3\delta}{1+2s}$. Then we can deduce that there
exists a constant $\varrho=\f{6s-3\delta+\epsilon}{2(3-\epsilon)}$
verifying $\varrho<s$ such that \beno
\f{1+2\varrho}{3}=\f{1+2s-\delta}{3-\epsilon}.\eeno   Then by
H\"{o}lder's inequality and Young's inequality, one has that \beno
\mathfrak{I}_1 &\lesssim& \|\mathcal {G}\|_{L^1(\R^3_v)}\|\mathcal
{F}\|_{L^2(\R^3_v)}^2+\|\mathcal {F}^2\|_{L^q(\R^3_v)}\|\mathcal
{G}*\chi\|_{L^p(\R^3_v)}), \eeno where $\f1p+\f1q=1$ and
$\f12=\f1{2q}+\f{ \varrho}3 $. By standard Sobolev's embedding
theorem, we obtain that  \beno \mathfrak{I}_1 &\lesssim& (\|\mathcal
{G}\|_{L^1(\R^3_v)}+\|\mathcal {G}\|_{L^\f32(\R^3_v)})\|\mathcal
{F}\|_{H^{\varrho}(\R^3_v)}^2.\eeno

Putting together the estimates for $\mathfrak{I}_1 $,  we arrive at
there exists a constant $\varrho$ verifying $\varrho<s$ such that
for the case of $\gamma+2s>0$, there holds \ben\label{eh10}
\mathfrak{I}_1\lesssim
 \|g\langle
v\rangle^{|\gamma|}\|_{L^1(\R^3_v)}
   \|f\langle v\rangle^\frac{\gamma}{2}\|_{H^{\varrho}(\R^3_v)}^2, \een
   and for
the case of  $\gamma+2s<0$, there holds \ben\label{eh8}
\mathfrak{I}_1\lesssim(\|g\langle
v\rangle^{|\gamma|}\|_{L^1(\R^3_v)}
 +\|g\langle
v\rangle^{|\gamma|}\|_{L^\f32(\R^3_v)})\|f\langle
v\rangle^\frac{\gamma}{2}\|_{H^{\varrho}(\R^3_v)}^2. \een

\medskip
\noindent{\bf Step 2: Lower bound for $\mathfrak{I}_2$. }  Recalling
that \beno 2\mathfrak{I}_2=\int_{\R^6}\int_{\SS^2}
 |v-v_*|^\gamma b(\cos\theta) g_*(f'-f)^2d\sigma dv_*dv,\eeno and setting $\mathcal{F}=f\langle
 v\rangle^\frac{\gamma}{2}$, we can rewrite $\mathfrak{I}_2$ as \beno  2\mathfrak{I}_2=\int_{\R^6}\int_{\SS^2}
 |v-v_*|^\gamma b(\cos\theta) g_*(\langle
 v'\rangle^{-\frac{\gamma}{2}}\mathcal{ F}'-\langle
 v\rangle^{-\frac{\gamma}{2}}\mathcal{F})^2d\sigma dv_*dv.\eeno
We shall   give the different lower bound for $\mathfrak{I}_2$ with
respect to the value $\gamma$.

{\bf Case 1:   $\gamma\le 0$}. Thanks to the fact $(A-B)^2\ge
\frac{A^2}{2}-B^2, $ one has \beno && (\langle
 v'\rangle^{-\frac{\gamma}{2}}\mathcal{ F}'-\langle
 v\rangle^{-\frac{\gamma}{2}}\mathcal{F})^2\\&&\quad=\bigg(\langle
 v'\rangle^{-\frac{\gamma}{2}}(\mathcal{ F}'- \mathcal{F} )+\mathcal{ F}( \langle
 v'\rangle^{-\frac{\gamma}{2}} -\langle
 v\rangle^{-\frac{\gamma}{2}})\bigg)^2\\&&\quad\ge\f12  \langle
 v'\rangle^{-\gamma}(\mathcal{ F}'- \mathcal{F} )^2-\mathcal{ F}^2( \langle
 v'\rangle^{-\frac{\gamma}{2}} -\langle
 v\rangle^{-\frac{\gamma}{2}})^2 .\eeno
Then one may obtain that \beno  2\mathfrak{I}_2&\ge& \f12
\int_{\R^6}\int_{\SS^2}
 |v-v_*|^\gamma b(\cos\theta) g_* \langle
 v'\rangle^{-\gamma}(\mathcal{ F}'- \mathcal{F} )^2d\sigma dv_*dv\\&&\quad- \int_{\R^6}\int_{\SS^2}
 |v-v_*|^\gamma b(\cos\theta) g_*\mathcal{ F}^2( \langle
 v'\rangle^{-\frac{\gamma}{2}} -\langle
 v\rangle^{-\frac{\gamma}{2}})^2  d\sigma dv_*dv\\&\eqdefa& \mathfrak{L}_1-\mathfrak{L}_2.\eeno

Noting that \beno  |v-v_*|^\gamma \langle
 v'\rangle^{-\gamma}\langle
 v_*\rangle^{-\gamma}\gtrsim 1,\eeno one has
 \beno \mathfrak{L}_1\gtrsim \int_{\R^6}
 b(\cos\theta) (g_* \langle
 v_*\rangle^{\gamma})(\mathcal{ F}'- \mathcal{F} )^2d\sigma dv_*dv.\eeno
Due to the well-known entropy dissipation inequality, we arrive at
\beno   \mathfrak{L}_1 \gtrsim C_g\|\mathcal{F}\|_{H^s(\R^3_v)}^2,
\eeno where $C_g$ depends on $\|g\langle
 v\rangle^{\gamma}\|_{L^1_1(\R^3_v)}$, $\|g\langle
 v\rangle^{\gamma}\|_{L\log L(\R^3_v)}$ and $b$.
Next, we turn to the estimate of $\mathfrak{L}_2$. One may have
\beno \mathfrak{L}_2\lesssim \int_{\R^6}\int_{\SS^2}
 |v-v_*|^\gamma b(\cos\theta) g_*\mathcal{ F}^2|v-v'|^2  \langle
 \gamma(\kappa)\rangle^{-\gamma-2}  d\sigma dv_*dv \eeno with $\kappa\in
 [0,1]$.
Noting the fact \beno |v-v'|=|v-v_*|\sin \frac{\theta}{2} ,\eeno we
deduce that  in the case of $\gamma+2\ge 0$, there holds \beno
\mathfrak{L}_2 &\lesssim& \int_{\R^6}  g_*\mathcal{ F}^2   \langle
 v_*\rangle^{\gamma+2}   dv_*dv\\&\lesssim& \|g\langle
 v\rangle^{|\gamma+2|}\|_{L^1(\R^3_v)} \|\mathcal{ F}  \|_{L^2(\R^3_v)}^2.\eeno

From which, we obtain that for $\gamma+2\ge0$, there holds
\ben\label{eh11} \mathfrak{I}_2 \gtrsim
\frac{C_g}{2}\|\mathcal{F}\|_{H^s(\R^3_v)}^2-  \|g\langle
 v\rangle^{|\gamma+2|}\|_{L^1(\R^3_v)}  \|\mathcal{ F}
 \|_{L^2(\R^3_v)}^2.\een

 While in the case of $\gamma+2<0$, thanks to the Assumption A, one has \beno
 \gamma+2>-1, s>\f12. \eeno Then for any $\epsilon$, there holds
 \beno \mathfrak{L}_2 &\lesssim& \int_{\R^6} (1+|v-v_*|^{\gamma+2}\mathbf{1}_{|v-v_*|\le1}) g_*\mathcal{ F}^2   \langle
 v_*\rangle^{|\gamma+2|}   dv_*dv\\&\lesssim&  \|g\langle
 v\rangle^{|\gamma+2|}\|_{L^1(\R^3_v)}(\|\mathcal{ F}\|_{L^2(\R^3_v)}^2+\|\mathcal{ F}\|_{\dot{H}^{\f12}(\R^3_v)}^2)\\&\lesssim&  \|g\langle
 v\rangle^{|\gamma+2|}\|_{L^1(\R^3_v)}[(1+\epsilon^{-(2s-1)})\|\mathcal{ F}\|_{L^2(\R^3_v)}^2 +
  \epsilon\|\mathcal{ F}\|_{\dot{H}^{s}(\R^3_v)}^2]. \eeno
 From which, we obtain that for $\gamma+2<0$,
there holds \ben\label{eh12} \mathfrak{I}_2 \gtrsim
\frac{C_g}{2}\|\mathcal{F}\|_{H^s(\R^3_v)}^2-C_g^{1-2s} (\|g\langle
 v\rangle^{|\gamma+2|}\|_{L^1(\R^3_v)}+C_g)^{2s} \|\mathcal{ F}
 \|_{L^2(\R^3_v)}^2.\een

Thanks to \eqref{eh11} and \eqref{eh12}, we conclude that for
$\gamma<0$, $\mathfrak{I}_2$ can be estimated as \ben\label{eh13}
\mathfrak{I}_2 &\gtrsim&
\frac{C_g}{2}\|\mathcal{F}\|_{H^s(\R^3_v)}^2-(C_g^{1-2s}\|g\langle
 v\rangle^{|\gamma+2|}\|_{L^1(\R^3_v)}^{2s}+
 \|g\langle
 v\rangle^{|\gamma+2|}\|_{L^1(\R^3_v)})\nonumber\\&&\times
 \|\mathcal{ F}
 \|_{L^2(\R^3_v)}^2.\een

{\bf  Case 2:  $\gamma>0$}. Observing the fact that \beno \langle v
\rangle\sim |v|+\mathbf{1}_{|v|\le 1}, \eeno we may obtain that
\beno \mathfrak{I}_2&\gtrsim&  \int_{\R^6}\int_{\SS^2}
 \langle v-v_*\rangle^\gamma b(\cos\theta) g_*(f'-f)^2d\sigma dv_*dv\\&&\quad- \int_{\R^6}\int_{\SS^2}
 \mathbf{1}_{|v-v_*|\le 1} b(\cos\theta) g_*(f'-f)^2d\sigma dv_*dv. \eeno
 Noting \beno \langle v-v_*\rangle^\gamma\geq\langle v'-v_*\rangle^\gamma\ge \langle v_*\rangle^{-\gamma}\langle v'\rangle^\gamma,  \eeno
and following the similar trick as the one in the case of
$\gamma<0$,
 we arrive at
  \beno \mathfrak{I}_2&\gtrsim& \f12
\int_{\R^6}\int_{\SS^2}
  b(\cos\theta) (g_* \langle
 v_*\rangle^{-\gamma})(\mathcal{ F}'- \mathcal{F} )^2d\sigma dv_*dv\\&&\quad- \int_{\R^6}\int_{\SS^2}
  b(\cos\theta)g_* \langle
 v_*\rangle^{-\gamma}f^2( \langle
 v'\rangle^{\frac{\gamma}{2}} -\langle
 v\rangle^{\frac{\gamma}{2}})^2  d\sigma dv_*dv  \\&&\quad- \int_{\R^6}\int_{\SS^2}
 \mathbf{1}_{|v-v_*|\le 1} b(\cos\theta) g_*(f'-f)^2d\sigma dv_*dv\\&\eqdefa&\mathfrak{L}_3-\mathfrak{L}_4-\mathfrak{L}_5.\eeno
Here we also set $\mathcal{F}=f\langle
 v\rangle^\frac{\gamma}{2}$.

It is easy to check that   \beno \mathfrak{L}_3 \gtrsim
C_g\|\mathcal{F}\|_{H^s(\R^3_v)}^2,  \eeno where $C_g$ depends on
$\|g\langle
 v\rangle^{-\gamma}\|_{L^1_1(\R^3_v)}$, $\|g\langle
 v\rangle^{-\gamma}\|_{L\log L(\R^3_v)}$ and $b$. As for the term
 $\mathcal{L}_4$, one has \beno \mathfrak{L}_4\lesssim \int_{\R^6}\int_{\SS^2}
   b(\cos\theta) g_*f^2|v-v'|^2 \langle v_*\rangle^{-\gamma}\langle
 \gamma(\kappa)\rangle^{\gamma-2}  d\sigma dv_*dv \eeno with $\kappa\in
 [0,1]$.
From which, we deduce that \beno \mathfrak{L}_4\lesssim  \|g\langle
 v\rangle^{|\gamma-2|}\|_{L^1(\R^3_v)}  \|\mathcal{ F}  \|_{L^2(\R^3_v)}^2. \eeno
As for the term $\mathcal{L}_5$, noting that \beno &&\int_{\R^6}
 \mathbf{1}_{|v-v_*|\le 1} b(\cos\theta) g_*(f'-f)^2d\sigma dv_*dv\\&&\quad=-2\int_{\R^6}
 \mathbf{1}_{|v-v_*|\le 1} b(\cos\theta) g_*f(f'-f)d\sigma dv_*dv\\&&\qquad+\int_{\R^6}
 \mathbf{1}_{|v-v_*|\le 1} b(\cos\theta) g_*(f'^2-f^2)d\sigma dv_*dv,\eeno
we conclude that by the proof of  the Theorem \ref{ub} and the
result of the step 1, there holds \beno \mathfrak{L}_5\lesssim  \|g
\|_{L^1(\R^3_v)}( \|f \|_{L^2(\R^3_v)}^2+ \|f
\|_{H^s(\R^3_v)}^2).\eeno From which, we deduce that \beno
\mathfrak{I}_2\gtrsim C_g\|\mathcal{F}\|_{H^s(\R^3_v)}^2-\|g \langle
 v\rangle^{|\gamma-2|} \|_{L^1 (\R^3_v)}\|\mathcal{F} \|^2_{L^2(\R^3_v)}- \|g \|_{L^1(\R^3_v)}\|f \|_{\dot{H}^s(\R^3_v)}^2.\eeno
Noticing that for any smooth function $\chi_R$ defined as
$\chi_R=\chi(\frac{\cdot}{R})$ with $0\le \chi\le 1, \chi=1$ on
$B_1$ and $\mathrm{supp}\,(\chi)\subset B_2$, there holds \beno
\|f\|_{\dot{H}^s(\R^3_v)}\le \|f\chi_R\|_{\dot{H}^s( \R^3_v)}+
R^{-\f{\gamma}{2}}\|\mathcal{F}\|_{\dot{H}^s(\R^3_v)}+\|\mathcal{F}\|_{L^2(\R^3_v)}
\eeno (one may check the proof in the Appendix), which implies that
\beno \mathfrak{I}_2&\gtrsim& (C_g-R^{-\gamma}\|g
\|_{L^1(\R^3_v)})\|\mathcal{F}\|_{H^s(\R^3_v)}^2\\&&-\|g \langle
 v\rangle^{|\gamma-2|} \|_{L^1 (\R^3_v)}\|\mathcal{F} \|^2_{L^2(\R^3_v)}- \|g \|_{L^1(\R^3_v)}\|f\chi_R\|_{\dot{H}^s( \R^3_v)}^2.\eeno

Thanks to the entropy dissipation inequality, we can also  deduce
that for any $R>0$, there holds \beno   \mathfrak{I}_2+ \|g
\|_{L^1(\R^3_v)}\|f^2\|_{L^1(\R^3_v)}\gtrsim
C_{g,R}\|f\chi_R\|_{\dot{H}^s( \R^3_v)}^2,  \eeno where $C_{g,R}$
depends on $\|g \|_{L^1(\R^3_v)}$, $\|g \|_{L\log L(\R^3_v)}$, $b$
and $R$. Choose $R_1$ such that \beno  \|g \|_{L^1(\R^3_v)}\le
\f12C_gR_1^{\gamma}, \eeno then we finally obtain that for
$\gamma>0$, there holds \ben\label{eh14} \mathfrak{I}_2&\gtrsim&
\frac{C_gC_{g,R_1}}{2(\|g \|_{L^1 (\R^3_v)}+C_{g,R_1})}
 \|\mathcal{F}\|_{H^s(\R^3_v)}^2\nonumber\\&&-\|g \langle
 v\rangle^{|\gamma-2|} \|_{L^1 (\R^3_v)}\|\mathcal{F} \|^2_{L^2(\R^3_v)}.\een

Now we can conclude that \eqref{eh10}, \eqref{eh13} and \eqref{eh14}
imply the coercivity estimate for the case of $\gamma+2s>0$. And
\eqref{eh8} and \eqref{eh13} imply the coercivity estimate for the
case of $\gamma+2s\le0$ which completes the proof to the Theorem
\ref{coc}. \ef

As a direct application, we obtain the entropy dissipation estimate:
\begin{thm}  Let the collision kernel $B(|v-v_*|, \sigma)$
satisfies the Assumption A. Suppose the function $g$ satisfies \beno
\|g\|_{L^1_2(\R^3_v)}+\|g\|_{L\log L( \R^3_v)}<\infty. \eeno Then
there exists a constant $C_g$ depending on $\|g\|_{L^1_1(\R^3_v)}$
and $\|g\|_{L\log L( \R^3_v)}$  and a constant $\varrho<s$ such that
in the case of $\gamma+2s>0$ and $\gamma\le2$, there holds
\ben\label{coc3} &&\langle D(g,f), f\rangle_v +
(C_g^{1-2s}\|g\|_{L^1_2(\R^3_v)}^{2s}+
 \|g \|_{L^1_2(\R^3_v)}+C_g^{-\f{
\varrho}{s- \varrho}}\|g\|_{L^1_2(\R^3_v)}^{\f{ s}{s-
\varrho}})\|f\|_{L^1_2}\nonumber \\&\gtrsim& C_g \|\sqrt{f}\langle
v\rangle^{\frac{\gamma}{2}}\|_{H^s(\R^3_v)}^2.\een
\end{thm}

\noindent{\bf Proof:} Direct calculation gives that \beno D(g,f)
&\ge&- \int_{\R^6}dv_*dv\int_{\SS^2}
B(|v-v_*|,\sigma)g_*(f'-f)d\sigma \\
&&\quad+ \int_{\R^6}dv_*dv
\int_{\SS^2}B(|v-v_*|,\sigma)g_*(\sqrt{f'}-\sqrt{f})^2d\sigma. \eeno
We stress that the estimates for righthand side of the above
inequality are exactly as the same as the ones  for $\mathfrak{I}_1$
and $\mathfrak{I}_2$. Then we arrive at \ben\label{eh15}  D(g,f)
 &\gtrsim& C_g \|\sqrt{f}\langle
v\rangle^{\frac{\gamma}{2}}\|_{H^s(\R^3_v)}^2- (C_g^{1-2s}\|g\langle
 v\rangle^{\tilde{\gamma}}\|_{L^1(\R^3_v)}^{2s}+
 \|g\langle
 v\rangle^{\tilde{\gamma}}\|_{L^1(\R^3_v)})\nonumber\\&&\quad\times  \|\sqrt{f}\langle
v\rangle^{\frac{\gamma}{2}} \|_{L^2(\R^3_v)}^2-\|g\langle
v\rangle^{|\gamma|}\|_{L^1(\R^3_v)}
   \|\sqrt{f}\langle v\rangle^\frac{\gamma}{2}\|_{H^{\varrho}(\R^3_v)}^2 ,\een   where
$\varrho<s$ and \beno \tilde{\gamma}=|\gamma+2|\mathbf{1}_{\gamma\le
0}+|\gamma-2|\mathbf{1}_{\gamma>0}.\eeno Since now $\gamma+2s>0$ and
$\gamma\le2$, we deduce that $\tilde{\gamma}<2$ and $|\gamma| \le2$.
Thanks to the Young's inequality \beno \|f\|_{H^{\varrho}}\lesssim
\epsilon\|f\|_{H^s} +\epsilon^{-\f{ \varrho}{s- \varrho}}\|f\|_{L^2}
,\eeno we can rewrite \eqref{eh15} as

\beno  D(g,f)
 &\gtrsim& C_g \|\sqrt{f}\langle
v\rangle^{\frac{\gamma}{2}}\|_{H^s(\R^3_v)}^2- (C_g^{1-2s}\|g
\|_{L^1_2(\R^3_v)}^{2s}+
 \|g \|_{L^1_2(\R^3_v)})\|f\|_{L^1_2(\R^3_v)}\nonumber\\&&\quad
 -C_g^{-\f{
\varrho}{s- \varrho}}\|g\|_{L^1_2(\R^3_v)}^{\f{ s}{s-
\varrho}}\|f\|_{L^1_2(\R^3_v)} ,\eeno which completes the proof to
the Theorem. \ef

\section{Smoothing effect for the homogeneous Boltzmann equation}

In this section, we shall give the proof to the Theorem \ref{smh}.

\noindent{\bf Proof of the Theorem \ref{smh}:} We first assume that
infinite  $L^2$ moment estimate \eqref{am2} holds true for all the
collision kernel. The inductive argument will be applied  to prove
the smoothing effect of the homogeneous Boltzmann equation.

Let us assume that for some $m\in \N$ and all $l\in N$, there hold
\ben\label{in} \sup_{[t_0, \infty)}\|f\|_{H^m_l(\R^3_v)}<\infty.\een

Noting that \beno \pa_t\pa_v^{\alpha}f=\pa_v^{\alpha}
Q(f,f)=Q(f, \pa^\alpha_v
f)+\sum_{{\alpha_1+\alpha_2=\alpha}\atop{|\alpha_1|\ge 1}}C^{\alpha}_{\alpha_1}Q(\pa^{\alpha_1}_v f,
\pa^{\alpha_2}_v f).\eeno Set $g=\pa_v^\alpha f\langle v\rangle^l$,
then g solves \beno \pa_t g=Q(f, g)+[Q(f, \pa^{\alpha}_v f)\langle
v\rangle^l-Q(f, \pa^{\alpha}_v f\langle
v\rangle^l)]+Q\sum_{{\alpha_1+\alpha_2=\alpha}\atop{|\alpha_1|\ge 1}}C^{\alpha}_{\alpha_1}
(\pa^{\alpha_1}_v f, \pa^{\alpha_2}_v f)   \langle v\rangle^l.
\eeno

Suppose $|\alpha|=m+1$. Thanks to the  Theorem \ref{coc} and
standard interpolation inequality, one has \beno \langle Q(f, g)
,g\rangle_v \lesssim -\frac{C_f}{2}\|g\langle
v\rangle^{\f{\gamma}{2}}\|_{H^s}^2+(C_f^{1-2s}\|f\|_{L^2_{|\gamma|+4}}^{2s}+\|f\|_{L^2_{|\gamma|+4}})\|
f\|_{H^{m+\varrho+1}_{l+\frac{\gamma}{2} }}^2, \eeno
with $\varrho<s$.

Due to  Corollary 2.2, taking $N_2=l+\frac{\gamma}{2}$ and
$N_3=\frac{\gamma}{2}$ in the case of $ s<\f12$, otherwise taking
$N_2=l+\frac{\gamma}{2}$ and $N_3=\frac{\gamma}{2}+2s-1$, one may
arrive at \beno &&\langle Q(f, \pa^{\alpha}_v f)\langle
v\rangle^l-Q(f, \pa^{\alpha}_v f\langle v\rangle^l) ,g\rangle_v\\
&\lesssim &  \|f\|_{L^2_{|\gamma|+2s+2l+4}}(\|g\langle
v\rangle^{\f{\gamma}{2}}\|_{H^{\varrho}}^2+\|g\langle
v\rangle^{\f{\gamma}{2}}\|_{H^s}
\|f\|_{H^{m+1}_{l+\frac{\gamma}{2}+2s}}^2). \eeno

Applying the Theorem \ref{ub} with $N_2=2s+\f{\gamma}{2}+l$ and
$N_3=-l+\f{\gamma}{2}$, one may have \beno \langle
Q(\pa^{\alpha_1}_v f, \pa^{\alpha_2}_v f) , \langle
v\rangle^lg\rangle_v \lesssim \| \pa^{\alpha_1}_v
f\|_{L^2_{|\gamma|+2s+3l+4}}\|  \pa^{\alpha_2}_v
f\|_{H^s_{2s+\f{\gamma}{2}+l }}\|g\langle
v\rangle^{\f{\gamma}{2}}\|_{H^s}. \eeno

Patching together all the above estimates, by standard energy
estimate, we easily deduce that \beno
\f{d}{dt}\|g\|_{L^2}^2+\frac{C_f}{3}\|g\langle
v\rangle^{\f{\gamma}{2}}\|_{H^s}^2\lesssim  \|
f\|_{H^{m+1}_{l+\frac{\gamma}{2}+2s}}^2+\|
f\|_{H^{m+\varrho+1}_{l+\frac{\gamma}{2} }}^2,\eeno where we use the
Young's inequality and \eqref{in}.

From which, one obtains that  \ben\label{in1}
\f{d}{dt}\|f\|_{H^{m+1}_l}^2+
\|f\|_{H^{m+1+s}_{l+\frac{\gamma}{2}}}^2\lesssim \|
f\|_{H^{m+1}_{l+\frac{\gamma}{2}+2s}}^2+\|
f\|_{H^{m+\varrho+1}_{l+\frac{\gamma}{2} }}^2.\een  Thanks to the
interpolation inequality  \beno \|f\|_{H^k_p}^2\lesssim
\|f\|_{H^{k-\epsilon}_{2p}}\|f\|_{H^{k+\epsilon}_0},\eeno  by
iteration argument, one has that there exists a constant $r_p$ and
$\delta\in (0,1)$ such that \ben\label{in2} \|f\|_{H^k_p} \lesssim
\|f\|_{H^{k-1}_{r_p}}^{\delta}\|f\|_{H^{k+\epsilon}_0}^{1-\delta}.\een

Denote $\mathfrak{C}_m$ as the quantity depending only on
$\sup\limits_{[t_0, \infty)}\|f\|_{H^m_l(\R^3_v)}$ with $l\in \R^+$.
Then by using \eqref{in2}, \eqref{in1} can be rewritten as \beno
\f{d}{dt}\|f\|_{H^{m+1}_l}^2+
\|f\|_{H^{m+1+s}_{l+\frac{\gamma}{2}}}^2\lesssim  \mathfrak{C}_m.
\eeno Thanks to \eqref{in2} again, we also can derive that there
exists a constant $\eta>0$ such that \beno
\f{d}{dt}\|f\|_{H^{m+1}_l}^2+ \|f\|_{H^{m+1}_l}^{2+\eta} \lesssim
\mathfrak{C}_m.\eeno Using a standard argument(used by Nash for
parabolic equations), we see that for $t_1>t_0$, there holds \beno
f\in L^\infty([t_1,\infty]; H^{m+1}_l)\eeno with $l\in \R^+$. This
gives the proof to the Theorem \ref{smh}.

Now we only need to check that infinite $L^1$ moment estimate
\eqref{am1} will imply  infinite $L^2$ moment estimate \eqref{am2}
in the case of  $\gamma+2s>0$. Set $\alpha=0$, then $g=f\langle
v\rangle^{l}$.  Following the similar procedure, one may obtain that
 \beno
&&\f{d}{dt}\|g\|_{L^2}^2+\frac{C_f}{3}\|g\langle
v\rangle^{\f{\gamma}{2}}\|_{H^s}^2\\
&\lesssim & (C_f^{1-2s}\|f\|_{L^1_{|\gamma|+2}}^{2s}+\|f\|_{L^1_{|\gamma|+2}})
\|g\langle
v\rangle^{\f{\gamma}{2}}\|_{H^\varrho}^2+\|f\|_{L^1_{|\gamma|+2s+2l+2}}\|g\langle
v\rangle^{\f{\gamma}{2}}\|_{H^s}\|g\|_{L^2_{\frac{\gamma}{2}+2s}}.\eeno

Thanks to the interpolation inequality \beno
\|f\|_{H^\varrho}\lesssim
\|f\|_{L^1}^{\frac{2s-2\varrho}{2s+3}}\|f\|_{H^s}^{\frac{3+2\varrho}{2s+3}}\eeno
and \beno \|f\|_{L^2}\le \|f\langle
v\rangle^{\frac{-3r}{2s+3}}\|_{L^1}^{\frac{2s}{2s+3}}\|f\langle
v\rangle^r\|_{H^s}^{\frac{3}{2s+3}}  \eeno with $r\in \R$, one may
deduce that there exists a constant $\eta>0$ such that \beno
\f{d}{dt}\|f\langle v\rangle^{l}\|_{L^2}^2+ \|f\langle
v\rangle^{l}\|_{L^2}^{2+\eta} \lesssim \mathfrak{C}_0^1,\eeno where
  $\mathfrak{C}_0^1$ represents the quantity depending only on
$\sup\limits_{[t_0, \infty)}\|f\langle v\rangle^l\|_{L^1(\R^3_v)}$
with $l\in \R^+$. From which, we complete the proof to the Theorem
\ref{smh}. \ef

\section{Smoothing effect for the inhomogeneous Boltzmann equation}
\subsection{Hypoelliptic estimate for the transport equation}

In this section we study the transport equation which reads:
\begin{equation}
\partial_t f(t,x,v)+v\cdot\nabla_x f(t,x,v)=g(t,x,v)\label{treq}
\end{equation}
and show the following hypoelliptic estimate.
\begin{lem}\label{hypoellipticity}
Suppose $g\in L^2([0,T]\times\mathbb{T}^3;H^{-1}(\mathbb{R}^3_v))$.
Let $f\in L^2([0,T]\times\mathbb{T}^3\times\mathbb{R}^3)$ be a weak
solution of the transport equation \eqref{treq}. If we assume
$f(0,x,v)\in L^2(\mathbb{T}^3\times\mathbb{R}^3)$ and $f\in
L^2([0,T]\times\mathbb{T}^3;\dot{H}^s(\mathbb{R}^3_v))$ for some
$0<s<1$, then we have that for any $l<-\frac{3}{2}$,
\begin{equation}
\langle v\rangle^l f\in
L^2([0,T]\times\mathbb{R}^3;\dot{H}^{\frac{s}{4(4+s)}}(\mathbb{T}^3_x)).
\label{xre}
\end{equation}
\end{lem}

\noindent{\bf Proof:} Let $\tau_k$ be the translation operator in
the $x$ variable by $k$, then one has
\begin{equation*}
\tau_k f(t,x,v)=f(t,x+k,v)-f(t,x,v).
\end{equation*}
We denote the finite difference of $f$ in the $x$ variable by
\begin{equation*}
\triangle_k f(t,x,v)=\tau_k f(t,x,v)-f(t,x,v).
\end{equation*}
Using these notations, we observe that
\begin{eqnarray}
\eqref{xre}&\Leftrightarrow
&\int_0^T\int_{\mathbb{R}^3}\int_{\mathbb{T}^3}\int_{\mathbb{T}^3}\langle
v\rangle^{2l}|\triangle_k
f|^2|k|^{-3-\frac{s}{2(4+s)}}dtdvdxdk<+\infty
\nonumber\\
&\Leftrightarrow &\int_0^T\int_{\mathbb{R}^3}\langle
v\rangle^{2l}\left(\sum_{m\in\mathbb{Z}^3}|m|^{\frac{s}{2(4+s)}}\left|\hat{f}(t,m,v)\right|^2\right)
dtdv<+\infty. \label{equv}
\end{eqnarray}
\par
We now turn to prove \eqref{equv}. Let $\chi (v)\in
C^{\infty}_0(\mathbb{R}^3)$ be a test function which satisfies
$\chi(v)\geq 0$ and $\int_{\mathbb{R}^3}\chi(v)dv=1$. For any
$\epsilon>0$, we denote the regularizing sequence $\chi_{\epsilon}$
by
$\chi_{\epsilon}(v)=\epsilon^{-3}\chi\left(\frac{v}{\epsilon}\right)$
and write
\begin{equation}
\hat f(t,m,v)=[\hat f(t,m,v)- (\hat f(t,m,\cdot)\ast_v
\chi_{\epsilon })(v)]+ (\hat f(t,m,\cdot)\ast_v \chi_{\epsilon })
(v).\label{dec}
\end{equation}
We point out here that $\epsilon$ in the above equality will be
chosen later and will depend on $|m|$.
\par
We use Minkowski's inequality and Cauchy-Schwartz inequality to get
\begin{eqnarray*}
&&\int_{\mathbb{R}^3}\langle v\rangle^{2l}|\hat f(t,m,v)- (\hat f(t,m,\cdot)\ast_v \chi_{\epsilon })(v)|^2 dv\\
&\lesssim& \int_{\mathbb{R}^3}\left|\int_{\mathbb{R}^3}[\hat f(t,m,v)-\hat f(t,m,v-u)]\,\chi_{\epsilon }(u)du\right|^2 dv\\
&\lesssim& \left(\int_{\mathbb{R}^3}\left(\int_{\mathbb{R}^3}|\hat f(t,m,v)-\hat f(t,m,v-u)|^2dv\right)^{1/2}\chi_{\epsilon }(u)du\right)^2 \\
&\lesssim& \left(\int_{\mathbb{R}^3}\chi_{\epsilon}^2(u)|u|^{3+2s}du\right)\left(\int_{\mathbb{R}^6}\frac{|\hat f(t,m,v)-\hat f(t,m,v-u)|^2}{|u|^{3+2s}}dudv\right)\\
&\lesssim& \epsilon
^{2s}\left\|\hat{f}(t,m,\cdot)\right\|^2_{\dot{H}^s(\mathbb{R}^3_v)},
\end{eqnarray*}
and we then obtain
\begin{eqnarray}
&&\int_0^T\int_{\mathbb{R}^3}\langle v\rangle^{2l}\left(\sum_{m\in\mathbb{Z}^3}|m|^{\frac{s}{2(4+s)}}\left|\hat f(t,m,v)- (\hat f(t,m,\cdot)\ast_v \chi_{\epsilon })(v)\right|^2\right)dtdv\nonumber\\
&\lesssim&\int_0^T\left(\sum_{m\in\mathbb{Z}^3}|m|^{\frac{s}{2(4+s)}}\epsilon^{2s}\left\|\hat{f}(t,m,\cdot)\right\|^2_{\dot{H}^s(\mathbb{R}^3_v)}\right)dt.\label{est1}
\end{eqnarray}
\par
For the second term of the right-hand side of \eqref{dec}, we shall
use the averaging lemma introduced by \cite{bd}. We first recall
that $g\in L^2([0,T]\times\mathbb{T}^3;H^{-1}(\mathbb{R}^3_v))$
implies that
\begin{equation*}
g(t,x,v)=g_0(t,x,v)+\sum_{j=1}^3\partial_{v_j}h_j(t,x,v),
\end{equation*}
where $g_0(t,x,v)=\mathcal{F}_v^{-1}[(1+|\xi|)^{-1}\mathcal{F}_v g](t,x,v)$ and $h_j(t,x,v)=-R_j g_0(t,x,v)$, $j=1,2,3$. Here $R_j$ is Riesz transform in
$v$ variable. Then, one has $g_0,h_j\in
L^2([0,T]\times\mathbb{T}^3\times\mathbb{R}^3)$,\,\,$j=1,2,3.$
According to (2.16) in  Theorem 2.1 (averaging lemma) of \cite{bd},
we can deduce
\begin{eqnarray*}
&& \int_0^T|(\hat f(t,m, \cdot)\ast_v\chi_{\epsilon})(v)|^2dt\\
&\lesssim&|m|^{-\frac{1}{2}}\left(\|\chi_{\epsilon}(v-u)(1+|u|^2)\|_{L^{\infty}_u}
+\|\nabla \chi_{\epsilon}(v-u)(1+|u|^2)\|_{L^{\infty}_u}\right)^2\\
&&\times\left(\|\hat f(0,m, \cdot)\|^2_{L^2(\mathbb{R}^3_v)}+\|\hat f(\cdot,m, \cdot)\|^2_{L^2([0,T];L^2(\mathbb{R}^3_v))}\right.\\
&&\left.+\|\hat
g_0(\cdot,m,\cdot)\|^2_{L^2([0,T];L^2(\mathbb{R}^3_v))} +\sum_{j=1}^3\|\hat
h_j(\cdot,m,\cdot)\|^2_{L^2([0,T];L^2(\mathbb{R}^3_v))}\right).
\end{eqnarray*}
Thanks to the fact
$\|\chi_{\epsilon}(v-u)(1+|u|^2)\|_{L^{\infty}_u}\lesssim\epsilon^{-3}(1+|v|^2)$,
we get
\begin{eqnarray}
&&\int_0^T\int_{\mathbb{R}^3}\langle v\rangle^{2l}\left(\sum_{m\in\mathbb{Z}^3}|m|^{\frac{s}{2(4+s)}}\left|(\hat f(t,m,\cdot)\ast_v \chi_{\epsilon })(v)\right|^2\right)dtdv\nonumber\\
&\lesssim&\sum_{m\in\mathbb{Z}^3}|m|^{\frac{s}{2(4+s)}-\frac{1}{2}}(\epsilon^{-6}+\epsilon^{-8})\left(\|\hat f(0,m, \cdot)\|^2_{L^2(\mathbb{R}^3_v)}+\|\hat f(\cdot,m, \cdot)\|^2_{L^2([0,T];L^2(\mathbb{R}^3_v))}\right.\nonumber\\
&&\left.+\|\hat
g_0(\cdot,m,\cdot)\|^2_{L^2([0,T];L^2(\mathbb{R}^3_v))} +\sum_{j=1}^3\|\hat
h_j(\cdot,m,\cdot)\|^2_{L^2([0,T];L^2(\mathbb{R}^3_v))}\right).\label{est2}
\end{eqnarray}
\par
Now we choose $\epsilon=|m|^{-\frac{1}{4(4+s)}}$, and we can bound
\eqref{est1} since we have $f\in
L^2([0,T]\times\mathbb{T}^3;\dot{H}^s(\mathbb{R}^3_v))$. As for
\eqref{est2}, we can also bound it if we notice the fact that
$f(0,x,v)\in L^2(\mathbb{T}^3\times\mathbb{R}^3)$ and $g_0,h\in
L^2([0,T]\times\mathbb{T}^3\times\mathbb{R}^3)$. This completes the
proof of Lemma \ref{hypoellipticity}.{\hfill $ \blacksquare $ \vskip
3mm}
\par
\begin{rmk} The idea of the proof of Lemma \ref{hypoellipticity} has been
used in that of Lemma 4.2 of \cite{cdh} which is
devoted to the smoothing effects for classical solutions of the full
Landau equation. We point out that we actually proved Lemma
\ref{hypoellipticity} in the case $s=1$ there.
\end{rmk}

Now we start to prove Theorem \ref{smih}
that gives the smoothing effect for the inhomogeneous Boltzmann
equation \eqref{bol}. We first note that if no confuse occurs, we
omit the domains $\mathbb{T}^3$ and $\mathbb{R}^3$, which correspond
to variables $x$ and $v$ respectively for simplicity, and we use the shorthand $\partial^\alpha_\beta=\partial^\alpha_x\partial^\beta_v$
for any multi-indices $\alpha$ and $\beta$ hereafter.
\par
In order to prove our main result, we shall use an induction on the
number of derivatives (in variables $x$ and $v$) that can be
controlled. One step of this induction is given by
\begin{prop}\label{prop}
Let $N\geq 5$ be a given integer, and let $f$ be a smooth solution
to the Boltzmann equation \eqref{bol}. For any $l\geq 0$, we
set $h=(\partial^\alpha_\beta f) \langle v\rangle^l$ with
$|\alpha|+|\beta|\leq N$. We assume that for any $T>0$ and any
$l\geq 0$, $h\in L^{\infty}_t([0,T];L^2_{x,v})$. Then we have that
$h\in L^{\infty}_t([\tau,T];H^1_{x,v})$ for any time $\tau\in(0,T)$.
\end{prop}
\par
We only prove the above proposition in the case $\gamma+2s>0$, the
proofs for the case $\gamma+2s\leq 0$ are analogous if we use the
estimates of this case in Theorem \ref{ub}, Corollary \ref{comm} and
Theorem \ref{coc} instead of those of the case $\gamma+2s>0$. We
start the proof with improving regularity in $x$ variable.

\subsection{Gain regularity in $x$ variable}\label{gain x}
We note the domain of variable $t$ should be $[0,T]$ if we omit it
hereafter. We present following
\begin{lem}\label{lemma1}
Let $N\geq 5$ be a given integer, and let $f$ be a smooth solution
to the Boltzmann equation \eqref{bol}. We suppose that for any $T>0$
and any $l\geq 0$, $h\in L^2_{t,x,v}$ and $h(0)\in L^2_{x,v}$, where
$h$ is defined in Proposition \ref{prop}. We further suppose that
$(\partial^\alpha_\beta f)\langle v\rangle^l\in
L^{\infty}_t(L^2_{x,v})$ for $|\alpha|+|\beta|\leq N-1$. Then we
have $h\in L^{\infty}_t(L^2_{x,v})\cap L^2_{t,x}(H^s_v)$.
\end{lem}

\noindent{\bf Proof:} Using Einstein's convention for repeated
indices, we have that $h$ satisfies the equation as follows: for
$|\beta|=0$,
\begin{equation}
\partial_t h+v\cdot\nabla_x h=[\partial^\alpha_x Q(f,f)]\langle v\rangle^l,\label{h1}
\end{equation}
and for $|\beta|\geq 1$,
\begin{equation}\label{eq1}
\partial_t h+v\cdot\nabla_x h=-\beta_i(\partial^{\alpha+e_i}_{\beta-e_i}f)\langle v\rangle^l+[\partial^\alpha_\beta Q(f,f)]\langle v\rangle^l,
\end{equation}
where $e_1=(1,0,0),e_2=(0,1,0)$ and $e_3=(0,0,1)$.
\par
We only consider the case $|\beta |\geq 1$, because the estimates
for the case $|\beta |=0$ are similar (and easier). Multiplying
equation \eqref{eq1} by $h$, and then integrating on $(t,x,v)$, we
shall estimate the resulting equation term by term.
\par
It is easy to see that
\begin{equation}
\int_0^T\int_{\mathbb{T}^3}\int_{\mathbb{R}^3}(\partial _t h)
hdtdxdv=\frac{1}{2}\left(\|h(T)\|_{L^2_{x,v}}^2-\|h(0)\|_{L^2_{x,v}}^2\right).\label{e1}
\end{equation}
Since $f$ is a spatially periodic function, we get that
\begin{equation}
\int_0^T\int_{\mathbb{T}^3}\int_{\mathbb{R}^3}(v\cdot\nabla _x h) h
dtdxdv=\frac{1}{2}\int_0^T\int_{\mathbb{R}^3}v\cdot\left(\int_{\mathbb{T}^3}\nabla
_x (h^2) dx\right)dtdv =0.\label{e2}
\end{equation}
Cauchy-Schwartz gives
\begin{equation}
\left|\int_0^T\int_{\mathbb{T}^3}\int_{\mathbb{R}^3}-\beta_i(\partial^{\alpha+e_i}_{\beta-e_i}f)\langle
v\rangle^l h dtdxdv\right| \lesssim
\|(\partial^{\alpha+e_i}_{\beta-e_i}f)\langle
v\rangle^l\|_{L^2_{t,x,v}}\|h\|_{L^2_{t,x,v}}.\label{e3}
\end{equation}
\par
We write
\begin{eqnarray}
&&\int_0^T\int_{\mathbb{T}^3}\int_{\mathbb{R}^3}[\partial^\alpha_\beta
Q(f,f)]\langle v\rangle^l h dtdxdv
=\int_0^T\int_{\mathbb{T}^3}\langle[\partial^\alpha_\beta Q(f,f)]\langle v\rangle^l, h\rangle_v dtdx\nonumber\\
&=&\int_0^T\int_{\mathbb{T}^3}\langle Q(f,h), h\rangle_v dtdx
+\int_0^T\int_{\mathbb{T}^3}\langle Q(f,\partial^\alpha_\beta f)\langle v\rangle^l-Q(f,h), h\rangle_v dtdx\nonumber\\
&&+\int_0^T\int_{\mathbb{T}^3}\sum_{{{\alpha _1+\alpha _2=\alpha
}\atop{\beta _1+\beta_2=\beta }}\atop {|\alpha _1|+|\beta _1|\geq
1}}C^\alpha_{\alpha_1}C^\beta_{\beta_1}
\langle Q(\partial^{\alpha_1}_{\beta_1}f,(\partial^{\alpha_2}_{\beta_2}f)\langle v\rangle^l), h\rangle_v dtdx\nonumber\\
&&+\int_0^T\int_{\mathbb{T}^3}\sum_{{{\alpha _1+\alpha _2=\alpha
}\atop{\beta _1+\beta_2=\beta }}\atop {|\alpha _1|+|\beta _1|\geq
1}}C^\alpha_{\alpha_1}C^\beta_{\beta_1}
\langle Q(\partial^{\alpha_1}_{\beta_1}f,\partial^{\alpha_2}_{\beta_2}f)\langle v\rangle^l-Q(\partial^{\alpha_1}_{\beta_1}f,(\partial^{\alpha_2}_{\beta_2}f)\langle v\rangle^l), h\rangle_v dtdx\nonumber\\
&{\buildrel\hbox{\footnotesize def}\over
=}&(J_1)+(J_2)+(J_3)+(J_4).\label{e4}
\end{eqnarray}
Theorem \ref{coc} gives
\begin{eqnarray*}
\langle Q(f,h), h\rangle_v&\lesssim& -C_f\|h\langle
v\rangle^{\frac{\gamma}{2}}\|_{H^s_v}+(C^{1-2s}_f\|f\langle
v\rangle^{\tilde{\gamma}}\|_{L^1_v}^{2s}+ \|f\langle
v\rangle^{\tilde{\gamma}}\|_{L^1_v})\|h\langle
v\rangle^{\frac{\gamma}{2}}\|^2_{L^2_v}\\
&&+\|f\langle
v\rangle^{|\gamma|}\|_{L^1_v}\|h\langle
v\rangle^{\frac{\gamma}{2}}\|_{H^{\varrho}_v},
\end{eqnarray*}
where $\varrho<s$. We remark here that the constant $C_f$ is
uniformly with respect to $x$ variables due to the Proposition 3 of
\cite{advw} and the assumption \eqref{am5}. It is easy to get by
interpolation and Young's inequality that
\begin{equation}
\|h\langle
v\rangle^{\frac{\gamma}{2}}\|^2_{H^{\varrho}_v}\leq
\epsilon\|h\langle
v\rangle^{\frac{\gamma}{2}}\|^2_{H^s_v}+C_{\epsilon}\|h\langle
v\rangle^{\frac{\gamma}{2}}\|^2_{L^2_v},\label{ga1}
\end{equation}
which implies
\begin{eqnarray}
(J_1)&\lesssim& -(C_f-\epsilon\|f\langle v\rangle^{|\gamma|+2}\|_{L^{\infty}_t(H^2_x(L^2_v))})\|h\langle v\rangle^{\frac{\gamma}{2}}\|^2_{L^2_{t,x}(H^s_v)}+(C^{1-2s}_f\|f\langle v\rangle^{\tilde{\gamma}+2}\|_{L^{\infty}_t(H^2_x(L^2_v))}^{2s}\nonumber\\
&&+\|f\langle v\rangle^{\tilde{\gamma}+2}\|_{L^{\infty}_t(H^2_x(L^2_v))}+C_\epsilon\|f\langle
v\rangle^{|\gamma|+2}\|_{L^{\infty}_t(H^2_x(L^2_v))})\|h\langle
v\rangle^{\frac{\gamma}{2}}\|^2_{L^2_{t,x,v}},\label{e5}
\end{eqnarray}
where we use Sobolev's embedding theorem and the inequality
$\|g\|_{L^1_v}\lesssim\|g\langle v\rangle^2\|_{L^2_v}$. As for the
term $(J_2)$, in the case $s\geq \frac{1}{2}$, Corollary \ref{comm} implies
\begin{equation*}
|(J_2)|\lesssim\int_0^T\int_{\mathbb{T}^3}\|f\langle
v\rangle^{N_1}\|_{L^1_v}\|(\partial^\alpha_\beta f)\langle
v\rangle^{N_2}\|_{H^s_v}\|h\langle v\rangle^{N_3}\|
_{L^2_v}dtdx.
\end{equation*}
If we choose $N_1=|l+\frac{\gamma}{2}|+|\frac{\gamma}{2}+2s-1|+\max\{|l-2|,|l-1|\}$,
$N_2=l+\frac{\gamma}{2}$ and $N_3=\frac{\gamma}{2}+2s-1$, we get by Sobolev's embedding theorem and Young's
inequality that
\begin{equation}
|(J_2)|\lesssim \epsilon\|h\langle
v\rangle^{\frac{\gamma}{2}}\|^2_{L^2_{t,x}(H^s_v)}+C_\epsilon\|f\langle
v\rangle^{N_1+2}\|_{L^{\infty}_t(H^2_x(L^2_v))}^2 \|h\langle
v\rangle^{N_3}\|_{L^2_{t,x,v}}^2.\label{e8}
\end{equation}
In the case $s<\frac{1}{2}$, we again use Corollary \ref{comm} to get
\begin{equation*}
|(J_2)|\lesssim\int_0^T\int_{\mathbb{T}^3}\|f\langle v\rangle^{N_1}\|_{L^1_v}\|(\partial^\alpha_\beta f)\langle v\rangle^{N_2}\|_{H^\rho_v}\|h\langle v\rangle^{N_3}\|_{H^\rho_v}dtdx,
\end{equation*}
where $\rho<s$. Taking $N_1=|l+\frac{\gamma}{2}|+|\frac{\gamma}{2}|+\max\{|l-2|,|l-1|\}$, $N_2=l+\frac{\gamma}{2}$ and $N_3=\frac{\gamma}{2}$, we obtain by using
\eqref{ga1} that
\begin{equation}
|(J_2)|\lesssim
\epsilon\|f\langle
v\rangle^{N_1+2}\|_{L^{\infty}_t(H^2_x(L^2_v))}\|h\langle
v\rangle^{\frac{\gamma}{2}}\|^2_{L^2_{t,x}(H^s_v)}+C_\epsilon\|f\langle
v\rangle^{N_1+2}\|_{L^{\infty}_t(H^2_x(L^2_v))}^2 \|h\langle
v\rangle^{\frac{\gamma}{2}}\|_{L^2_{t,x,v}}^2.
\end{equation}
To deal with the term $(J_3)$, we shall consider two cases. In the
case $1\leq |\alpha_1|+|\beta_1|\leq \left[\frac{N}{2}\right]$ (for
$r\in\mathbb{R}, [r]$ denotes the maximum integer which is less than
or equal to $r$), we have $|\alpha_1|+|\beta_1|+2\leq N-1$ and
$|\alpha_2|+|\beta_2|\leq N-1$. Then Theorem \ref{ub} gives
\begin{eqnarray*}
&&\left|\int_0^T\int_{\mathbb{T}^3}\langle Q(\partial^{\alpha_1}_{\beta_1}f,(\partial^{\alpha_2}_{\beta_2}f)\langle v\rangle^l), h\rangle_v dtdx\right|\nonumber\\
&\lesssim&\|(\partial^{\alpha_1}_{\beta_1}f)\langle
v\rangle^{N_1+2}\|_{L^{\infty}_{t,x}(L^2_v)}\|(\partial^{\alpha_2}_{\beta_2}f)\langle
v\rangle^{l+N_2}\|
_{L^2_{t,x}(H^s_v)}\|h\langle v\rangle^{N_3}\|_{L^2_{t,x}(H^s_v)},
\end{eqnarray*}
where we take
$N_1=|\frac{\gamma}{2}+2s|+|\frac{\gamma}{2}|$, $N_2=\frac{\gamma}{2}+2s$,
$N_3=\frac{\gamma}{2}$.  And Young's inequality gives
\begin{equation}
|(J_3)|\lesssim\epsilon\|h\langle
v\rangle^{\frac{\gamma}{2}}\|_{L^2_{t,x}(H^s_v)}^2+C_\epsilon\|(\partial^{\alpha_1}_{\beta_1}f)\langle
v\rangle^{N_1+2}\|^2_{L^{\infty}_t(H^2_x(L^2_v))}\|(\partial^{\alpha_2}_{\beta_2}f)\langle
v\rangle^{l+N_2}\|^2_{L^2_{t,x}(H^s_v)}\label{e10}
\end{equation}
In the case $|\alpha_1|+|\beta_1|\geq \left[\frac{N}{2}\right]+1$,
we see $|\alpha_2|+|\beta_2|+2+s\leq N-1$. Again by Theorem
\ref{ub},
\begin{eqnarray*}
&&\left|\int_0^T\int_{\mathbb{T}^3}\langle Q(\partial^{\alpha_1}_{\beta_1}f,(\partial^{\alpha_2}_{\beta_2}f)\langle v\rangle^l), h\rangle_v dtdx\right|\nonumber\\
&\lesssim&\|(\partial^{\alpha_1}_{\beta_1}f)\langle
v\rangle^{N_1+2}\|_{L^2_{t,x}(L^2_v)}\|(\partial^{\alpha_2}_{\beta_2}f)\langle
v\rangle^{l+N_2}\|
_{L^{\infty}_{t,x}(H^s_v)}\|h\langle v\rangle^{N_3}\|_{L^2_{t,x}(H^s_v)},
\end{eqnarray*}
where we take
$N_1=|\frac{\gamma}{2}+2s|+|\frac{\gamma}{2}|$, $N_2=\frac{\gamma}{2}+2s$,
$N_3=\frac{\gamma}{2}$. And Young's inequality gives
\begin{equation}
|(J_3)|\lesssim\epsilon\|h\langle
v\rangle^{\frac{\gamma}{2}}\|_{L^2_{t,x}(H^s_v)}^2+C_\epsilon\|(\partial^{\alpha_1}_{\beta_1}f)\langle
v\rangle^{N_1+2}\|^2_{L^2_{t,x,v}}\|(\partial^{\alpha_2}_{\beta_2}f)\langle
v\rangle^{l+N_2}\|^2_{L^{\infty}_t(H^2_x(H^s_v))}.\label{e11}
\end{equation}
Similarly, using Corollary \ref{comm}, we can bound $(J_4)$ by the quantities which have been controlled and $\epsilon\|h\langle
v\rangle^{\frac{\gamma}{2}}\|^2_{L^2_{t,x}(H^s_v)}$ if we study the case $1\leq|\alpha_1|+|\beta_1|\leq \left[\frac{N}{2}\right]$ and the case
$|\alpha_1|+|\beta_1|\geq \left[\frac{N}{2}\right]+1$ respectively.
\par
If we choose $\epsilon$ appeared small enough and collect all the above estimates, we arrive at $h\in
L^{\infty}_t(L^2_{x,v})\cap L^2_{t,x}(H^s_v)$. This ends up
the proof. {\hfill $ \blacksquare $ \vskip 3mm}
\par
Now we begin to improve the regularity in $x$ variable. We have
\begin{lem}\label{lemma2}
Let $N\geq 5$ be a given integer, and let $f$ be a smooth solution
to the Boltzmann equation \eqref{bol}. We suppose that for any
$T>0$ and any $l\geq 0$, $h\in L^{\infty}_t(L^2_{x,v})$ and $h(0)\in
L^2_{x,v}$, where $h$ is defined in Proposition \ref{prop}. Then we
have $h\in L^2_{t,v}(\dot{H}^{\frac{s}{4(4+s)}}_x)$.
\end{lem}

\noindent{\bf Proof:} We still only consider the case $|\beta|\geq
1$. Recall that $h$ satisfies the equation \eqref{eq1}, then we get
\begin{eqnarray}
\partial_t h+v\cdot\nabla_x h&=&-\beta_i(\partial^{\alpha+e_i}_{\beta-e_i}f)\langle v\rangle^l+\sum_{{\alpha _1+\alpha _2=\alpha }\atop{\beta _1+\beta_2=\beta }}C^\alpha_{\alpha_1}C^\beta_{\beta_1}Q(\partial^{\alpha_1}_{\beta_1}f,(\partial^{\alpha_2}_{\beta_2}f)\langle v\rangle^l)\nonumber\\
&&+\sum_{{\alpha _1+\alpha _2=\alpha }\atop{\beta _1+\beta_2=\beta }}C^\alpha_{\alpha_1}C^\beta_{\beta_1}\left[Q(\partial^{\alpha_1}_{\beta_1}f,\partial^{\alpha_2}_{\beta_2}f)\langle v\rangle^l-Q(\partial^{\alpha_1}_{\beta_1}f,(\partial^{\alpha_2}_{\beta_2}f)\langle v\rangle^l)\right]\nonumber\\
&{\buildrel\hbox{\footnotesize def}\over =}&(K_1)+(K_2)+(K_3).
\end{eqnarray}
\par
It is obvious that $(K_1)\in L^2_{t,x,v}$. Theorem \ref{ub}
gives (taking $N_1=N_2=\gamma+2s$, $N_3=0$)
\begin{eqnarray*}
\|Q(\partial^{\alpha_1}_{\beta_1}f,(\partial^{\alpha_2}_{\beta_2}f)\langle
v\rangle^l)\|_{H^{-s}_v}&\lesssim &
\|(\partial^{\alpha_1}_{\beta_1}f)\langle v\rangle^{\gamma+2s}\|_{L^1_v}\|(\partial^{\alpha_2}_{\beta_2}f)\langle v\rangle^{\gamma+2s+l}\|_{H^s_v}\\
&\lesssim & \|(\partial^{\alpha_1}_{\beta_1}f)\langle
v\rangle^{\gamma+2s+2}\|_{L^2_v}\|(\partial^{\alpha_2}_{\beta_2}f)\langle
v\rangle^{\gamma+2s+l}\|_{H^s_v}.
\end{eqnarray*}
Then in the case $|\alpha_1|+|\beta_1|\leq
\left[\frac{N}{2}\right]$, one has $|\alpha_1|+|\beta_1|+2\leq N-1$,
so that
\begin{equation*}
\|Q(\partial^{\alpha_1}_{\beta_1}f,(\partial^{\alpha_2}_{\beta_2}f)\langle
v\rangle^l)\|_{L^2_{t,x}({H^{-s}_v})}\lesssim\|(\partial^{\alpha_1}_{\beta_1}f)\langle
v\rangle^{\gamma+2s+2}\|_{L^{\infty}_t(H^2_x({L^2_v}))}\|(\partial^{\alpha_2}_{\beta_2}f)\langle
v\rangle^{\gamma+2s+l}\|_{L^2_{t,x}({H^s_v})},
\end{equation*}
while in the case $|\alpha_1|+|\beta_1|\geq
\left[\frac{N}{2}\right]+1$, one has $|\alpha_2|+|\beta_2|+2+s\leq
N-1$, so that
\begin{equation*}
\|Q(\partial^{\alpha_1}_{\beta_1}f,(\partial^{\alpha_2}_{\beta_2}f)\langle
v\rangle^l)\|_{L^2_{t,x}({H^{-s}_v})}\lesssim\|(\partial^{\alpha_1}_{\beta_1}f)\langle
v\rangle^{\gamma+2s+2}\|_{L^2_{t,x,v}}\|(\partial^{\alpha_2}_{\beta_2}f)\langle
v\rangle^{\gamma+2s+l}\|_{L^{\infty}_t(H^2_x({H^s_v}))}.
\end{equation*}
We thus obtain that $(K_2)\in L^2_{t,x}(H^{-s}_v)$ by Lemma
\ref{lemma1}. In the case $s\geq \frac{1}{2}$, Corollary \ref{comm} gives (taking
$N_1=|l-1+\gamma+2s|+\max\{|l-2|,|l-1|\}$, $N_2=l-1+\gamma+2s$ and $N_3=0$)
\begin{eqnarray*}
&&\|Q(\partial^{\alpha_1}_{\beta_1}f,\partial^{\alpha_2}_{\beta_2}f)\langle v\rangle^l-Q(\partial^{\alpha_1}_{\beta_1}f,(\partial^{\alpha_2}_{\beta_2}f)\langle v\rangle^l)\|_{L^2_v}\\
&\lesssim & \|(\partial^{\alpha_1}_{\beta_1}f)\langle
v\rangle^{N_1+2}\|_{L^2_v}\|(\partial^{\alpha_2}_{\beta_2}f)\langle
v\rangle^{N_2}\|_{H^s_v}.
\end{eqnarray*}
Considering the case $|\alpha_1|+|\beta_1|\leq
\left[\frac{N}{2}\right]$ and the case $|\alpha_1|+|\beta_1|\geq
\left[\frac{N}{2}\right]+1$ respectively just as the estimates for
$(K_2)$, we can deduce that $(K_3)\in L^2_{t,x,v}$. While in the case $s<\frac{1}{2}$, again by Corollary \ref{comm}, we get for $\rho\in(0,s)$,
\begin{eqnarray*}
&&\|Q(\partial^{\alpha_1}_{\beta_1}f,\partial^{\alpha_2}_{\beta_2}f)\langle v\rangle^l-Q(\partial^{\alpha_1}_{\beta_1}f,(\partial^{\alpha_2}_{\beta_2}f)\langle v\rangle^l)\|_{H^{-\rho}_v}\\
&\lesssim & \|(\partial^{\alpha_1}_{\beta_1}f)\langle
v\rangle^{N_1+2}\|_{L^2_v}\|(\partial^{\alpha_2}_{\beta_2}f)\langle
v\rangle^{N_2}\|_{H^\rho_v}
\end{eqnarray*}
where we take $N_1=|l+\gamma|+\max\{|l-2|,|l-1|\}$, $N_2=l+\gamma$ and $N_3=0$. Still considering the case $|\alpha_1|+|\beta_1|\leq
\left[\frac{N}{2}\right]$ and the case $|\alpha_1|+|\beta_1|\geq
\left[\frac{N}{2}\right]+1$ respectively, we may obtain that $(K_3)\in L^2_{t,x}(H^{-\rho}_v)\subset L^2_{t,x}(H^{-s}_v)$.
\par
Now we have that $(K_1)+(K_2)+(K_3)\in L^2_{t,x}(H^{-s}_v)\subset
L^2_{t,x}(H^{-1}_v)$, and hence $h\in
L^2_{t,v}(\dot{H}^{\frac{s}{4(4+s)}}_x)$ if we use Lemma
\ref{hypoellipticity}. {\hfill $ \blacksquare $ \vskip 3mm}
\par
Roughly speaking, Lemma \ref{lemma2} shows that if $h\in
L^{\infty}_t(L^2_{x,v})$, then $h$ has a $\frac{s}{4(4+s)}$
derivative gain in $x$ variable. The next step is to prove that $h$
can still gain a $\frac{s}{4(4+s)}$ derivative in $x$ variable, so
that we can finally reach that $h\in L^{\infty}_t(H^1_x(L^2_v))$ if
we repeat this step several times. To this end, we present following
\begin{lem}\label{lemma3}
We denote $\delta=\frac{s}{4(4+s)}$ for simplicity. Let $N\geq 5$ be
a given integer, and let $f$ be a smooth solution to the Boltzmann
equation \eqref{bol}. We suppose that for any $T>0$ and any
$l\geq 0$, $h\in L^{\infty}_t(L^2_{x,v})\cap
L^2_{t,v}(\dot{H}^\delta_x)$ and $h(0)\in L^2_v(\dot{H}^\delta_x)$
where $h$ is defined in Proposition \ref{prop}. Then we have
$g_{\delta,k}\in L^2_{t,x,k}(H^s_v)$ with
$g_{\delta,k}=\triangle_k h|k|^{-\delta-\frac{3}{2}}$, $k\in
\mathbb{T}^3$.
\end{lem}

\noindent{\bf Proof:} We still only consider the case $|\beta|\geq
1$. The equation for $g_{\delta,k}$ in this case reads:
\begin{eqnarray}
&&\partial_t g_{\delta,k}+v\cdot\nabla_x g_{\delta,k}\nonumber\\
&=&-\beta_i\left[\frac{\triangle_k(\partial^{\alpha+e_i}_{\beta-e_i}f)}{|k|^{\delta+\frac{3}{2}}}\right]\langle
v\rangle^l
+Q(f(x+k),g_{\delta,k})\nonumber\\
&&+\left[Q\left(f(x+k),\frac{\triangle_k(\partial^{\alpha}_{\beta}f)}{|k|^{\delta+\frac{3}{2}}}\right)\langle v\rangle^l-Q(f(x+k),g_{\delta,k})\right]\nonumber\\
&&+\sum_{{{\alpha _1+\alpha _2=\alpha }\atop{\beta _1+\beta_2=\beta }}\atop {|\alpha _1|+|\beta _1|\geq 1}}C^\alpha_{\alpha_1}C^\beta_{\beta_1} Q\left(\partial^{\alpha_1}_{\beta_1}f(x+k),\frac{\triangle_k((\partial^{\alpha_2}_{\beta_2}f)(x))\langle v\rangle^l}{|k|^{\delta+\frac{3}{2}}}\right)\nonumber\\
&&+\sum_{{{\alpha _1+\alpha _2=\alpha }\atop{\beta _1+\beta_2=\beta }}\atop {|\alpha _1|+|\beta _1|\geq 1}}C^\alpha_{\alpha_1}C^\beta_{\beta_1} \left[Q\left(\partial^{\alpha_1}_{\beta_1}f(x+k),\frac{\triangle_k(\partial^{\alpha_2}_{\beta_2}f(x))}{|k|^{\delta+\frac{3}{2}}}\right)\langle v\rangle^l\right.\nonumber\\
&&\qquad\qquad\qquad\left.-Q\left(\partial^{\alpha_1}_{\beta_1}f(x+k),\frac{\triangle_k((\partial^{\alpha_2}_{\beta_2}f)(x))\langle v\rangle^l}{|k|^{\delta+\frac{3}{2}}}\right)\right]\nonumber\\
&&+\sum_{{\alpha _1+\alpha _2=\alpha }\atop{\beta _1+\beta_2=\beta }}C^\alpha_{\alpha_1}C^\beta_{\beta_1} Q\left(\frac{\triangle_k(\partial^{\alpha_1}_{\beta_1}f)}{|k|^{\delta+\frac{3}{2}}},(\partial^{\alpha_2}_{\beta_2}f)\langle v\rangle^l\right)\nonumber\\
&&+\sum_{{\alpha _1+\alpha _2=\alpha }\atop{\beta _1+\beta_2=\beta }}C^\alpha_{\alpha_1}C^\beta_{\beta_1}\left[Q\left(\frac{\triangle_k(\partial^{\alpha_1}_{\beta_1}f)}{|k|^{\delta+\frac{3}{2}}},\partial^{\alpha_2}_{\beta_2}f\right)\langle v\rangle^l-Q\left(\frac{\triangle_k(\partial^{\alpha_1}_{\beta_1}f)}{|k|^{\delta+\frac{3}{2}}},(\partial^{\alpha_2}_{\beta_2}f)\langle v\rangle^l\right)\right]\nonumber\\
&{\buildrel\hbox{\footnotesize def}\over
=}&(L_1)+(L_2)+(L_3)+(L_4)+(L_5)+(L_6)+(L_7).\label{e13}
\end{eqnarray}
\par
Multiplying the above equation by $g_{\delta,k}$, then integrating
on $(t,x,v,k)$ in the domain $[0,T]\times
\mathbb{T}^3\times\mathbb{R}^3\times\mathbb{T}^3$. Similar to the
estimates \eqref{e1}-\eqref{e3}, we can get
\begin{equation}
\int_0^T\int_{\mathbb{T}^3}\int_{\mathbb{R}^3}\int_{\mathbb{T}^3}(\partial
_t g_{\delta,k})
g_{\delta,k}dtdxdvdk=\frac{1}{2}\left(\|g_{\delta,k}(T)\|_{L^2_{x,v,k}}^2-\|g_{\delta,k}(0)\|_{L^2_{x,v,k}}^2\right),
\end{equation}
\begin{equation}
\int_0^T\int_{\mathbb{T}^3}\int_{\mathbb{R}^3}\int_{\mathbb{T}^3}(v\cdot\nabla
_x g_{\delta,k}) g_{\delta,k}
dtdxdvdk=\frac{1}{2}\int_0^T\int_{\mathbb{R}^3}\int_{\mathbb{T}^3}v\cdot\left(\int_{\mathbb{T}^3}\nabla
_x (g_{\delta,k}^2) dx\right)dtdvdk =0
\end{equation}
and
\begin{equation}
\left|\int_0^T\int_{\mathbb{T}^3}\int_{\mathbb{R}^3}\int_{\mathbb{T}^3}(L_1)
g_{\delta,k} dtdxdvdk\right| \lesssim
\|(\partial^{\alpha+e_i}_{\beta-e_i}f)\langle
v\rangle^l\|_{L^2_{t,v}(\dot{H}^\delta_x)}\|g_{\delta,k}\|_{L^2_{t,x,v,k}}.
\end{equation}
If we do the estimates like those for $(J_1)$ and $(J_2)$ in
\eqref{e4}, we can use Theorem \ref{coc} and Corollary
\ref{comm} to get
\begin{eqnarray}
&&\left|\int_0^T\int_{\mathbb{T}^3}\int_{\mathbb{R}^3}\int_{\mathbb{T}^3}(L_2) g_{\delta,k} dtdxdvdk\right|\nonumber\\
&\lesssim& -(C_f-\epsilon\|f\langle v\rangle^{|\gamma|+2}\|_{L^{\infty}_t(H^2_x(L^2_v))})\|g_{\delta,k}\langle v\rangle^{\frac{\gamma}{2}}\|^2_{L^2_{t,x,k}(H^s_v)}
+(C^{1-2s}_f\|f\langle v\rangle^{\tilde{\gamma}+2}\|_{L^{\infty}_t(H^2_x(L^2_v))}^{2s}\nonumber\\
&&+\|f\langle v\rangle^{\tilde{\gamma}+2}\|_{L^{\infty}_t(H^2_x(L^2_v))}+C_\epsilon\|f\langle
v\rangle^{|\gamma|+2}\|_{L^{\infty}_t(H^2_x(L^2_v))})\|g_{\delta,k}\langle
v\rangle^{\frac{\gamma}{2}}\|^2_{L^2_{t,x,v,k}}\label{f1}
\end{eqnarray}
and in the case $s\geq \frac{1}{2}$ (taking $N_1=|l+\frac{\gamma}{2}|+|\frac{\gamma}{2}+2s-1|+\max\{|l-2|,|l-1|\}$, $N_2=l+\frac{\gamma}{2}$ and $N_3=\frac{\gamma}{2}+2s-1$),
\begin{eqnarray}
&&\left|\int_0^T\int_{\mathbb{T}^3}\int_{\mathbb{R}^3}\int_{\mathbb{T}^3}(L_3) g_{\delta,k} dtdxdvdk\right|\nonumber\\
&\lesssim& \epsilon\|f\langle
v\rangle^{N_1+2}\|_{L^{\infty}_t(H^2_x(L^2_v))}\|g_{\delta,k}\langle
v\rangle^{\frac{\gamma}{2}}\|^2_{L^2_{t,x,k}(H^s_v)}\nonumber\\
&&+C_\epsilon\|f\langle
v\rangle^{N_1+2}\|_{L^{\infty}_t(H^2_x(L^2_v))}\|g_{\delta,k}\langle
v\rangle^{N_3}\|_{L^2_{t,x,v,k}},
\end{eqnarray}
while in the case $s<\frac{1}{2}$ (taking $M_1=|l+\frac{\gamma}{2}|+|\frac{\gamma}{2}|+\max\{|l-2|,|l-1|\}$, $M_2=l+\frac{\gamma}{2}$ and $N_3=\frac{\gamma}{2}$),
\begin{eqnarray}
&&\left|\int_0^T\int_{\mathbb{T}^3}\int_{\mathbb{R}^3}\int_{\mathbb{T}^3}(L_3) g_{\delta,k} dtdxdvdk\right|\nonumber\\
&\lesssim& \epsilon\|f\langle
v\rangle^{M_1+2}\|_{L^{\infty}_t(H^2_x(L^2_v))}\|g_{\delta,k}\langle
v\rangle^{\frac{\gamma}{2}}\|^2_{L^2_{t,x,k}(H^s_v)}\nonumber\\
&&+C_\epsilon\|f\langle
v\rangle^{M_1+2}\|_{L^{\infty}_t(H^2_x(L^2_v))}\|g_{\delta,k}\langle
v\rangle^{M_3}\|_{L^2_{t,x,v,k}}.
\end{eqnarray}
We now turn to consider the term containing $(L_4)$. In the case
$1\leq |\alpha_1|+|\beta_1|\leq \left[\frac{N}{2}\right]$, one has
$|\alpha_1|+|\beta_1|+2\leq N-1$ and $|\alpha_2|+|\beta_2|\leq N-1$,
so that Theorem \ref{ub} gives (taking $N_1=|\frac{\gamma}{2}+2s|+|\frac{\gamma}{2}|$, $N_2=\frac{\gamma}{2}+2s$ and
$N_3=\frac{\gamma}{2}$)
\begin{eqnarray}
&&\left|\int_0^T\int_{\mathbb{T}^3}\int_{\mathbb{R}^3}\int_{\mathbb{T}^3}(L_4) g_{\delta,k} dtdxdvdk\right|\nonumber\\
&\lesssim&\|(\partial^{\alpha_1}_{\beta_1}f)\langle
v\rangle^{N_1}\|_{L^{\infty}_{t,x}(L^1_v)}\left\|\frac{\triangle_k(\partial^{\alpha}_{\beta}f)}{|k|^{\delta+\frac{3}{2}}}\langle
v\rangle^{l+N_2}\right\|
_{L^2_{t,x,k}(H^s_v)}\|g_{\delta,k}\langle v\rangle^{N_3}\|_{L^2_{t,x,k}(H^s_v)}\nonumber\\
&\lesssim&\|(\partial^{\alpha_1}_{\beta_1}f)\langle
v\rangle^{N_1+2}\|_{L^{\infty}_t(H^2_x(L^2_v))}\|(\partial^{\alpha_2}_{\beta_2}f)\langle
v\rangle^{l+N_2}\|
_{L^2_{t}(\dot{H}^\delta_x(H^s_v))}\|g_{\delta,k}\langle v\rangle^{N_3}\|_{L^2_{t,x,k}(H^s_v)}\nonumber\\
&\lesssim&\epsilon\|g_{\delta,k}\langle
v\rangle^{\frac{\gamma}{2}}\|_{L^2_{t,x,k}(H^s_v)}^2+C_\epsilon\|(\partial^{\alpha_1}_{\beta_1}f)\langle
v\rangle^{N_1+2}\|^2_{L^{\infty}_t(H^2_x(L^2_v))}\|(\partial^{\alpha_2}_{\beta_2}f)\langle
v\rangle^{l+N_2}\|^2_{L^2_{t}(\dot{H}^\delta_x(H^s_v))}.
\end{eqnarray}
In the case $|\alpha_1|+|\beta_1|\geq \left[\frac{N}{2}\right]+1$,
one has $|\alpha_2|+|\beta_2|+2+\delta+s\leq N$. Again by Theorem
\ref{ub} (taking $N_1=|\frac{\gamma}{2}+2s|+|\frac{\gamma}{2}|$, $N_2=\frac{\gamma}{2}+2s$ and
$N_3=\frac{\gamma}{2}$),
\begin{eqnarray}
&&\left|\int_0^T\int_{\mathbb{T}^3}\int_{\mathbb{R}^3}\int_{\mathbb{T}^3}(L_4) g_{\delta,k} dtdxdvdk\right|\nonumber\\
&\lesssim&\|(\partial^{\alpha_1}_{\beta_1}f)\langle
v\rangle^{N_1}\|_{L^2_{t,x}(L^1_v)}\left\|\frac{\triangle_k(\partial^{\alpha}_{\beta}f)}{|k|^{\delta+\frac{3}{2}}}\langle
v\rangle^{l+N_2}\right\|
_{L^{\infty}_{t,x}(L^2_k(H^s_v))}\|g_{\delta,k}\langle v\rangle^{N_3}\|_{L^2_{t,x,k}(H^s_v)}\nonumber\\
&\lesssim&\|(\partial^{\alpha_1}_{\beta_1}f)\langle
v\rangle^{N_1+2}\|_{L^2_t(H^2_x(L^2_v))}\|(\partial^{\alpha_2}_{\beta_2}f)\langle
v\rangle^{l+N_2}\|
_{L^{\infty}_{t}(H^{2+\delta}_x(H^s_v))}\|g_{\delta,k}\langle v\rangle^{N_3}\|_{L^2_{t,x,k}(H^s_v)}\nonumber\\
&\lesssim&\epsilon\|g_{\delta,k}\langle
v\rangle^{\frac{\gamma}{2}}\|_{L^2_{t,x,k}(H^s_v)}^2+C_\epsilon\|(\partial^{\alpha_1}_{\beta_1}f)\langle
v\rangle^{N_1+2}\|^2_{L^2_t(H^2_x(L^2_v))}\|(\partial^{\alpha_2}_{\beta_2}f)\langle
v\rangle^{l+N_2}\|^2 _{L^{\infty}_{t}(H^{2+\delta}_x(H^s_v))}.
\end{eqnarray}
We point out that the estimates for
$\int_0^T\int_{\mathbb{T}^3}\int_{\mathbb{R}^3}\int_{\mathbb{T}^3}(L_i)
g_{\delta,k} dtdxdvdk$ ($i=5,6,7$) are similar to those for
$\int_0^T\int_{\mathbb{T}^3}\int_{\mathbb{R}^3}\int_{\mathbb{T}^3}(L_4)
g_{\delta,k} dtdxdvdk$, so we omit them.
\par
Therefore, if we notice $h(0)\in L^2_v(\dot{H}^\delta_x)$ and $h\in
L^2_{t,v}(\dot{H}^\delta_x)$ imply $g_{\delta,k}(0)\in L^2_{x,v,k}$
and $g_{\delta,k}\in L^2_{t,x,v,k}$, we finally get $g_{\delta,k}\in
L^2_{t,x,k}(H^s_v)$ by choosing $\epsilon$ appeared in the
above estimates sufficiently small. This completes the proof of
Lemma \ref{lemma3}.

{\hfill $ \blacksquare $ \vskip 3mm}
\par
The following lemma shows that $h$ can gain another
$\frac{s}{4(4+s)}$ derivative in $x$ variable based on Lemma
\ref{lemma2}.
\begin{lem}\label{lemma4}
Let $N\geq 5$ be a given integer, and let $f$ be a smooth solution
to the Boltzmann equation \eqref{bol}. We suppose that for any
$T>0$ and any $l\geq 0$, $h\in L^{\infty}_t(L^2_{x,v})\cap
L^2_{t,v}(\dot{H}^\delta_x)$ and $h(0)\in L^2_v(\dot{H}^\delta_x)$
where $h$ is defined in Proposition \ref{prop},
$\delta=\frac{s}{4(4+s)}$. Then we have $h\in
L^2_{t,v}(\dot{H}^{2\delta}_x)$.
\end{lem}

\noindent{\bf Proof:} Thanks to the fact
\begin{equation*}
\int_{\mathbb{T}^3}|\widehat{g_{\delta,k}}(m)|dk=C|m|^{2\delta}|\hat{h}(m)|^2,
\end{equation*}
we know that in order to get $h\in L^2_{t,v}(\dot{H}^{2\delta}_x)$,
we can equivalently prove $g_{\delta,k}\in
L^2_{t,v,k}(\dot{H}^{\delta}_x)$. We still only consider the case
$|\beta|\geq 1$. The equation for $g_{\delta,k}$ can be rewritten
as:
\begin{eqnarray}
&&\partial_t g_{\delta,k}+v\cdot\nabla_x g_{\delta,k}\nonumber\\
&=&-\beta_i\left[\frac{\triangle_k(\partial^{\alpha+e_i}_{\beta-e_i}f)}{|k|^{\delta+\frac{3}{2}}}\right]\langle
v\rangle^l
+\sum_{{\alpha _1+\alpha _2=\alpha }\atop{\beta _1+\beta_2=\beta }}C^\alpha_{\alpha_1}C^\beta_{\beta_1} Q\left(\partial^{\alpha_1}_{\beta_1}f(x+k),\frac{\triangle_k((\partial^{\alpha_2}_{\beta_2}f)(x))\langle v\rangle^l}{|k|^{\delta+\frac{3}{2}}}\right)\nonumber\\
&&+\sum_{{\alpha _1+\alpha _2=\alpha }\atop{\beta _1+\beta_2=\beta }}C^\alpha_{\alpha_1}C^\beta_{\beta_1} \left[Q\left(\partial^{\alpha_1}_{\beta_1}f(x+k),\frac{\triangle_k(\partial^{\alpha_2}_{\beta_2}f(x))}{|k|^{\delta+\frac{3}{2}}}\right)\langle v\rangle^l\right.\nonumber\\
&&\qquad\qquad\qquad\left.-Q\left(\partial^{\alpha_1}_{\beta_1}f(x+k),\frac{\triangle_k((\partial^{\alpha_2}_{\beta_2}f)(x))\langle v\rangle^l}{|k|^{\delta+\frac{3}{2}}}\right)\right]\nonumber\\
&&+\sum_{{\alpha _1+\alpha _2=\alpha }\atop{\beta _1+\beta_2=\beta }}C^\alpha_{\alpha_1}C^\beta_{\beta_1} Q\left(\frac{\triangle_k(\partial^{\alpha_1}_{\beta_1}f)}{|k|^{\delta+\frac{3}{2}}},(\partial^{\alpha_2}_{\beta_2}f)\langle v\rangle^l\right)\nonumber\\
&&+\sum_{{\alpha _1+\alpha _2=\alpha }\atop{\beta _1+\beta_2=\beta }}C^\alpha_{\alpha_1}C^\beta_{\beta_1}\left[Q\left(\frac{\triangle_k(\partial^{\alpha_1}_{\beta_1}f)}{|k|^{\delta+\frac{3}{2}}},\partial^{\alpha_2}_{\beta_2}f\right)\langle v\rangle^l-Q\left(\frac{\triangle_k(\partial^{\alpha_1}_{\beta_1}f)}{|k|^{\delta+\frac{3}{2}}},(\partial^{\alpha_2}_{\beta_2}f)\langle v\rangle^l\right)\right]\nonumber\\
&{\buildrel\hbox{\footnotesize def}\over
=}&(M_1)+(M_2)+(M_3)+(M_4)+(M_5).\label{e14}
\end{eqnarray}
It is easy to see that $(M_1)\in L^2_{t,x,v,k}$. Theorem
\ref{ub} gives (taking $N_1=N_2=\gamma+2s$, $N_3=0$)
\begin{eqnarray*}
&&\left\|Q\left(\partial^{\alpha_1}_{\beta_1}f(x+k),\frac{\triangle_k((\partial^{\alpha_2}_{\beta_2}f)(x))\langle v\rangle^l}{|k|^{\delta+\frac{3}{2}}}\right)\right\|_{H^{-s}_v}\\
&\lesssim & \|(\partial^{\alpha_1}_{\beta_1}f)(x+k)\langle
v\rangle^{\gamma+2s+2}\|_{L^2_v}\left\|\frac{\triangle_k((\partial^{\alpha_2}_{\beta_2}f)(x))}{|k|^{\delta+\frac{3}{2}}}\langle
v\rangle^{\gamma+2s+l}\right\|_{H^s_v}.
\end{eqnarray*}
Then in the case $|\alpha_1|+|\beta_1|\leq
\left[\frac{N}{2}\right]$, on has $|\alpha_1|+|\beta_1|+2\leq N-1$,
so that
\begin{eqnarray*}
&&\left\|Q\left(\partial^{\alpha_1}_{\beta_1}f(x+k),\frac{\triangle_k((\partial^{\alpha_2}_{\beta_2}f)(x))\langle v\rangle^l}{|k|^{\delta+\frac{3}{2}}}\right)\right\|_{L^2_{t,x,k}({H^{-s}_v})}\\
&\lesssim & \|(\partial^{\alpha_1}_{\beta_1}f)\langle
v\rangle^{\gamma+2s+2}\|_{L^{\infty}_t(H^2_x({L^2_v}))}\|(\partial^{\alpha_2}_{\beta_2}f)\langle
v\rangle^{\gamma+2s+l}\|_{L^2_t(\dot{H}^\delta_x({H^s_v}))}.
\end{eqnarray*}
while in the case $|\alpha_1|+|\beta_1|\geq
\left[\frac{N}{2}\right]+1$, one has
$|\alpha_2|+|\beta_2|+2+\delta+s\leq N$, so that
\begin{eqnarray*}
&&\left\|Q\left(\partial^{\alpha_1}_{\beta_1}f(x+k),\frac{\triangle_k((\partial^{\alpha_2}_{\beta_2}f)(x))\langle v\rangle^l}{|k|^{\delta+\frac{3}{2}}}\right)\right\|_{L^2_{t,x,k}({H^{-s}_v})}\\
&\lesssim & \|(\partial^{\alpha_1}_{\beta_1}f)\langle
v\rangle^{\gamma+2s+2}\|_{L^2_{t,x,v}}\|(\partial^{\alpha_2}_{\beta_2}f)\langle
v\rangle^{\gamma+2s+l}\|_{L^{\infty}_t(H^{2+\delta}_x({H^s_v}))}.
\end{eqnarray*}
We thus obtain that $(M_2)\in L^2_{t,x,k}(H^{-s}_v)$ by Lemma
\ref{lemma3}. Analogously, we can get $(M_3)$, $(M_4)$, $(M_5)\in L^2_{t,x,k}(H^{-s}_v)$. We now know that the
right-hand side of equation \eqref{e14} belongs to
$L^2_{t,x,k}(H^{-s}_v)$, and then we conclude $g_{\delta,k}\in
L^2_{t,v,k}(\dot{H}^{\delta}_x)$ thanks to Lemma
\ref{hypoellipticity}. {\hfill $ \blacksquare $ \vskip 3mm}
\par
Applying Lemmas \ref{lemma1} and \ref{lemma2}, we get $h\in
L^2([0,T];L^2_v(\dot{H}^\delta_x))$ with $\delta=\frac{s}{4(4+s)}$.
Then for any $t_*\in(0,T)$, we can find some time $t_1\in (0,t_*)$
such that $h(t_1)\in L^2_v(\dot{H}^\delta_x)$. So we can use Lemmas
\ref{lemma3} and \ref{lemma4} to obtain $h\in
L^2([t_1,T];L^2_v(\dot{H}^{2\delta}_x))$. As a consequence, we can
find some time $t_2\in(t_1,t_*)$ such that $h(t_2)\in
L^2_v(\dot{H}^{2\delta}_x)$. If we repeat this procedure (by using
Lemmas \ref{lemma3} and \ref{lemma4}) $m-1$ times such that
$m\delta\geq 1$, we obtain $h\in L^2([t_{m-1},T];L^2_v(H^1_x))$ for
some time $t_{m-1}\in (t_{m-2},t_* )$, and we can find some time
$t_m\in(t_{m-1},t_*)$ such that $h(t_m)\in L^2_v(H^1_x)$. Therefore,
thanks to Lemma \ref{lemma1}, we finally get $h\in
L^{\infty}([t_*,T]; L^2_v(H^1_x))$.

\subsection{Gain regularity in $v$ variable}
In this subsection we shall improve regularity in $v$ variable. We
begin with
\begin{lem}\label{lemma5}
Let $N\geq 5$ be a given integer, and let $f$ be a smooth solution
to the Boltzmann equation \eqref{bol}. We suppose that for any
$T>0$ and any $l\geq 0$, $h\in L^{\infty}_t(H^1_x(L^2_v))\cap
L^2_{t,x}(H^s_v)$ and $h(0)\in L^2_{x}(H^s_v)$ where $h$
is defined in Proposition \ref{prop}. Then we have $h\in
L^{\infty}_t(L^2_x(H^s_v))\cap L^2_{t,x}(H^{2s}_v)$.
\end{lem}

\noindent{\bf Proof:} Let $\tau_u$ be the translation operator in
the $v$ variable by $u$, then one has
\begin{equation*}
\tau_u f(t,x,v)=f(t,x,v+u)-f(t,x,v).
\end{equation*}
We denote the finite difference of $f$ in the $v$ variable by
\begin{equation}
\triangle_u f(t,x,v)=\tau_u f(t,x,v)-f(t,x,v).\label{lu}
\end{equation}
If we define
$g_{s,u}(t,x,v)=\triangle_uh(t,x,v)|u|^{-s-\frac{3}{2}}$, then we know that in order to get $h\in
L^{\infty}_t(L^2_x(H^s_v))\cap L^2_{t,x}(H^{2s}_v)$, we
can equivalently prove $g_{s,u}\in L^{\infty}_t(L^2_{x,v,u})\cap
L^2_{t,x,u}(H^s_v)$. Since $h\in L^{\infty}_t(H^1_x(L^2_v))\cap
L^2_{t,x}(H^s_v)$, we can restrict the integral domain of variable $u$
to $\mathbb{B}_1=\{u\in\mathbb{R}^3:|u|\leq 1\}$.
\par
We still only consider the case $|\beta|\geq 1$. Firstly, for any
function $p(v)$ and $q(v)$, one has
\begin{equation*}
\triangle_u(p(v)q(v))=p(v+u)\triangle_u q(v)+p(v)\triangle_u q(v).
\end{equation*}
Secondly, the translation invariance of the collision operator with
respect to the variable $v$ gives that (see
\cite{dw} for instance)
\begin{equation*}
\tau_uQ(f,g)=Q(\tau_u f,\tau_u g).
\end{equation*}
Applying these two equalities, we get the equation for $g_{s,u}$ as
follows:
\begin{eqnarray}
&&\partial_t g_{s,u}+v\cdot \nabla_x g_{s,u}\nonumber\\
&=&-u\cdot\nabla_x h(v+u)|u|^{-s-\frac{3}{2}}-\beta_i\left[(\partial^{\alpha+e_i}_{\beta-e_i}f)(v+u)\right]\frac{\triangle_u\langle v\rangle^l}{|u|^{s+\frac{3}{2}}}\nonumber\\
&&-\beta_i\left[\frac{\triangle_u(\partial^{\alpha+e_i}_{\beta-e_i}f)}{|u|^{\delta+\frac{3}{2}}}\right]\langle
v\rangle^l
+\left[\partial^\alpha_\beta Q(f(v+u),f(v+u))\right]\frac{\triangle_u\langle v\rangle^l}{|u|^{s+\frac{3}{2}}}\nonumber\\
&&+\left[\partial^\alpha_\beta Q\left(f(v+u),\frac{\triangle_u
f}{|u|^{s+\frac{3}{2}}}\right)\right]\langle v\rangle^l
+\left[\partial^\alpha_\beta Q\left(\frac{\triangle_u f}{|u|^{s+\frac{3}{2}}},f\right)\right]\langle v\rangle^l\nonumber\\
&{\buildrel\hbox{\footnotesize def}\over
=}&(P_1)+(P_2)+(P_3)+(P_4)+(P_5)+(P_6).
\end{eqnarray}
\par
Multiplying the above equation by $g_{s,u}$, then integrating on
$(t,x,v,u)$ in the domain $[0,T]\times
\mathbb{T}^3\times\mathbb{R}^3\times\mathbb{B}_1$. Similar to the
estimates \eqref{e1} and \eqref{e2}, we yield
\begin{equation}
\int_0^T\int_{\mathbb{T}^3}\int_{\mathbb{R}^3}\int_{\mathbb{B}_1}(\partial
_t g_{s,u})
g_{s,u}dtdxdvdu=\frac{1}{2}\left(\|g_{s,u}(T)\|_{L^2_{x,v,u}}^2-\|g_{s,u}(0)\|_{L^2_{x,v,u}}^2\right)
\end{equation}
and
\begin{equation}
\int_0^T\int_{\mathbb{T}^3}\int_{\mathbb{R}^3}\int_{\mathbb{B}_1}(v\cdot\nabla
_x g_{s,u}) g_{s,u}
dtdxdvdu=\frac{1}{2}\int_0^T\int_{\mathbb{R}^3}\int_{\mathbb{B}_1}v\cdot\left(\int_{\mathbb{T}^3}\nabla
_x (g_{s,u}^2) dx\right)dtdvdu =0.
\end{equation}
Cauchy-Schwartz inequality gives
\begin{eqnarray}
\left|\int_0^T\int_{\mathbb{T}^3}\int_{\mathbb{R}^3}\int_{\mathbb{B}_1}(P_1)g_{s,u}dtdxdvdu\right|
&\lesssim&\int_0^T\int_{\mathbb{T}^3}\left(\int_{\mathbb{B}_1}|u|^{-s-\frac{1}{2}}\|h\|_{\dot{H}^1_x}\|g_{s,u}\|_{L^2_v}du\right)dtdx\nonumber\\
&\lesssim&\|h\|_{L^2_{t,v}(\dot{H}^1_x)}\|g_{s,u}\|_{L^2_{t,x,u,v}}.
\end{eqnarray}
We have that
\begin{equation}
|\triangle_u\langle v\rangle^l|=|\langle v+u\rangle^l-\langle
v\rangle^l|\lesssim\int_0^1\langle v+\theta
u\rangle^{l-1}d\theta|u|\lesssim\langle v+u\rangle^{l-1}|u|,
\label{ff}
\end{equation}
where the last inequality holds due to the fact $|u|\leq 1$. Then,
\begin{eqnarray}
&&\left|\int_0^T\int_{\mathbb{T}^3}\int_{\mathbb{R}^3}\int_{\mathbb{B}_1}(P_2)g_{s,u}dtdxdvdu\right|\nonumber\\
&\lesssim&\int_0^T\int_{\mathbb{T}^3}\left(\int_{\mathbb{B}_1}|u|^{-s-\frac{1}{2}}\|(\partial^{\alpha+e_i}_{\beta-e_i}f)(v+u)\langle v+u\rangle^{l-1}\|_{L^2_v}\|g_{s,u}\|_{L^2_v}du\right)dtdx\nonumber\\
&\lesssim&\|(\partial^{\alpha+e_i}_{\beta-e_i}f)\langle
v\rangle^{l-1}\|_{L^2_{t,x,v}}\|g_{s,u}\|_{L^2_{t,x,u,v}}.
\end{eqnarray}
Again by Cauchy-Schwartz inequality, we get
\begin{equation}
\left|\int_0^T\int_{\mathbb{T}^3}\int_{\mathbb{R}^3}\int_{\mathbb{B}_1}(P_3)g_{s,u}dtdxdvdu\right|\lesssim\|\partial^{\alpha+e_i}_{\beta-e_i}f\|_{L^2_{t,x}(H^s_v)}
\|g_{s,u}\langle v\rangle^l\|_{L^2_{t,x,u,v}}.
\end{equation}
We now turn to treat the term
$\int_0^T\int_{\mathbb{T}^3}\int_{\mathbb{R}^3}\int_{\mathbb{B}_1}(P_5)g_{s,u}dtdxdvdu$.
We write
\begin{eqnarray}
(P_5)&=&Q(f(v+u),g_{s,u})-Q\left(f(v+u),(\partial^\alpha_\beta f)(v+u)\frac{\triangle_u\langle v\rangle^l}{|u|^{s+\frac{3}{2}}}\right)\nonumber\\
&&+\left[Q\left(f(v+u),\frac{\triangle_u(\partial^\alpha_\beta f)}{|u|^{s+\frac{3}{2}}}\right)\langle v\rangle^l-Q\left(f(v+u),\frac{\triangle_u(\partial^\alpha_\beta f)\langle v\rangle^l}{|u|^{s+\frac{3}{2}}}\right)\right]\nonumber\\
&&+\sum_{{{\alpha _1+\alpha _2=\alpha }\atop{\beta _1+\beta_2=\beta }}\atop {|\alpha _1|+|\beta _1|\geq 1}}C^\alpha_{\alpha_1}C^\beta_{\beta_1} Q\left((\partial^{\alpha_1}_{\beta_1}f)(v+u),\frac{\triangle_u((\partial^{\alpha_2}_{\beta_2} f)\langle v\rangle^l)}{|u|^{s+\frac{3}{2}}}\right)\nonumber\\
&&-\sum_{{{\alpha _1+\alpha _2=\alpha }\atop{\beta _1+\beta_2=\beta }}\atop {|\alpha _1|+|\beta _1|\geq 1}}C^\alpha_{\alpha_1}C^\beta_{\beta_1} Q\left((\partial^{\alpha_1}_{\beta_1}f)(v+u),(\partial^{\alpha_2}_{\beta_2} f)(v+u)\frac{\triangle_u\langle v\rangle^l}{|u|^{s+\frac{3}{2}}}\right)\nonumber\\
&&+\sum_{{{\alpha _1+\alpha _2=\alpha }\atop{\beta _1+\beta_2=\beta }}\atop {|\alpha _1|+|\beta _1|\geq 1}}C^\alpha_{\alpha_1}C^\beta_{\beta_1} \left[Q\left((\partial^{\alpha_1}_{\beta_1}f)(v+u),\frac{\triangle_u(\partial^{\alpha_2}_{\beta_2} f)}{|u|^{s+\frac{3}{2}}}\right)\langle v\rangle^l\right.\nonumber\\
&&\qquad\qquad\qquad\left.-Q\left((\partial^{\alpha_1}_{\beta_1}f)(v+u),\frac{\triangle_u(\partial^{\alpha_2}_{\beta_2}
f)\langle v\rangle^l}{|u|^{s+\frac{3}{2}}}\right)
\right]\nonumber\\
&{\buildrel\hbox{\footnotesize def}\over
=}&(P_5)_1+(P_5)_2+(P_5)_3+(P_5)_4+(P_5)_5+(P_5)_6.
\end{eqnarray}
Applying Theorem \ref{coc}, we get
\begin{eqnarray}
&&\left|\int_0^T\int_{\mathbb{T}^3}\int_{\mathbb{R}^3}\int_{\mathbb{B}_1}(P_5)_1 g_{s,u} dtdxdvdu\right|\nonumber\\
&\lesssim& -(C_f-\epsilon\|f\langle v\rangle^{|\gamma|+2}\|_{L^{\infty}_t(H^2_x(L^2_v))})\|g_{s,u}\langle v\rangle^{\frac{\gamma}{2}}\|^2_{L^2_{t,x,u}(H^s_v)}+(C^{1-2s}_f\|f\langle v\rangle^{\tilde{\gamma}+2}\|_{L^{\infty}_t(H^2_x(L^2_v))}^{2s}\nonumber\\
&&+\|f\langle v\rangle^{\tilde{\gamma}+2}\|_{L^{\infty}_t(H^2_x(L^2_v))}+C_\epsilon\|f\langle
v\rangle^{|\gamma|+2}\|_{L^{\infty}_t(H^2_x(L^2_v))})\|g_{s,u}\langle
v\rangle^{\frac{\gamma}{2}}\|^2_{L^2_{t,x,v,u}}.
\end{eqnarray}
Remembering $|u|\leq 1$, We then get by Theorem \ref{ub} and
\eqref{ff} that
\begin{eqnarray}
&&\left|\int_0^T\int_{\mathbb{T}^3}\int_{\mathbb{R}^3}\int_{\mathbb{B}_1}(P_5)_2 g_{s,u} dtdxdvdu\right|\nonumber\\
&\lesssim&\epsilon\|g_{s,u}\langle
v\rangle^{\frac{\gamma}{2}}\|^2_{L^2_{t,x,u}(H^s_v)}+C_\epsilon\|f\langle
v\rangle^{|\frac{\gamma}{2}+2s|+|\frac{\gamma}{2}|+2}\|^2_{L^{\infty}_t(H^2_x(L^2_v))}\|(\partial^\alpha_\beta
f)\langle v\rangle^{l-1+\frac{\gamma}{2}+2s}\|^2_{L^2_{t,x}(H^s_v)}.
\end{eqnarray}
As for the terms
$\int_0^T\int_{\mathbb{T}^3}\int_{\mathbb{R}^3}\int_{\mathbb{B}_1}(P_5)_i
g_{s,u} dtdxdvdu$ $(i=3,4)$, we claim that we can bound them by the quantities
which have been controlled and $\epsilon\|g_{s,u}\langle
v\rangle^{\frac{\gamma}{2}}\|^2_{L^2_{t,x,u}(H^s_v)}$. The
proof is very close to the estimates for the terms
$\int_0^T\int_{\mathbb{T}^3}\int_{\mathbb{R}^3}\int_{\mathbb{T}^3}(L_i)
g_{\delta,k} dtdxdvdk$ $(i=3,4)$ in Lemma \ref{lemma3}, so we omit it. And we
can analogously bound the terms
$\int_0^T\int_{\mathbb{T}^3}\int_{\mathbb{R}^3}\int_{\mathbb{B}_1}(P_5)_j
g_{s,u} dtdxdvdu$ for $j=5,6$. The estimates for
$\int_0^T\int_{\mathbb{T}^3}\int_{\mathbb{R}^3}\int_{\mathbb{B}_1}(P_k)
g_{s,u} dtdxdvdu$ $(k=4,6)$ are similar.
\par
Therefore, the lemma is proved by taking $\epsilon$ appeared small
enough. {\hfill $ \blacksquare $ \vskip 3mm}
\par
We now finish the proof of Proposition \ref{prop}. According to the result of Subsection \ref{gain x}, one has for any $t_*\in(0,\tau)$, $h\in
L^{\infty}([t_*,T];L^2_v(H^1_x))\cap
L^2([t_*,T];L^2_x(H^s_v))$. Then we
can find some time $\tau_1\in(t_*,\tau)$ such that $h(\tau_1)\in
L^2_x(H^s_v)$. Thanks to Lemma \ref{lemma5}, we get $h\in
L^{\infty}([\tau_1,T];L^2_x(H^s_v))\cap
L^2([\tau_1,T];L^2_x(H^{2s}_v))$. As a consequence, we can
find some time $\tau_2\in(\tau_1,\tau)$ such that $h(\tau_2)\in
L^2_x(H^{2s}_v))$. Repeating this procedure (by using Lemma
\ref{lemma5}) $m$ times such that $ms\geq 1$, we finally obtain
$h\in L^{\infty}([\tau_m,T];L^2_x(H^1_v))\cap
L^2([\tau_m,T];L^2_x(H^{(m+1)s}_v))$ (for some time $\tau_m\in
(\tau_{m-1},\tau)$), which implies $h\in
L^{\infty}([\tau,T];L^2_x(H^1_v))$. Since we already have
$h\in L^{\infty}([\tau,T];L^2_v(H^1_x))$, Proposition \ref{prop} is
therefore proven.
{\hfill $ \blacksquare $ \vskip 3mm}

\subsection{Proof of Theorem \ref{smih}}
We now end up the proof of Theorem \ref{smih}. By applying
Proposition \ref{prop} repeatedly, we get that for any $l\geq 0$ and
any $0<\tau<T<+\infty$,
\begin{equation}
f\langle v\rangle^l\in
L^{\infty}([\tau,T];H^{\infty}_{x,v}).\label{g1}
\end{equation}
\par
We claim that for any nonnegative integer $n$,
\begin{equation}
(\partial_t^n f)\langle v\rangle^l\in
L^{\infty}([\tau,T];H^{\infty}_{x,v}).\label{g2}
\end{equation}
We shall prove this by induction on $n$. Thanks to \eqref{g1}, this
is true for $n=0$. If we assume that \eqref{g2} holds for any
integer $k\leq n$, then for any $l\geq 0$ and any multi-indices
$\alpha$ and $\beta$,
\begin{eqnarray*}
(\partial^\alpha_\beta\partial_t^{n+1}f)\langle v\rangle^l&=&-\partial^\alpha_\beta[v\cdot\nabla_x(\partial_t^n f)]\langle v\rangle^l\\
&&+\sum_{{\alpha _1+\alpha _2=\alpha }\atop{\beta _1+\beta_2=\beta
}}\sum_{m=0}^n C^\alpha_{\alpha_1}C^\beta_{\beta_1}C^n_m
Q(\partial^{\alpha_1}_{\beta_1}\partial_t^m
f,\partial^{\alpha_2}_{\beta_2}\partial_t^{n-m}f)\langle v\rangle^l\\
&{\buildrel\hbox{\footnotesize def}\over =}&(R_1)+(R_2).
\end{eqnarray*}
It is obvious from the
induction hypothesis that $(R_1)\in L^{\infty}([\tau,T];L^2_{x,v})$.
For any nonnegative function $W\in L^1([\tau,T];L^2_{x,v})$, we shall consider the quantity
\begin{equation*}
\langle Q(U,V)\langle v\rangle^l,W\rangle_v,
\end{equation*}
where we set $U=\partial^{\alpha_1}_{\beta_1}\partial_t^m
f$ and $V=\partial^{\alpha_2}_{\beta_2}\partial_t^{n-m}f$ for simplicity. We write
\begin{eqnarray}
&&\langle Q(U,V)\langle v\rangle^l,W\rangle_v\nonumber\\
&=&\int_{\mathbb{R}^6}dvdv_*\int_{\mathbb{S}^2}B(v-v_*,\sigma)U_*V(W'\langle v'\rangle^l-W\langle v\rangle^l)d\sigma\nonumber\\
&=&\int_{\mathbb{R}^6}dvdv_*\int_{\mathbb{S}^2}B(v-v_*,\sigma)U_*(V'W'\langle v'\rangle^l-VW\langle v\rangle^l)d\sigma\nonumber\\
&&+\int_{\mathbb{R}^6}dvdv_*\int_{\mathbb{S}^2}B(v-v_*,\sigma)U_*W'(V\langle v\rangle^l-V'\langle v'\rangle^l)d\sigma\nonumber\\
&&+\int_{\mathbb{R}^6}dvdv_*\int_{\mathbb{S}^2}B(v-v_*,\sigma)U_*V'W'(\langle v'\rangle^l-\langle v\rangle^l)d\sigma\nonumber\\
&&+\int_{\mathbb{R}^6}dvdv_*\int_{\mathbb{S}^2}B(v-v_*,\sigma)U_*W'(V-V')(\langle v'\rangle^l-\langle v\rangle^l)d\sigma\nonumber\\
&{\buildrel\hbox{\footnotesize def}\over =}&(S_1)+(S_2)+(S_3)+(S_4).
\end{eqnarray}
We only give the estimates for $(S_2)$, and those for $(S_1), (S_3)$ and $(S_4)$ are analogous. Applying Taylor expansion formula up to order 2, we get
\begin{equation*}
V\langle v\rangle^l-V'\langle v'\rangle^l=(v-v')\cdot\nabla_v(V'\langle v'\rangle^l)+\int_0^1(v-v')\otimes(v-v'):\nabla^2_v(V(\gamma(\kappa))\langle \gamma(\kappa)
\rangle^l)d\kappa,
\end{equation*}
where $\gamma(\kappa)=\kappa v'+(1-\kappa)v$. If we change the variables from $v$ to $v'$, and then use \eqref{van}, we deduce that
\begin{equation}
\int_{\mathbb{R}^6}dvdv_*\int_{\mathbb{S}^2}B(v-v_*,\sigma)U_*W'(v-v')\cdot\nabla_v(V'\langle v'\rangle^l)d\sigma=0.
\end{equation}
So we only need to study the term
\begin{equation*}
\int_0^1d\kappa\int_{\mathbb{R}^6}dvdv_*\int_{\mathbb{S}^2}B(v-v_*,\sigma)U_*W'(v-v')\otimes(v-v'):\nabla^2_v(V(\gamma(\kappa))\langle \gamma(\kappa)
\rangle^l)d\sigma{\buildrel\hbox{\footnotesize def}\over =}(Y).
\end{equation*}
Thanks to the fact $|v-v'|^2=|v-v_*|^2\sin^2\frac{\theta}{2}$, we get by Cauchy-Schwartz inequality that
\begin{eqnarray*}
&&|(Y)|\nonumber\\
&\lesssim&\left(\int_{\mathbb{R}^6}dvdv_*\int_{\mathbb{S}^2}B(v-v_*,\sigma)|v-v_*|^{-\gamma}U_*W'^2\sin^2\frac{\theta}{2}d\sigma\right)^{\frac{1}{2}}\nonumber\\
&&\,\times\left(\int_0^1d\kappa\int_{\mathbb{R}^6}dvdv_*\int_{\mathbb{S}^2}B(v-v_*,\sigma)|v-v_*|^{4+\gamma}U_*|\nabla^2_v(V(\gamma(\kappa))\langle \gamma(\kappa)
\rangle^l)|^2\sin^2\frac{\theta}{2}d\sigma\right)^{\frac{1}{2}}\nonumber\\
&{\buildrel\hbox{\footnotesize def}\over =}&(Y_1)^{\frac{1}{2}}(Y_2)^{\frac{1}{2}}.
\end{eqnarray*}
Changing the variables from $v$ to $v'$, we get from \eqref{ja} with $\kappa=1$ that
\begin{equation*}
(Y_1)\lesssim\int_0^{\frac{\pi}{2}}\theta^{1-2s}d\theta\int_{\mathbb{R}^6}U_*W'^2dv'dv_*\lesssim\|U\langle v\rangle^2\|_{L^2_v}\|W\|_{L^2_v}^2.
\end{equation*}
Similarly, changing the variables from $v$ to $\gamma(\kappa)=u$, we get from \eqref{ja} that
\begin{equation*}
(Y_2)\lesssim\int_0^{\frac{\pi}{2}}\theta^{1-2s}d\theta\int_{\mathbb{R}^6}U_*|u-v_*|^{4+2\gamma}|\nabla^2_v(V(u)\langle u\rangle^l)|^2dudv_*.
\end{equation*}
Noticing $\gamma+2s>-1$ implies $4+2\gamma>-2$, we get in the case $4+2\gamma\geq 0$,
\begin{equation*}
(Y_2)\lesssim\|U\langle v\rangle^{6+2\gamma}\|_{L^2_v}\|V\langle v\rangle^{l+4+2\gamma}\|_{H^2_v}^2,
\end{equation*}
while in the case $-2<4+2\gamma<0$,
\begin{eqnarray*}
(Y_2)&\lesssim&\int_{\mathbb{R}^6}U_*(1+|u-v_*|^{4+2\gamma}\mathbf{1}_{|v-v_*|\le
1})|\nabla^2_v(V(u)\langle u\rangle^l)|^2dudv_*\\
&\lesssim&(\|U\|_{L^1_v}+\|U\|_{L^{\infty}_v})\|V\langle v\rangle^l\|_{H^2_v}\lesssim\|U\langle v\rangle^2\|_{H^2_v}\|V\langle v\rangle^l\|_{H^2_v}.
\end{eqnarray*}
Then we arrive at
\begin{equation}
|(Y)|\lesssim(\|U\langle v\rangle^{6+2\gamma}\|_{L^2_v}+\|U\langle v\rangle^2\|_{H^2_v})\|V\langle v\rangle^{l+6+2\gamma}\|_{H^2_v}\|W\|_{L^2_v},
\end{equation}
which together with the induction hypothesis implies $(S_2)\in  L^{\infty}([\tau,T];L^2_{x,v})$.
We now conclude that \eqref{g2} holds true, that is, for any $l\geq
0$ and any $0<\tau<T<+\infty$,
\begin{equation*}
f\langle v\rangle^l\in W^{\infty,\infty}([\tau,T];H^{\infty}_{x,v}).
\end{equation*}
\par
We point out Theorem 1 of \cite{gs2} shows that for
any $l\geq 0$ and any integer $N\geq 5$, $f\langle v\rangle^l\in
L^{\infty}([0,+\infty);H^N_{x,v})$ as long as $\|F_0\langle
v\rangle^l\|_{H^N_{x,v}}\leq \eta_0$ for some $\eta_0>0$. Then it is
easy to check that our estimates established up to now can be made
independent on the time $T$, so that we actually obtain
\begin{equation*}
f\langle v\rangle^l\in
W^{\infty,\infty}([\tau,\infty);H^{\infty}_{x,v}).
\end{equation*}
This completes the proof of Theorem \ref{smih}. {\hfill $
\blacksquare $ \vskip 3mm}
\par
\begin{rmk}
We note that for the purpose of having a completely rigorous proof of Theorem \ref{smih},
all the estimates in the proof should actually be made on a version
of the Boltzmann equation \eqref{bol} with smooth data and
then extended to the solution under consideration by a passage to
the limit. This creates no difficulty.
\end{rmk}

\section{Appendix}
In this appendix, we give the proof to the following proposition:
\begin{prop}
Suppose $\gamma>0$ and $0<s<1$. Then for any $R\geq 1$ and   any smooth
function $\chi_R$ defined as $\chi_R=\chi(\frac{\cdot}{R})$ with
$0\le \chi\le 1, \chi=1$ on $B_1$ and $\mathrm{supp}\,(\chi)\subset
B_2$, there holds \beno \|f\|_{\dot{H}^s(\R^3_v)}\le
\|f\chi_R\|_{\dot{H}^s( \R^3_v)}+ R^{-\f{\gamma}{2}}\|f\langle
v\rangle^{\frac{\gamma}{2}}\|_{\dot{H}^s(\R^3_v)}+\|f\langle
v\rangle^{\frac{\gamma}{2}}\|_{L^2(\R^3_v)}\eeno
\end{prop}
\noindent{\bf Proof:} We first recall that for $0<s<1$ and smooth
function $\phi$, there hold \beno
\|f\phi\|_{\dot{H}^s}^2=\int_{\R^6}
\frac{|f\phi(v)-f\phi(w)|^2}{|v-w|^{3+2s}}dvdw,\eeno and \beno
\|f\|_{\dot{H}^s(\R^3_v)}\le \|f\chi_R\|_{\dot{H}^s( \R^3_v)}+
 \|f(1-\chi_R)\|_{\dot{H}^s(\R^3_v)}+\|f \|_{L^2(\R^3_v)}. \eeno

We also note that \beno
|f\phi(v)-f\phi(w)|=|\phi(v)(f(v)-f(w))+f(w)(\phi(v)-\phi(w))|.
\eeno Then on the one hand, one may deduce that \beno
\|f\phi\|_{\dot{H}^s}^2&\lesssim &\|f\phi\|_{L^2}^2+\int_{\R^6}
\frac{|f(v)-f(w) |^2}{|v-w|^{3+2s}}|\phi(v)|^2\mathbf{1}_{|v-w|\le
1}dvdw\\&&\quad + \|\na \phi\|^2_{L^\infty}\|f\|_{L^2}^2, \eeno
which implies that \ben\label{ap1} \|
f(1-\chi_R)\|_{\dot{H}^s}^2&\lesssim
&\|f(1-\chi_R)\|_{L^2}^2+\int_{\R^6} \frac{|f(v)-f(w)
|^2}{|v-w|^{3+2s}}|(1-\chi_R)|^2\mathbf{1}_{|v-w|\le
1}dvdw\nonumber\\&&\quad + \f1{R} \|\na
\chi\|^2_{L^\infty}\|f\|_{L^2}^2. \een

On the other hand, one has \beno \|f\phi\|_{\dot{H}^s}^2&\gtrsim &
\f12\int_{\R^6} \frac{|f(v)-f(w)
|^2}{|v-w|^{3+2s}}|\phi(v)|^2\mathbf{1}_{|v-w|\le 1}dvdw\\&&\quad
-\int_{\R^6} \frac{ |\phi(v)-\phi(w)
|^2}{|v-w|^{3+2s}}|f(w)|^2\mathbf{1}_{|v-w|\le 1}dvdw.\eeno
Observing that \beno  |\langle v\rangle^{\frac{\gamma}{2}}-\langle
w\rangle^{\frac{\gamma}{2}}|^2\mathbf{1}_{|v-w| \le 1} \lesssim
|v-w|^2\langle w\rangle^{\gamma-2},\eeno one may obtain that
\ben\label{ap2} \|f\langle
v\rangle^{\frac{\gamma}{2}}\|_{\dot{H}^s}^2\gtrsim \f12\int_{\R^6}
\frac{|f(v)-f(w) |^2}{|v-w|^{3+2s}}\langle
v\rangle^{\gamma}\mathbf{1}_{|v-w|\le 1}dvdw-\|f\langle
v\rangle^{\frac{\gamma}{2}-1}\|_{L^2}^2.  \een

Thanks to the fact \beno |1-\chi_R|^2 \le \f{|v|^\gamma}{R^\gamma},\eeno
\eqref{ap1} and \eqref{ap2} will lead to the proposition. \ef



\begin{thebibliography}{99}
\bibitem{al} R. Alexandre, Integral estimates for a linear singular operator linked with the
Boltzmann operator. I. Small singularities $0<\mu<1$, {\it Indiana
Univ. Math. J.} 55 (2006), no. 6, 1975-2021.

\bibitem{al1} R. Alexandre, A Review of Boltzmann Equation with Singular
Kernels, {\it Kinet. Relat. Models} 2 (2009), no. 4, 551-646.


\bibitem{al2} R. Alexandre, Sur le taux de dissipation d'entropie
sans troncature angulaire.  {\it C. R. Acad. Sci. Paris Serie.
Math.} 326, 3 (1998), 311-315.

\bibitem{al3} R. Alexandre, From Boltzmann to Landau, preprint,
1989.

\bibitem{advw} R. Alexandre, L. Desvillettes, C. Villani, and B. Wennberg, Entropy dissipation and long-
range interactions, {\it Arch. Ration. Mech. Anal.} 152 (2000), no.
4, 327-355.

\bibitem{alhe} R. Alexandre and L. He, Integral estimates for a linear singular operator linked with the
Boltzmann operator. II. High singularities $1\le\mu<2$ (2008), no.4,
491--513.

\bibitem{amuxy1} R. Alexandre, Y. Morimoto, S. Ukai, C.-J. Xu, and T. Yang, Uncertainty principle and kinetic
equations, {\it J. Funct. Anal.} 255 (2008), no. 8, 2013-2066.

\bibitem{amuxy2} R. Alexandre, Y. Morimoto, S. Ukai, C.-J. Xu, and T. Yang, Regularizing effect and local existence for
non-cutoff Boltzmann equation, to appear in {\it Arch. Ration. Mech.
Anal.}

\bibitem{amuxy3} R. Alexandre, Y. Morimoto, S. Ukai, C.-J. Xu, and T.
Yang, Global existence and full regularity of the Boltzmann equation
without angular cutoff, preprint 2010.

\bibitem{amuxy4} R. Alexandre, Y. Morimoto, S. Ukai, C.-J. Xu, and T.
Yang, The Boltzmann equation without angular cutoff in the whole
space I: an essential coercivity estimate, preprint 2010.

\bibitem{amuxy5} R. Alexandre, Y. Morimoto, S. Ukai, C.-J. Xu, and T.
Yang, The Boltzmann equation without angular cutoff in the whole
space I: Soft potential, preprint 2010.

\bibitem{alsa} R. Alexandre and M. El Safadi, Littlewood-Paley
theory and regularity issues in Boltzmann homogeneous equations. I.
Non-cutoff case and Maxwellian molecules, {\it Math. Models Methods
Appl. Sci.} 15 (2005), no. 6, 907-920.

 \bibitem {av} R. Alexandre and C. Villani, On the Boltzmann equation for long-range interactions, {\it Comm. Pure Appl.
              Math.}, {\bf 55}, 1 (2002), 30-70.


\bibitem{av2} R. Alexandre and C. Villani,
On the Landau approximation in plasma physics,
 {\it Ann. Inst. H. Poincar\'e, Anal. Non Lin\'eaire}, {\bf 21}, 1 (2004), 61-95.

\bibitem{bob} A. V. Bobylev, The theory of the nonlinear spatially uniform
Boltzmann equation for Maxwell molecules, Mathematical physics
reviews, Vol. 7, Soviet Sci. Rev. Sect. C Math. Phys. Rev., vol. 7,
Harwood Academic Publ., Chur, 1988, pp. 111-233.

\bibitem{bou} F. Bouchut, Hypoelliptic regularity in kinetic equations, {\it J.
Math. Pure Appl.} 81 (2002) 1135-1159.

\bibitem {bd1} F. Bouchut and L. Desvillettes, A proof of the smoothing properties of the positive part of Boltzmann's kernel,
{\it Rev.  Mat. Iberoamericana}, {\bf 14}, 1 (1998), 47-61.

\bibitem {bd3} F. Bouchut and L. Desvillettes, Averaging lemmas without time Fourier transform and application to
              discretized kinetic equations, {\it Proc. Roy. Soc. Edinburgh}, {\bf 129A}, 1 (1999), 19-36.


\bibitem{bd}  L. Bernis and L. Desvillettes, Propagation of
singularities for classical solutions of the
Vlasov-Poisson-Boltzmann equation, {\it Discrete Contin. Dyn. Syst.}
24 (2009), no. 1, 13-33.

\bibitem {bd2} L. Boudin and L. Desvillettes, On the Singularities of the Global Small Solutions of the full Boltzmann
              Equation, {\it Monatschefte f\"{u}r Mathematik}, {\bf 131}, 2 (2000), 91-108.

\bibitem {ce1} Carlo Cercignani, {\em The Boltzmann equation and its applications,
Applied Mathematical Sciences}, vol. 67, Springer-Verlag, New York,
1988.

\bibitem {ce2} Carlo Cercignani, Reinhard Illner, and Mario Pulvirenti, {\em The
mathematical theory of dilute gases}, Applied Mathematical Sciences,
vol. 106, Springer-Verlag, New York, 1994.

\bibitem{CC} S. Chapman and T. G. Cowling, {\em The mathematical theory of non-uniform gases}, Cambridge Univ. Press., London, 1970.


\bibitem {cdh} Y.~Chen, L.~Desvillettes and L.~He. Smoothing effects for classical solutions of the full Landau equation,
{\it Arch. Ration. Mech. Anal.}  193 (2009), no. 1,
21-55.

\bibitem{DTTSP1} L. Desvillettes, On asymptotics of the {B}oltzmann equation when the collisions become grazing,
 {\it Transp. Th. Stat. Phys.}, {\bf 21}, 3 (1992), 259-276.

\bibitem{ld1} L. Desvillettes, About the regularizing properties of the
non-cut-off Kac equation, {\it Comm. Math. Phys.} 168 (1995), no. 2,
417-440.

\bibitem{parme} L. Desvillettes, About the Use of the Fourier Transform for the Boltzmann Equation,
 {\it Riv. Mat. Univ. Parma}, {\bf 7}, 2 (2003), 1-99 (special issue).

\bibitem{dw} L. Desvillettes and B. Wennberg, Smoothness of the
solution of the spatially homogeneous Boltzmann equation without
cutoff, {\it Comm. Partial Differential Equations } 29 (2004), no.
1-2, 133-155.

\bibitem{dg2} L. Desvillettes and  F. Golse,
On a model Boltzmann equation without angular cutoff, {\it Diff. Int.  equation}, {\bf{13}}, 4-6 (2000), 567-594.

\bibitem{dm} Laurent Desvillettes and Clement Mouhot, Stability and uniqueness
for the spatially homogeneous Boltzmann equation with long-range
interactions, {\it Arch. Ration. Mech. Anal.} 193 (2009), no. 2,
227-253.

\bibitem {dl} R. DiPerna and P. L. Lions, On the Cauchy problem for the Boltzmann equation: Global existence and weak
              stability, {\it Ann. Math.}, {\bf 130}, 2 (1989), 312-366.

\bibitem{dly} Renjun Duan, Meng-Rong Li, and Tong Yang, Propagation of
singularities in the solutions to the Boltzmann equation near
equilibrium, {\it Math. Models Methods Appl. Sci.} 18 (2008), no. 7,
1093-1114.

\bibitem {glps} F. Golse, P. L. Lions, B. Perthame and R. Sentis, Regularity of the moments of the solution of a transport equation,
{\it J. Funct. Anal.}, {\bf 76}, 1 (1988), 110-125.


\bibitem{fs} F. Golse and L. Saint-Raymond,  The Navier Stokes limit of the
Boltzmann equation for bounded collision kernels. {\it Invent.
Math.} 155, 81-161 (2004)

\bibitem{goudon} T. Goudon, On the Boltzmann equations and
Fokker-Planck asymptotic: influence of grazing collisions. {\it
J.Sta.Phys.} 89,3-4(1997), 751-776.

\bibitem{gs1} Gressman and Robert M. Strain, Global Strong Solutions of the
Boltzmann Equation without Angular Cut-off, preprint 2009.

\bibitem{gs2} Gressman and Robert M. Strain, Global Classical Solutions of the
Boltzmann Equation with Long-Range Interactions and Soft-Potentials,
preprint 2010.

\bibitem{gs3} Gressman and Robert M. Strain, Sharp anisotropic estimates for the Boltzmann collision operator and
its entropy production, preprint 2010.

\bibitem {guo1} Y. Guo, Classical solutions to the Boltzmann equation for molecules with an angular cutoff,
{\it Arch. Ration. Mech. Anal.} 169 (2003), no. 4, 305-353.

\bibitem {guo2} Y. Guo, The Boltzmann equation in the whole space, {\it Indiana Univ. Math. J.} 53 (2004), no. 4,
1081-1094.

\bibitem{hmuy} Z. Huo, Y. Morimoto, S. Ukai and T. Yang,
Regularity of solutions for spatially homogeneous Boltzmann equation
without angular cutoff. {\it Kinet. Relat. Models } 1 (2008), no. 3,
453-489.

\bibitem{LL} E. M. Lifshitz and L. P. Pitaevskii, {\em {Physical kinetics}}, Perg. Press., Oxford, 1981.


\bibitem{lions} P. Lions, Regularite et compacite pour des noyaux
de collision de Boltzmann sans troncature angulaire, {\it C. R.
Acad. Sci. Paris Serie.  Math.} 326 (1998), no. 1, 37-41.

\bibitem{lm1} P.-L. Lions  and N. Masmoudi,  From the Boltzmann equations to the
equations of incompressible fluid mechanics, I. {\it Archive Rat.
Mech. Anal.}  158 (2001), 173-193.

\bibitem{lm2}
P.-L. Lions  and N. Masmoudi, From the Boltzmann equations to the
equations of incompressible fluid mechanics II. {\it Archive Rat.
Mech. Anal.}
 158 (2001), 195-211.


\bibitem{lyy}  T. Liu, T. Yang, and Shih-Hsien Yu, Energy method for
Boltzmann equation, {\it Phys. D }188 (2004), no. 3-4, 178-192.

\bibitem{lu} X. Lu, A direct method for the regularity of the gain term in the {B}oltzmann equation,
 {\it J. Math. Anal. Appl.}, {\bf 228}, 2 (1998), 409-435.


\bibitem{muxy} Y. Morimoto, S. Ukai, C.-J. Xu, and T. Yang,
Regularity of solutions to the spatially homogeneous Boltzmann
equation without angular cutoff, {\it Discrete Contin. Dyn. Syst. }
24 (2009), no. 1, 187-212.

\bibitem{vi3}  C. Villani, On a new class of weak solutions for the
spatially homogeneous Boltzmann and Landau equations. {\it Arch.
Rat. Mech. Anal.} 143 (1998), 273-307.

\bibitem{vi1} C. Villani, Regularity estimates via the entropy
dissipation for the spatially homogeneous Boltzmann equaiton without
cut-off, {\it Rev.Mat. Iberoam.} 15, 2 (1999), 335-352.

\bibitem{vi2} C. Villani, A review of mathematical topics in collisional kinetic theory,
North-Holland, Amsterdam, Handbook of mathematical fluid dynamics,
Vol. I, 2002, pp. 71-305.


\bibitem {w} B. Wennberg, Regularity in the Boltzmann equation and the Radon transform, {\it Comm. Partial Differential Equations},
              {\bf 19}, 11-12 (1994), 2057-2074.



\end{thebibliography}
\end{document}